\documentclass[12pt,reqno]{amsart}

\usepackage{hyperref}
\usepackage{amsthm}
\usepackage[letterpaper]{geometry}
\usepackage{geometry}
\usepackage{amssymb,mathrsfs}
\usepackage{xcolor}
\usepackage{galois}
\usepackage{enumitem}
\usepackage{mathabx}
\usepackage{makecell}
\usepackage{longtable}
\newtheorem{theorem}{Theorem}[section]
\newtheorem{proposition}[theorem]{Proposition}

\newtheorem{lemma}[theorem]{Lemma}

\newtheorem{corollary}[theorem]{Corollary}

\theoremstyle{remark}

\theoremstyle{definition}
\newtheorem{definition}[theorem]{Definition}
\newtheorem{remark}[theorem]{Remark}

\newtheorem{assumption}[theorem]{Assumption}

 \title[Rank Two Irregular Higgs Bundles]{Asymptotic Geometry of the Moduli Space of Rank Two Irregular Higgs Bundles over the Projective Line}

\author{Gao Chen} \thanks{The first named author is supported by the Project of Stable Support for Youth Team in Basic Research Field, Chinese Academy of Sciences, YSBR-001.}
  \address{Institute of Geometry and Physics, University of Science and Technology of China, Shanghai, China, 201315}
 \email{chengao1@ustc.edu.cn}

 \author{Nianzi Li}
 \address{Department of Mathematics, University of Wisconsin-Madison, Madison, WI 53706, USA}
 \email{nli62@wisc.edu}
 
\begin{document}

 \begin{abstract}
   We study the asymptotic behavior of Hitchin's hyperkähler metric on the moduli space of rank two irregular Higgs bundles over $\mathbb{C}P^1$. Along a generic curve, we prove that the Hitchin metric is asymptotic to the semiflat metric at an arbitrary polynomial order. When there are no weakly parabolic singularities, the rate is exponential. In the case of four-dimensional moduli spaces, we prove that the semiflat metric is asymptotic to an ALG/ALG$^\ast$ model metric.
\end{abstract}

 \maketitle
 \section{Introduction}
 Many moduli spaces from gauge theory admit natural hyperkähler metrics and are noncompact. It is intriguing to study the compactification and the end behavior of the metrics on these moduli spaces. Some of these moduli spaces are expected to be gravitational instantons, i.e., noncompact complete hyperkähler 4-manifold with $\int |\mathrm{Rm}|^2<\infty$. \cite{chen2021gra1} and \cite{sun_zhang_2023collapsing} completely classified gravitational instantons into types ALE, ALF, ALG, ALH, ALG$^\ast$, and ALH$^\ast$. The moduli space of centered charge two $\mathrm{SU}(2)$-monopoles with $k$ singularities is ALF-$D_k$ \cite{CherkisKapustin99, CherkisHitchin05}, and conversely every ALF-$D_k$ gravitational instanton can be realized as such a moduli space \cite{chen2019gra2}. However, little is known about the ALG/ALG$^\ast$ case. Boalch's modularity conjecture \cite{boalch2012hyperkahler, modularity} asserts that every 4d Hitchin moduli space is ALG/ALG$^\ast$, and conversely every ALG/ALG$^\ast$ gravitational instanton can be realized as a Hitchin moduli space. For each type of ALG/ALG$^\ast$ metric, Boalch \cite[Sec.~4.1]{boalch2018wild} proposed a minimal Hitchin moduli space realization. As noncompact analogues of K3 surfaces, they are called H3 surfaces in honor of Higgs, Hitchin and Hodge. Note that some ALG/ALG$^\ast$ metrics can only be realized by allowing irregular singularities (higher order poles). It is proved in \cite{fredrickson2022asymptotic} that the moduli space of Higgs bundles over $\mathbb{C}P^1$ with four strongly parabolic singularities is ALG-$D_4$. We aim to prove similar results for the irregular case.

 Another motivation comes from the GMN (Gaiotto-Moore-Neitzke) conjecture \cite{gaiotto2010four, Neitzke14} on the detailed asymptotic formula of the Hitchin metric, which is a complete hyperkähler metric on the Hitchin moduli space. A consequence of the GMN conjecture is that $g_{L^2}-g_{\mathrm{sf}}=O(\mathrm{e}^{-c t})$ along a generic ray parametrized by $t\in \mathbb{R}^+$, where $g_{L^2}$ is the Hitchin metric, and $g_{\mathrm{sf}}$ is the semiflat metric \cite{freed1999special} arising from the algebraic integrability. With tools from geometric analysis, substantial progress has been made towards proving this result. The pioneering work \cite{mazzeo2019asymptotic} used the gluing construction of harmonic metrics in \cite{mazzeo_swoboda_weiss_witt_2016} to establish the polynomial decay of the above difference. In \cite{dumas2019asymptotics}, Dumas and Neitzke used the local biholomorphic flow to show that the decay is exponential on the Hitchin section. These techniques were combined in \cite{Fredrickson:2018fun, fredrickson_2019} to prove the exponential decay for $\mathrm{SL}(n,\mathbb{C})$-Higgs bundles under a mild assumption on the ramification of the spectral curve if $n\geqslant 3$. Using descriptions of $g_{L^2}$ and $g_{\mathrm{sf}}$ by harmonic 1-forms, Mochizuki \cite{mochizuki_2023} established the exponential decay on the smooth locus of the moduli space of $\mathrm{SL}(n,\mathbb{C})$-Higgs bundles, based on the work \cite{mochizuki_szabo_2023} about large-scale solutions of Hitchin's equations. In \cite{fredrickson2022asymptotic}, the authors considered the moduli space of parabolic $\mathrm{SL}(2,\mathbb{C})$-Higgs bundles, where the Higgs fields admit simple poles. As a further step, it is natural to study the asymptotic behavior of the Hitchin metric for the irregular case.

 In \cite{biquard_boalch_2004}, Biquard and Boalch established the wild non-abelian Hodge correspondence, which implies that the moduli space of irregular Higgs bundles (with untwisted singularities) is a complete hyperkähler manifold. In this paper we consider the moduli space $\mathcal{M}$ of rank two irregular Higgs bundles over $\mathbb{C}P^1$. Using the gluing construction similar to \cite{mazzeo_swoboda_weiss_witt_2016, fredrickson2022asymptotic}, we can control harmonic metrics on Higgs bundles near the ends of $\mathcal{M}$. From this, under the extra assumption that the parabolic structure at every weakly parabolic point is trivial, we obtain the following main result on the metric comparison.

 \begin{theorem}\label{Main_thm}
   Fix a generic curve $[(\bar{\partial}_E,\varphi_t)]$ in $\mathcal{M}$, and an infinitesimal deformation $[(\dot{\eta},\dot{\varphi})]\in T_{[(\bar{\partial}_E,\varphi_t)]}\mathcal{M}$. As $t\to \infty$, for any $N>0$ we have
 \[\lVert [(\dot{\eta},\dot{\varphi})]\rVert_{g_{L^2}}^2-\lVert [(\dot{\eta},\dot{\varphi})]\rVert_{g_{\mathrm{sf}}}^2=O(t^{-N}).\]
 If, moreover there are no weakly parabolic points, there exist positive constants $c,\sigma$ independent of $t$ and $[(\dot{\eta},\dot{\varphi})]$ such that the above difference is $O(\mathrm{e}^{-ct^{\sigma}})$.
 \end{theorem}
 Here $O(t^{-N})$ means that it is bounded by $c_Nt^{-N}\lVert [(\dot{\eta},\dot{\varphi})]\rVert_{g_{\mathrm{sf}}}^2$ for some $c_N>0$ independent of $t$ and $[(\dot{\eta},\dot{\varphi})]$. The precise definition of the curve $[(\bar{\partial}_E,\varphi_t)]$ will be given in Section \ref{approxsol_sec}.

 Later in Section 6.3, we will specialize to the case when the moduli space is four dimensional. With this result, we are able to show that $g_{L^2}$ is polynomially close to an ALG/ALG$^\ast$ model metric $g_{\mathrm{model}}$. Moreover, the constants in the difference $g_{L^2}-g_{\mathrm{model}}$ can be shown to be independent of the choice of the curve $[(\bar{\partial}_E,\varphi_t)]$, i.e., uniform near the ends of $\mathcal{M}$. The correspondence between the configuration of singularities and the type of the model metric is listed below.
 \renewcommand{\arraystretch}{1.3}
 \begin{table}[h]
   \centering
   \begin{tabular}{|c|c|c|c|}
     \hline
  Kodaira type&$\uppercase\expandafter{\romannumeral2}^\ast$ & $\uppercase\expandafter{\romannumeral3}^\ast$ & $\uppercase\expandafter{\romannumeral4}^\ast$\\
 Dynkin diagram &$A_0$&$A_1$&$A_2$\\
  $\beta$&$\frac{5}{6}$&$\frac{3}{4}$&$\frac{2}{3}$\\
 $\tau$&$\mathrm{e}^{2\pi \mathrm{i}/3}$&$\mathrm{i}$&$\mathrm{e}^{2\pi \mathrm{i}/3}$\\
  $D$&$4\cdot\{\tilde{0}\}$&$4\cdot\{0\}$ or $3\cdot\{\tilde{0}\}+\{\infty\}$&$3\cdot\{0\}+\{\infty\}$\\
  $\mathrm{dim}_{\mathcal{M}_{\mathrm{ALG}}}$&$2$&$5$&$8$\\
  parameters &$2$&$5$&$8$\\
  \hline
   \end{tabular}
   \vspace{0.5em}
   \caption{ALG}
   \label{ALG_tab}
 \end{table}

 The tilde $(\tilde{~})$ on the pole indicates a twisted irregular type (defined in Section 2). To the best of our knowledge, the analytic construction of $\mathcal{M}$ with twisted singularities as a complete hyperkähler manifold has not been rigorously established in the existing literature, so our discussion of the hyperkähler metric on such $\mathcal{M}$ should be regarded as formal. Nevertheless, in the rank two case, there should be no essential technical obstacles in extending \cite{biquard_boalch_2004} to include twisted poles, since locally by passing to a double cover, a twisted pole becomes untwisted. This is not true in higher rank, and we refer to \cite{Szabo16} for local normal forms of the Higgs field near twisted poles.

 In the previous table, $(\beta, \tau)$ means that $\mathcal{M}$ is asymptotic to the standard ALG model $(X,G_{\beta,\tau})$ of type $(\beta,\tau)$. Here $X$ is obtained by identifying two boundary components of \[\{u\in \mathbb{C}\,|\,\mathrm{Arg}(u)\in [0,2\pi\beta]\text{ and }|u|\geq R\}\times \mathbb{C}_v/(\mathbb{Z}\oplus \mathbb{Z}\tau)\] via the gluing map $(u,v)\sim (\mathrm{e}^{2\pi \mathrm{i}\beta}u,\mathrm{e}^{-2\pi \mathrm{i}\beta}v)$. $G_{\beta,\tau}$ is a flat hyperkähler metric on $X$ such that $\omega^1=\mathrm{i}/2 (\mathrm{d}u\wedge \mathrm{d}\bar{u}+\mathrm{d}v\wedge \mathrm{d}\bar{v})$ and $\omega^2+\mathrm{i}\omega^3=\mathrm{d}u\wedge \mathrm{d}v$. Kodaira type means that $\mathcal{M}$ is biholomorphic to a rational elliptic surface minus a fiber with the given Kodaira type \cite{chen_chen_2020}. Dynkin diagram means that $H^2(\mathcal{M})$ is generated by the given extended Dynkin diagram. This makes sense because ALG gravitational instantons with the same $\beta$ are diffeomorphic to each other \cite{chen2021gravitational}. $\mathrm{dim}\mathcal{M}_{\mathrm{ALG}}$ is the dimension of the deformation space of ALG gravitational instantons of that type. The parameters consist of singularity data and parabolic weights which determine the Hitchin moduli space $\mathcal{M}$. The counting of parameters (will be explained in Section \ref{ALG_subsubsec}) gives evidence of the modularity conjecture. The moduli space of Higgs bundles with 4 parabolic singularities (toy models) has 12 parameters, we refer to \cite[Sec.~6.3]{swoboda21} for a more detailed discussion of this case.

 \begin{table}[h]
   \centering
   \begin{tabular}{|c|c|c|c|c|}
     \hline
  Kodaira type&$\uppercase\expandafter{\romannumeral1}_4^\ast$&$\uppercase\expandafter{\romannumeral1}_3^\ast$&$\uppercase\expandafter{\romannumeral1}_2^\ast$&$\uppercase\expandafter{\romannumeral1}_1^\ast$\\
 Dynkin diagram &$D_0$&$D_1$&$D_2$&$D_3=A_3$\\
  $D$&$2\cdot\{\tilde{0}\}+2\cdot\{\tilde{\infty}\}$&$2\cdot\{0\}+2\cdot\{\tilde{\infty}\}$&\makecell{$2\cdot\{0\}+2\cdot\{\infty\}$ \\or $\{0\}+\{1\}+2\cdot\{\tilde{\infty}\}$}&$\makecell{\{0\}+\{1\}\\+2\cdot\{\infty\}}$\\
  $\mathrm{dim}_{\mathcal{M}_{\mathrm{ALG}^\ast}}$&$0$&$3$&$6$&$9$\\
  parameters&$2$&$5$&$8$&$11$\\
  \hline
   \end{tabular}
   \vspace{0.5em}
   \caption{ALG$^\ast$}
 \end{table}
 \vspace{-1em}

 In the previous table,
 Kodaira type means that $\mathcal{M}$ is biholomorphic to a rational elliptic surface minus a fiber with the given Kodaira type \cite{chen2021gravitational}. Dynkin diagram means that $H^2(\mathcal{M})$ is generated by the given extended Dynkin diagram. This makes sense because ALG$^*$ gravitational instantons with the same Kodaira type at infinity are diffeomorphic to each other \cite{chen2021gravitational}. The ALG$^*$ model metric is a $\mathbb{Z}_2$ quotient of the Gibbons-Hawking metric \[V(\mathrm{d}x^2+\mathrm{d}y^2+\mathrm{d}\theta_2^2)+V^{-1}\Theta^2,\]
 where $\Theta=(\nu/\pi)(\mathrm{d}\theta_3-\theta_2 \mathrm{d}\theta_1), V=\kappa_0+(\nu/\pi)\log\,r$. See \cite[Sec.~2]{chen2021gravitational} for more details. The ALG$^\ast$ cases have two redundant parameters, indicating that there are some extra symmetries among these moduli spaces.

 To prove Theorem \ref{Main_thm}, we use a similar strategy as that in \cite{fredrickson2022asymptotic}. Their strategy is:
 \begin{itemize}
 \item	Construct the harmonic metric $h_t$ for $(\bar{\partial}_E,t\varphi)$ by a gluing method, with the following building blocks.
 \begin{itemize}
   \item Limiting metric $h_\infty$, which is singular at the points where the spectral cover is ramified.
   \item Fiducial solutions $h_{t}^{\mathrm{model}}$ in some fixed disks around the ramification points.
   \item Desingularize $h_\infty$ by $h_{t}^{\mathrm{model}}$ to obtain $h_t^{\mathrm{app}}$, and perturb $h_t^{\mathrm{app}}$ to $h_t$.
 \end{itemize}
 \item Interpolate $g_{\mathrm{app}}$ between $g_{L^2}$ and $g_{\mathrm{sf}}$, compare tangent vectors in Coulomb gauges using $L^2$ descriptions of these metrics.
 \end{itemize}
 However, there are several new features and difficulties in our case. The first is the lack of a natural $\mathbb{C}^\ast$ action on $\mathcal{M}$ since the singularity types are fixed. Fortunately, irregular Higgs bundles over $\mathbb{C}P^1$ can be explicitly described, and we can find some curve $[(\bar{\partial}_E,\varphi_t)]$ in $\mathcal{M}$ tending to infinity in place of the ray. Another problems is that if we just use the fiducial solutions as in \cite{fredrickson2022asymptotic}, the approximate harmonic metric we construct will not make the error arbitrarily small as $t\to\infty$. Instead, we will use another model solution near weakly parabolic points, which is not an explicit ODE solution. This model metric is the harmonic metric on the Higgs bundle over $\mathbb{C}P^1$ with a weakly parabolic pole and a twisted order two pole, which exists by \cite{harland22parabolic, mochizuki21good}. Near the boundary a disk we show that the harmonic metric is asymptotic to the decoupled harmonic metric at an arbitrary polynomial rate in $t$. Recently, we are informed by Mochizuki that this can be possibly improved to an exponential rate using his recent work \cite{mochizuki_2023} with an assumption on parabolic weights. Next, the analysis of the linearized operator requires modifications because of the presence of irregular singularities. We use the analytic tools from \cite{biquard_boalch_2004} to handle these singularities. Finally, to compare $g_{\mathrm{app}}$ and $g_{\mathrm{sf}}$, we should modify the proof of \cite{fredrickson2022asymptotic} for weakly parabolic points. With the aforementioned assumption that the two parabolic weights are equal at each weakly parabolic point, we can apply the idea of \cite{dumas2019asymptotics}, but we will lose a parameter (weight) for each weakly parabolic point.

 This paper is organized as follows. In Section 2, we give some definitions relating to irregular Higgs bundles, and describe the Hitchin base and generic Hitchin fibers explicitly. In Section 3, we introduce the building blocks, which are the model solutions and the decoupled harmonic metric, and we construct an approximate solution by gluing these together. Note that in our case the zeros of $\mathrm{det}\,\varphi_t$ are moving with $t$, so the gluing regions also depend on $t$. The decaying rates of the error terms in different regions depend on some local parameters, which are called \emph{local masses} by analogy with those appearing in \cite{foscolo2017gluing}. We analyze the linearized Hitchin equation in Section 4, and nonlinear terms in Section 5, then a genuine solution to the Hitchin equation is obtained. Finally in Section 6, we prove Theorem \ref{Main_thm}, and analyze the asymptotics of semiflat metrics for the four-dimensional moduli spaces.
 \paragraph{\bfseries Acknowledgements} The authors benefited from discussions with Laura Fredrickson on twisted singularities, and with Bin Xu on quadratic differentials. We are grateful to Takuro Mochizuki for explaining his recent work and for many insightful discussions. We would like to thank Mao Sheng, Song Sun, Alex Waldron, and Ruobing Zhang for discussions and their interests in this work.

 \section{Preliminaries}\label{Prelim_sec}
 \subsection{Irregular Higgs bundles and the moduli space}\label{irreg_subsec}
 Let $E$ be a complex vector bundle of rank 2 over $C=\mathbb{C}P^1$, $S$ be a set of points on $C$. Cover $C$ by the usual coordinate charts $U=\mathbb{C}_z$, $V=\mathbb{C}_w$ with $w=1/z$, and such that $\{w=0\}\notin S$. For $x\in S$, we also denote its $z$-coordinate by $x\in \mathbb{C}$. Fix the \emph{parabolic structure} at each $x\in S$, which consists of a one-dimensional subspace $L_x\subset E_x$ (a full flag) and parabolic weights $1/2<\alpha_{x,2}<1$ (resp. $\alpha_{x,1}=1-\alpha_{x,2}$) asscociated to $L_x$ (resp. $E_x$). When $\alpha_{x,1}=\alpha_{x,2}$, we get a trivial parabolic structure $\{0\}\subset E_x$ with weight $1/2$ asscociated to $E_x$. We allow this case for weakly parabolic singularities.
 \begin{definition}[\cite{biquard_boalch_2004}]
   An $\mathrm{SL}(2,\mathbb{C})$-irregular Higgs bundle  is a pair $(\bar{\partial}_E,\varphi)$, where $\bar{\partial}_E$ is a holomorphic structure on $E$, and the Higgs field $\varphi\in H^0(C,\mathfrak{sl}(\mathcal{E})\otimes K(D))$. Here $\mathcal{E}:=(E,\bar{\partial}_E)$ and $D$ is a divisor $D=\sum_{x\in S} m_x\cdot \{x\}$ with $m_x\in \mathbb{Z}^+$.
 \end{definition}

  Write $S=I\cup P$, where $m_x>1$ for $x\in I$ and $m_x=1$ for $x\in P$. Then $\varphi$ is a meromorphic section of $\Omega^{1,0}(\mathfrak{sl}(E))$ with singularities at points in $S$, which are called \emph{irregular} (resp. \emph{parabolic}) singularities when they belong to $I$ (resp. $P$). We require that $\varphi$ has a fixed \emph{singularity type} at each $x\in S$, meaning that there exists a holomorphic trivialization of $(E,\bar{\partial}_E)$ over a neighborhood $U_x$ of $x$ such that
 \[\varphi=\left(\frac{\phi_{x,m_x}}{z_x^{m_x}}+\cdots+\frac{\phi_{x,1}}{z_x}+\text{holomorphic terms}\right)\,\mathrm{d}z_x,\]
 where $z_x=z-x$, and $\phi_{x,m_x},\ldots,\phi_{x,1}\in \mathfrak{sl}(2,\mathbb{C})$ are given. We decompose $I=I_u\cup I_t$, and $P=P_w\cup P_s$. For $x\in I_u$ (resp. $P_w$), the leading term $\phi_{x,m_x}$ is diagonalizable with opposite nonzero eigenvalues $\pm \rho_{x,m_x}$, and the singularity $x$ is called \emph{untwisted} (resp. \emph{weakly parabolic}). For $x\in I_t$ (resp. $P_s$), $\phi_{x,m_x}$ is nilpotent, and $x$ is called \emph{twisted} (resp. \emph{strongly parabolic}). Note that the choice of $\phi_{x,1}$ for $x\in P_s$ does not impose any restriction, since there always exists a local holomorphic frame in which the residue of $\varphi$ is $\phi_{x,1}$. We require that $\varphi$ is adapted to the parabolic structure, i.e., $\phi_{x,m_x}(L_x)=L_x$ for $x\in I_u\cup P_w$, and $\phi_{x,m_x}(L_x)=\{0\}$ for $x\in I_t\cup P_s$.
 \begin{definition}
 An irregular Higgs bundle $(\bar{\partial}_E,\varphi)$ is \emph{stable} if for every $\varphi$-invariant holomorphic line subbundle $\mathcal{F}$ of $\mathcal{E}$,
 \[\mathrm{pdeg}\,\mathcal{F}<\frac{1}{2}\mathrm{pdeg}\,\mathcal{E}, \]
 where $\mathrm{pdeg}\,\mathcal{E}=\mathrm{deg}\,\mathcal{E}+|S|$, and $\mathrm{pdeg}\,\mathcal{F}=\mathrm{deg}\,\mathcal{F}+\sum_{x\in S}\tilde{\alpha}_x$, $\tilde{\alpha}_x=\alpha_{x,2}$ if $\mathcal{F}_x=L_x$, otherwise $\tilde{\alpha}_x=\alpha_{x,1}$.
 \end{definition}

 \begin{definition}
   The \emph{moduli space} of irregular Higgs bundles is
 \begin{align*}
   \mathcal{M}=\{(\bar{\partial}_E,\varphi) \text{ stable and compatible with the parabolic }&\\\text{structure and the singularity type at each }x\in S\}&\,/\,\mathscr{G}_{\mathbb{C}},
 \end{align*}
 where $\mathscr{G}_{\mathbb{C}}=\Gamma(\mathrm{ParEnd}(E)\cap \mathrm{SL}(E))$, which consists of sections of $\mathrm{SL}(E)$ preserving $L_x$ for each $x\in S$. Alternatively,
  $\mathcal{M}=\{(\bar{\partial}_E,\varphi,L)\}\,/\,\Gamma(\mathrm{SL}(E))$, where $L=\{L_x\}_{x\in S}$, $L_x\in \mathbb{CP}^1_x$ which parametrizes one-dimensional subspaces of $E_x$.
   The gauge group acts as
 \[g\cdot (\bar{\partial}_E,\varphi,L)=(g^{-1}\comp \bar{\partial}_E\comp g,g^{-1}\varphi g,g^{-1}\cdot L).\]
 By \cite[Th.~3.1]{inaba_michi_2016}, $\mathrm{dim}_{\mathbb{C}}\,\mathcal{M}=2(N-3)$, where $N=\mathrm{deg}\,D$. We assume that $N\geqslant 4$.
 \end{definition}

 In the triple $(\bar{\partial}_E,\varphi,L)$, as $\varphi$ is compatible with $L$, we can determine $L_x$ from $\varphi$ as the kernel of $\phi_{x,m_x}$ for $x\in I_t\cup P_s$, and as one of the eigenspace of $\phi_{x,m_x}$ for $x\in I_u\cup P_w$. There are $2^{|I_u|+|P_w|}$ choices, depending on whether $L_x$ corresponds to $\rho_{x,m_x}$ or $-\rho_{x,m_x}$, where $0\leqslant \mathrm{Arg}(\rho_{x,m_x})<\pi$. These $2^{|I_u|+|P_w|}$ components of $\mathcal{M}$ are equivalent, so we only need to consider the one where $L_x$ corresponds to $\rho_{x,m_x}$, and still denote this component by $\mathcal{M}$. As $L$ is determined by $\varphi$, we will write $[(\bar{\partial}_E,\varphi,L)]\in\mathcal{M}$ as $[(\bar{\partial}_E,\varphi)]$.

 \begin{definition}
   The \emph{Hitchin fibration} is defined as \[H:\mathcal{M}\to\mathcal{B},~[(\bar{\partial}_E,\varphi)]\mapsto \det\varphi.\]
   Here $\mathcal{B}\subset H^0(C,K(D)^2)$ consists of quadratic differentials $\nu$ which has the expansion
   \[\nu=\left(\frac{\mu_{x,2m_x}}{z_x^{2m_x}}+\cdots+\frac{\mu_{x,m_x+1}}{z_x^{m_x+1}}+\text{higher order terms}\right)\,\mathrm{d}z_x^2,\]
   near each $x\in I\cup P_w$. $\mu_{x,j}$'s are complex numbers fixed by the singularity type at $x$, and clearly $\mu_{x,2m_{x}}\neq 0$ for $x\in I_u\cup P_w$, while $\mu_{x,2m_{x}}= 0$ for $x\in I_t$. Moreover, we assume that $\mu_{x,2m_{x}-1}\neq 0$ for $x\in I_t$ (when it is zero, by choosing a different extension of the holomorphic structure over $x$, this case can be reduced to the nonzero case or the case where the leading coefficient is diagonalizable \cite[p.~53]{witten2008gauge}).
 \end{definition}
 \subsection{Explicit descriptions}
 \begin{lemma}
   The \emph{Hitchin base} $\mathcal{B}$ is a dimension $N-3$ affine subspace of $H^0(C,K(D)^2)$, which can be explicitly written as
   \begin{align}
     \mathcal{B}=\left\{\,\left(\sum_{x\in I_{\geqslant 3}}\sum_{a=m_x+1}^{2m_x}\frac{\mu_{x,a}}{(z-x)^a}+\sum_{x\in I_{=2}}\left( \frac{\mu_{x,4}}{(z-x)^4}+\frac{\mu_{x,3}(x-y_x)}{(z-x)^3(z-y_x)}\right)\hspace{2.5cm}\right.\right.\notag\\+\left.\left.\left.\sum_{x\in P_w}\frac{\mu_{x,2}(x-y_x')(x-y_x'')}{(z-x)^2(z-y_x')(z-y_x'')}+\frac{\sum_{b=0}^{N-4}\nu_{b}z^b}{\prod_{x\in S}(z-x)^{m_x}}\right)\,\mathrm{d}z^2\,\right|\, \nu_0,\ldots,\nu_{N-4}\in \mathbb{C} \right\},\label{HitBase_eq}
   \end{align}
     where $I_{\geqslant 3}$ consists of poles with order $\geqslant 3$, and $I_{=2}$ consists of order two poles. If $I_{=2}\neq \varnothing$, then we can choose $y_x\in S\backslash\{x\}\neq \varnothing$ since $N\geqslant 4$. Similarly, we choose $y_x',y_x''\in S\backslash\{x\}$ for $x\in P_w$, and $y_x'=y_x''$ only when $y_x'\in I$. Choosing a different $y_x$ (or $y_x'$, $y_x''$) for $x\in I_{=2}$ (or $P_w$) amounts to shifting $(\nu_0,\ldots,\nu_{N-4})$ by an element of $\mathbb{C}^{N-3}$.
 \end{lemma}
 \begin{proof}
 By considering the difference of two elements in $\mathcal{B}$, one can see that $\mathcal{B}$ is an affine space modeled on $H^0(C,K^2\otimes \mathcal{O}(D))$, which has dimension $N-3$. Let the right hand side of \eqref{HitBase_eq} be $\mathcal{B}_1$, then $\mathcal{B}_1\subset \mathcal{B}$. Since $\prod_{x\in S}(z-x)^{-m_x}z^i$ ($i=0,\ldots,N-4$) are linearly independent, we have $\mathrm{dim}\,\mathcal{B}_1=N-3$ and $\mathcal{B}_1=\mathcal{B}$. For $x\in I_{=2}$, and $y_{x,1},y_{x,2}\in S\backslash \{x\}$, \[\frac{\mu_{x,3}(x-y_{x,1})}{(z-x)^3(z-y_{x,1})}-\frac{\mu_{x,3}(x-y_{x,2})}{(z-x)^3(z-y_{x,2})}=\frac{\mu_{x,3}(y_{x,2}-y_{x,1})}{(z-x)^2(z-y_{x,1})(z-y_{x,2})},\]
 the last statement follows. The case for $x\in P_w$ is similar.
 \end{proof}
 \begin{remark}
   Let $\widetilde{\mathcal{M}}$ be a moduli space of $\mathrm{GL}(2,\mathbb{C})$-irregular Higgs bundles, defined in the same way as for $\mathcal{M}$ except that $\varphi$ belongs to $H^0(C,\mathfrak{gl}(E)\otimes K(D))$. Then for two different Higgs fields $\varphi_1,\varphi_2$ in $\widetilde{\mathcal{M}}$, since the singularity types are fixed, we have $\mathrm{tr}(\varphi_1)-\mathrm{tr}(\varphi_2)\in H^0(C,K)$ which is trivial. This means that Higgs fields in $\widetilde{\mathcal{M}}$ have constant traces, and there is no loss in considering $\mathrm{SL}(2,\mathbb{C})$-irregular Higgs bundles over $C$.
 \end{remark}
 Let $\pi:K(D)\to C$ be the projection, and $\lambda$ be the tautological section of $\pi^\ast K(D)$. Each point $\nu$ in $\mathcal{B}$ defines a section of $(\pi^\ast K(D))^2$ of the form $\lambda^2+\pi^\ast \nu$ whose zero locus $S_\nu$ is called the \emph{spectral curve} determined by $\nu$. From now on we assume that $S_\nu$ is smooth, which is true for generic $\nu$ by Bertini's theorem. $\pi: S_\nu\to C$ is a two-fold ramified cover with ramification divisor $R_\nu$. By the Riemann-Hurwitz formula, the genus of $S_\nu$ is $N-3$. For $\varphi$ with $\mathrm{det}\,\varphi=\nu$, one can associate the \emph{spectral line bundle} $\mathcal{L}_\varphi$ over $S_\nu$ such that $\mathcal{L}_\varphi(-R_\nu)\subset \pi^{\ast} E$ and $\pi_\ast \mathcal{L}_\varphi=E$. Away from the support of $R_\nu$, the fiber $\mathcal{L}_{\varphi,q}$ is the eigenspace in $E_{\pi(q)}$ of $\varphi$ with the eigenvalue $q$. By the Grothendieck–Riemann–Roch theorem \cite[p.~96]{logares_martens_2010},
 \begin{equation}\label{SpecLinedeg_eq}
   \mathrm{deg}\,\mathcal{L}_\varphi=\mathrm{deg}\,E+N-2.
 \end{equation}

 Assume for simplicity that $\mathrm{pdeg}\,E=0$ (if $\mathrm{pdeg}\,E=2d$ for some $d\in \mathbb{Z}$, one can replace $E$ by $E\otimes \mathcal{O}(-d)$ \cite[p.~624]{fredrickson2021moduli}). Then $\mathrm{deg}\,E=-|S|$ and $\mathcal{E}=(E,\bar{\partial}_E)\cong \mathcal{O}(m)\oplus \mathcal{O}(-|S|-m)$, $m\geqslant -|S|/2$, by the Birkhoff–Grothendieck theorem. Now \[\mathrm{End}\,\mathcal{E}\cong \begin{pmatrix}
   \mathcal{O}&\mathcal{O}(2m+|S|)\\ \mathcal{O}(-|S|-2m)&\mathcal{O}
 \end{pmatrix}.\]
 In a trivialization of $\mathrm{End}\,\mathcal{E}$ over $U$, we have
 \begin{equation}\label{GlobalHiggs_eq}
   \varphi=\frac{\mathrm{d}z}{\prod_{x\in S}(z-x)^{m_x}}\begin{pmatrix}
     a(z)&b(z)\\ c(z)&-a(z)
   \end{pmatrix},\text{ where }a(z),b(z),c(z)\in \mathbb{C}[z].
 \end{equation}
 The following lemma shows that the singularity data are essentially equivalent to the data of the Hitchin base.
 \begin{lemma}
 Every $\varphi$ in \eqref{GlobalHiggs_eq} with $\det\varphi\in\mathcal{B}$ of \eqref{HitBase_eq} is compatible with some fixed singularity data at points in $I\cup P_w$, i.e.,
 \begin{enumerate}[label=(\roman*)]
   \item there exists a local holomorphic frame around $x\in I_u\cup P_w$ in which
     \begin{equation}\label{DiagLocHiggs_eq}
       \varphi=\sum_{j=1}^{m_{x}}\rho_{x,j}z_{x}^{-j}\sigma_3\,\mathrm{d}z_{x}+\text{holomorphic terms},
     \end{equation}
   where $\sigma_3$ is the Pauli matrix $\mathrm{diag}(1,-1)$ and $\rho_{x,j}$'s are determined by $\mu_{x,a}$'s in \eqref{HitBase_eq};
   \item there exists a local holomorphic frame around $x\in I_t$ in which
     \begin{equation}\label{NilLocHiggs_eq}
       \varphi=\begin{pmatrix}
         0&-\sum_{j=1}^{m_{x}-1}\mu_{x,m_{x}+j}z_{x}^{-j}\\ z_{x}^{-m_{x}}&0
       \end{pmatrix}\,\mathrm{d}z_{x}+\text{holomorphic terms}.
     \end{equation}
 \end{enumerate}
 \end{lemma}
 \begin{proof}
   (\romannumeral1) This is essentially \cite[Lem.~1.1]{biquard_boalch_2004}. By \eqref{HitBase_eq} and \eqref{GlobalHiggs_eq}, in a neighborhood of $x$, we have the expansion
   \begin{align}
     \det\varphi(z_{x})&=\bigg(\sum_{j=m_{x}+1}^{2m_{x}}\mu_{x,j}z_{x}^{-j}+\text{higher order terms}\bigg)\,\mathrm{d}z_{x}^2,\label{DiagLocDetExp_eq}\\
     \varphi(z_{x})&=\bigg(\sum_{j=1}^{m_{x}}\phi_{x,j}z^{-j}+\text{holomorphic terms}\bigg)\,\mathrm{d}z_{x},\label{DiagLocHiggsExp_eq}
   \end{align}
  $\mu_{x,2m_{x}}=\det\phi_{x,m_{x}}\neq 0$, so we can find a frame where $\phi_{x,m_{x}}=\rho_{x,m_{x}}\sigma_3$ with $\rho_{x,m_{x}}=\mu_{x,2m_{x}}^{1/2}$. Then we use gauge transformations of the form $g=1-g_jz_{x}^j$ ($j=1,\ldots,m_{x}-1$) to successively cancel the off-diagonal terms of $\phi_{x,m_{x}-j}$. In each stage, if necessary, we shrink the neighborhood of $0$ so that $g$ is invertible. For example, for $g=1-g_1z$, \[g^{-1}\varphi g=\left(\phi_{x,m_{x}}z^{-m_{x}}+(\phi_{x,m_{x}-1}-[\phi_{x,m_{x}},g_1])z^{-m_{x}+1}+O\left(z^{-m_{x}+2}\right)\right)\,\mathrm{d}z_{x},\]
   and $\phi_{x,m_{x}-1}-[\phi_{x,m_{x}},g_1]$ becomes diagonal if we choose
   \[g_1=\begin{pmatrix}
   0&\phi_{x,m_{x}-1,12}/2\\-\phi_{x,m_{x}-1,21}/2&0
 \end{pmatrix},\text{ where }\phi_{x,m_{x}-1,ij}\text{ is the }(i,j)-\text{entry of }\phi_{x,m_{x}-1}.\]
 Therefore, in some local holomorphic frame around $x$, $\varphi$ has the form in \eqref{DiagLocHiggs_eq}. Now $\rho_{x,m_{x}}=\mu_{x,2m_{x}}^{1/2}\neq 0$, and $\rho_{x,j}$ ($j=m_{x}-1,\ldots,1$) can be recursively determined as
 \[\rho_{x,j}=(-2\rho_{x,m_{x}})^{-1}\bigg(\mu_{x,j+m_{x}}+\sum_{k=1}^{m_{x}-j-1}\rho_{x,j+k}\rho_{x,m_{x}-k}\bigg).\]
 (\romannumeral2) The proof is similar to (\romannumeral1). Around $x$, we can expand $\det\varphi$ and $\varphi$ as in \eqref{DiagLocDetExp_eq}, \eqref{DiagLocHiggsExp_eq}. Since $\mu_{x,2m_{x}}=0$ and $\mu_{x,2m_{x}-1}\neq 0$, then $\phi_{x,m_{x}}\neq 0$ and is nilpotent, so we can find a local frame in which $\phi_{x,m_{x}}=\left(\begin{smallmatrix}
     0&0\\1&0
   \end{smallmatrix}\right)$. Then we use gauge transformations of the form $g=1-g_jz_{x}^j$ ($j=1,\ldots,m_{x}-1$) to successively cancel the lower-triangular part of $\phi_{x,m_{x}-j}$. For example if we choose \[g_1=\begin{pmatrix}
     \phi_{x,m_{x}-1,21}&\phi_{x,m_{x}-1,22}\\0&0
   \end{pmatrix}\]
   in $g=1-g_1z$, then $\phi_{x,m_{x}-1}-[\phi_{x,m_{x}},g_1]$ becomes strictly upper-triangular. Therefore, in some local holomorphic frame around $x$, $\varphi$ has the form in \eqref{NilLocHiggs_eq}.
 \end{proof}
 In $V$, we can write $\varphi$ as
 \[\varphi=\frac{-w^{N-2}\,\mathrm{d}w}{\prod_{x\in S}(1-xw)^{m_x}}\begin{pmatrix}
   a(w^{-1})&w^{2m+|S|}b(w^{-1})\\ w^{-2m-|S|}c(w^{-1})&-a(w^{-1})
 \end{pmatrix}.\]
 $\varphi$ is regular at $\infty$, so $\mathrm{deg}\,a\leqslant N-2$, $\mathrm{deg}\,b\leqslant 2m+|S|+N-2$, $\mathrm{deg}\,c\leqslant -2m-|S|+N-2$. Since $ \mathrm{deg}\,c\geqslant 0$, and recall that $m\geqslant -|S|/2$, we have
 \[-|S|\leqslant 2m\leqslant -|S|+N-2,\]
 which yields a stratification of $\mathcal{M}$:
 \begin{equation}\label{Strata_eq}
   \mathcal{M}=\bigsqcup_{m=\lceil -|S|/2\rceil}^{\lfloor (-|S|+N-2)/2 \rfloor}\mathcal{M}_m,
 \end{equation}
 where $[(\bar{\partial}_E,\varphi)]\in \mathcal{M}_m$ if $(E,\bar{\partial}_E)\cong \mathcal{O}(m)\oplus \mathcal{O}(-|S|-m)$. Let $H_m$ be the restriction of $H$ to $\mathcal{M}_m$. Note that
 \begin{equation}\label{detHiggs_eq}
   -a(z)^2-b(z)c(z)=\prod_{x\in S}(z-x)^{2m_x}\nu(z):=\tilde{\nu}(z) \text{ where } \det\varphi=\nu=\nu(z)\,\mathrm{d}z^2,
 \end{equation}
  then $b(z)=(-a(z)^2-\tilde{\nu}(z))/c(z)$. The spectral cover $S_\nu\to C$ ramifies exactly at the $2N-4$ zeros (counted with multiplicity) of $\tilde{\nu}(z)$. We say that $\nu\in\mathcal{B}'$, the \emph{regular locus}, if $\tilde{\nu}(z)$ only has simple zeros. When $\nu=\det\varphi\in \mathcal{B}'$, the Higgs bundle is stable \cite[p.~7]{mazzeo_swoboda_weiss_witt_2016} and $S_\nu$ is smooth \cite[p.~10]{fredrickson2022asymptotic}. Fix $\nu\in\mathcal{B}'$, and next we study the fiber $H_m^{-1}(\nu)$. Parabolic Higgs bundles have been studied in \cite{fredrickson2022asymptotic}, so from now on we suppose that $I\neq\varnothing$. If $I_u\neq \varnothing$, by a coordinate change when necessary, we may assume that $0\in I_u$, or else we assume that $0\in I_t$.

 If $m=-|S|/2$, then $|S|$ is even, the holomorphic gauge transformations are of the form
 \[g=\begin{pmatrix}
   g_{11}&g_{12}\\g_{21}&g_{22}
 \end{pmatrix},\quad g_{11},g_{12},g_{21},g_{22}\in \mathbb{C},\det g=1.\]
 Let $a_i,b_i,c_i,\tilde{\nu}_i$ be the coefficients of $z^i$ in $a(z),b(z),c(z),\tilde{\nu}(z)$, and make the generic assumption that $\tilde{\nu}_{2N-4},\tilde{\nu}_{2N-5}\neq 0$. By some gauge transformation $g$ above, we can make $b_0=c_0=0$ when $0\in I_u$, or $a_0=c_0=0$ when $0\in I_t$. Now $c(z)\neq 0$, otherwise $\tilde{v}(z)=-a(z)^2$, contradicting that $\tilde{\nu}(z)$ only has simple zeros. If $\mathrm{deg}\,c=N-2$, using the diagonal gauge transformation $g=\mathrm{diag}(c_{N-2}^{-1/2},c_{N-2}^{1/2})$ we can make $c(z)$ monic. Then we apply \[g=\begin{pmatrix}
   1& a_{N-2}\\ 0&1
 \end{pmatrix}\]
 to make $a_{N-2}=0$. There is no gauge freedom left. The remaining $N-2$ (if $0\in I_u$, or $N-3$ if $0\in I_t$) coefficients of $a(z)$ can be determined (up to finitely many choices) by $c_1,c_2,\ldots,c_{N-3}$ and the relation \eqref{detHiggs_eq}. For example, if $c(z)$ has distinct roots $x_1=0,\ldots,x_{N-2}$, then $a_0,a_1,\ldots,a_{N-3}$ satisfy $N-2$ linear equations \[a(x_i)=y_i,~ i=1,\ldots, N-2,\text{ where }y_i=\pm (-\tilde{\nu}(x_i))^{1/2}.\] If $\mathrm{deg}\,c<N-2$, then $a_{N-2}^2= -\tilde{\nu}_{2N-4}$ and $\mathrm{deg}\,c=N-3$. In fact, if $a_{N-2}^2\neq -\tilde{\nu}_{2N-4}$, then $\mathrm{deg}(-a(z)^2-\tilde{\nu}(z))=2N-4>2N-5\geqslant\mathrm{deg}(b(z)c(z))$, a contradiction. Then we have $ \mathrm{deg}(-a(z)^2-\tilde{\nu}(z))=2N-5$, and $\mathrm{deg}\,c=2N-5-\mathrm{deg}\,b\geqslant N-3$, and $\mathrm{deg}\,c=N-3$. Applying gauge transformations as above we can make $c(z)$ monic and $a_{N-3}=0$. As before, the remaining coefficients of $a(z)$ are determined by $c_1,\ldots, c_{N-4}$. Therefore \[\mathrm{dim}(\mathcal{M}_{-|S|/2})=N-3.\]

 If $m>-|S|/2$, the holomorphic gauge transformations are of the form
 \[g=\begin{pmatrix}
   g_{11}&g_{12}(z)\\
   0&g_{11}^{-1}
 \end{pmatrix},\quad g_{11}\in \mathbb{C},g_{12}(z)\in \mathbb{C}[z], \mathrm{deg}\,g_{12}(z)\leqslant 2m+|S|.\]
 If $\mathrm{deg}\,c=-2m-|S|+N-2$, we may assume that $c(z)$ is monic after applying a diagonal gauge transformation as above. Then we successively apply
 \[g=\begin{pmatrix}
   1& a_{N-2-i}z^{2m+|S|-i}\\0&1
 \end{pmatrix},\quad i=0,1,\ldots, 2m+|S|,\]
 to make $a_{N-2}=a_{N-3}=\cdots=a_{N-2-2m-|S|}=0$. There is no gauge freedom left. The remaining $N-2-2m-|S|$ coefficients of $a(z)$ are determined by $c(z)$ and \eqref{detHiggs_eq}. If $\mathrm{deg}\,c<-2m-|S|+N-2$, then $a_{N-2}^2=-\tilde{\nu}_{2N-4}$ and $\mathrm{deg}\,c=-2m-|S|+N-3$. Again we can find a representative in each gauge orbit with $c(z)$ monic and $a_{N-3}=\cdots=a_{N-3-2m-|S|}=0$. The remaining coefficients of $a(z)$ are determined by $c(z)$. Therefore \[\mathrm{dim}(\mathcal{M}_{m})=-2m-|S|+N-2.\]

 \begin{remark}
   If there is a twisted order $m+3$ pole at $0$ with $\mu_{0,2m+5}\neq 0$ and $\mu_{0,a}=0$ for $m+4\leqslant a\leqslant 2m+4$, then the moduli space $\mathcal{M}$ is exactly $\mathcal{M}_{2,2m+1}$ in \cite{fredrickson2021moduli}. Such a moduli space admits a $\mathbb{C}^\ast$ action, and \eqref{Strata_eq} is a Białynicki–Birula-type stratification associated to this action \cite[Sec.s 1.4-1.5]{fredrickson2021moduli}.
 \end{remark}

 \section{Approximate Solutions}\label{approxsol_sec}
 \subsection{Harmonic metric}
 The hyperkähler metric $g_{L^2}$ on $\mathcal{M}$ is defined through the nonabelian Hodge correspondence, which realizes $\mathcal{M}$ as a hyperkähler quotient, i.e., unitary gauge equivalence classes of solutions to Hitchin's equations. To understand $g_{L^2}$, we need to study the solution of the following \emph{Hitchin’s equation} for a Higgs bundle $(\bar{\partial}_E,\varphi)$:
 \[F_h+[\varphi,\varphi^{\ast_h}]=0,\]
 where $F_h$ is the curvature of the Chern connection $D(\bar{\partial}_E,h)$. $h$ is called a \emph{harmonic metric} if it satisfies the above equation and is adapted to the parabolic structure, which means for each $x$ the filtration of the fiber $E_x$ given by $E_{x,\beta}:=\{s(0):s\text{ holomorphic},|s(z)|_{h}=O(|z|^\beta)\}$ coincides with that of the parabolic structure. When $I_t=\varnothing$, such a metric $h$ always exists on a stable Higgs bundle by \cite{biquard_boalch_2004}. However, if $I_t\neq \varnothing$, the existence of $h$ will impose further constraints on the parabolic weights $\alpha_{x,1}$, $\alpha_{x,2}$ for $x\in I_t$.
 \begin{lemma}\label{TIrregWeight_lem}
   Let $(\mathcal{E},\varphi)$ be an $\mathrm{SL}(2,\mathbb{C})$-Higgs bundle over a Riemann surface $C$, where $\varphi$ has an order $n\geqslant 2$ pole at $x$, with nilpotent leading coefficient $\phi_{-n}$. The kernel $L$ of $\phi_{-n}$ determines a filtration $0\subset L\subset E_x$ and let the associated parabolic weights be $1>\alpha_2>\alpha_1>0$, $\alpha_1+\alpha_2=1$. Suppose $\mathrm{pdeg}\,E=0$. If $h$ is a harmonic metric adapted to the parabolic structure, then $\alpha_1=1/4,\alpha_2=3/4$.
 \end{lemma}
 \begin{proof}
 As in \cite[Lem.~3.6]{mazzeo_swoboda_weiss_witt_2016}, we first find a standard form of $\varphi$ near $x$. Since $\phi_{-n}$ is nilpotent, $\det\varphi$ has an order $2n-1$ pole at $p_0$. By \cite[p.~28]{strebel_1984}, one can find a local coordinate $z$ centered at $p_0$, such that $\det\varphi=-z^{1-2n}\,\mathrm{d}z^2$. Again by the nilpotency of the leading term, in some local holomorphic frame
   \[z^n\phi(z)=\begin{pmatrix}
 a(z)&b(z)\\c(z)&-a(z)
 \end{pmatrix},\quad z^n\phi(z)|_{z=0}=\begin{pmatrix}
     0&0\\1&0
   \end{pmatrix}, \text{ where }\varphi(z)=\phi(z)\,\mathrm{d}z.\]
   Now $\sqrt{c(z)}$ is well-defined near $0$ since $c(0)=1$, and \[g^{-1}\varphi g=\begin{pmatrix}
     0&\frac{1}{z^{n-1}}\\\frac{1}{z^n}&0
   \end{pmatrix}\,\mathrm{d}z,\quad\text{where }g=\frac{1}{\sqrt{c(z)}}\begin{pmatrix}
     1&a(z)\\0&c(z)
   \end{pmatrix}.\]
   In the local holomorphic frame determined by $g$, we have \[h=\begin{pmatrix}
     h_{11}(z)&h_{12}(z)\\\widebar{{h}_{12}(z)}&h_{22}(z)
   \end{pmatrix},\quad h_{11}(z)\sim c_1r^{2\alpha_1},h_{22}(z)\sim c_2r^{2\alpha_2},\det h=Q(z)r^{2},h_{12}(z)=O(r),\]
   where $c_1,c_2$ are positive constants, $r=|z|$, $Q(z)$ is determined by $h_{\det E}$, the metric on $\det E$ induced by $h$, $Q(0)>0$. The $(1,1)$-entry of $F_h$ is $O(r^{-2})$ while that of $[\varphi,\varphi^{\ast_h}]$ is $O(r^{4\alpha_1-2n})+O(r^{4\alpha_2-2-2n})$ (the leading terms in these $O(\cdot)$'s are nonzero), which grows faster than $r^{-2}$ near $0$ if the sum of the leading terms is nonzero. Therefore $F_h+[\varphi,\varphi^{\ast_h}]=0$ implies that $4\alpha_1-2n=4\alpha_2-2-2n$, and $\alpha_1=1/4,\alpha_2=3/4$.
 \end{proof}

 Henceforth, we assume that $\alpha_{x,1}=1/4, \alpha_{x,2}=3/4$, for $x\in I_t$. When these conditions are fulfilled, $(\bar{\partial}_E,\varphi)$, which represents an element in $\mathcal{M}$, gives rise to a stable good filtered Higgs bundle in the sense of \cite{mochizuki21good}, and there exists a harmonic metric $h$ adapted to the parabolic structure by \cite[Th.~1.1]{mochizuki21good}. To get a more explicit description of the hyperkähler metric near the ends of $\mathcal{M}$, we need more information on $h$, which can be obtained via gluing constructions as in \cite{fredrickson2022asymptotic}, and we also follow this strategy.

 \subsection{A curve in $\mathcal{M}$}
 We consider the case $I_u\neq \varnothing$, and $0\in I_u$. Fix $\nu\in\mathcal{B}'$ satisfying the generic conditions $\tilde{\nu}_{2N-4}, \tilde{\nu}_{2N-5}\neq 0$, and $\nu_{N-4}\neq 0$. By \eqref{HitBase_eq}, we can represent $\nu$ by $(\nu_0,\ldots,\nu_{N-4})\in \mathbb{C}^{N-3}$. Consider a curve $\nu_t=(f_0(t)\nu_0,\ldots,f_{N-4}(t)\nu_{N-4})$ in $\mathcal{B}$ with $\nu_{t=1}=\nu$ and its norm tending to $\infty$, where $f_i(t)$ is positive and nondecreasing in $t$. For simplicity, let $f_{N-4}(t)=t$ and $f_i(t)=o(t^{(m_0+i)/(m_0+N-4)})$ for $i=0,\ldots,N-5$, so that the following holds.
 \begin{lemma}\label{HiggsDetRoot_lem}
   For $t$ sufficiently large, $\tilde{\nu}_t(z)$ has $2N-4$ distinct roots, with leading terms given by
 \begin{align*}
   z_{x,j}(t)&\sim x+ \bigg(-\frac{\mu_{x,2m_x}}{\nu_{N-4}}x^{4-N}\prod_{y\in S\backslash \{x\}}(x-y)^{m_y}\bigg)^{\frac{1}{j(x)}}\mathrm{e}^{\frac{2j\pi \mathrm{i}}{j(x)}}t^{-\frac{1}{j(x)}},~j=0,\ldots,j(x)-1,
 \end{align*}
 for $x\in I_u\cup P_w$ (replace $\mu_{x,2m_x}$ by $\mu_{x,2m_x-1}$ for $x\in I_t$), where
  \[j(x)=\begin{cases}
       m_0+N-4,&\text{ if }x=0,\\
       m_x,&\text{ if }x\in I_u\backslash\{0\}\cup P_w,\\
       m_x-1,&\text{ if }x\in I_t.
     \end{cases}\]
   $z_{x,m_x-1}(t)=x$ for $x\in I_t\cup P_s$.
 Here the notation $A\sim B$ means $A=B+o(B)$. Therefore $\nu_t\in \mathcal{B}'$ for $t$ large enough.
 \end{lemma}
 \begin{proof}
   By \eqref{HitBase_eq}, we can write
   \[\tilde{\nu}_t(z)=\bigg(\sum_{b=0}^{N-4}f_b(t)\nu_b z^b\bigg)z^{m_0}\prod_{y\in S\backslash \{0\}}(z-y)^{m_y}+\mu_{0,2m_0}\prod_{y\in S\backslash\{0\}}(z-y)^{2m_y}+zg_0(z),\]
   where $g_0(z)$ is some polynomial independent of $t$. By Rouché's theorem, there is a constant $\kappa$ such that in the disk $B_{\kappa t^{-1/(m_0+N-4)}}(0)$, $\tilde{\nu}_t(z)$ has the same number of zeros as that for
   \[\tilde{\nu}_t^{(0)}(z):=\nu_{N-4}tz^{m_0+N-4}\prod_{y\in S\backslash\{0\}}(-y)^{m_y}+\mu_{0,2m_0}\prod_{y\in S\backslash\{0\}}(-y)^{2m_y},\]
   for $t$ large enough. Choose $\kappa$, so that $\tilde{\nu}_t^{(0)}$ has $m_0+N-4$ distinct roots in the disk. Using the scaling $z\mapsto t^{1/(m_0+N-4)}z$, we see that the roots $z_{0,j}(t)$ ($j=0,\ldots,m_0+N-5$) of $\tilde{\nu}_t(z)$ in the disk are asymptotic to those of $\tilde{\nu}_t^{(0)}(z)$, which are given in the statement of the lemma. The analysis near $x\in S\backslash\{0\}$ is similar, using the coordinate $z_x$ in pace of $z$.
 \end{proof}

  Fix $[(\bar{\partial}_E,\varphi)]\in H^{-1}(\nu)$, and choose $[(\bar{\partial}_E,\varphi_t)]\in H^{-1}(\nu_t)$, such that $\varphi_t(z)$ and $\varphi(z)$ have the same $c(z)$ when written as in \eqref{GlobalHiggs_eq} (by the previous section, there are only finitely many gauge inequivalent choices of $\varphi_t$). Note that if $[(\bar{\partial}_E,\varphi)]$ belongs to $\mathcal{M}_m$ of \eqref{Strata_eq}, then the whole curve $[(\bar{\partial}_E,\varphi_t)]$ remains in $\mathcal{M}_m$ as $t$ varies.

 Now we aim to find a solution $h_t$ of $F_{h_t}+[\varphi_t,\varphi_t^{\ast_{h_t}}]=0$ for large $t$. First we construct an approximate solution $h_t^{\mathrm{app}}$, and later we will deform it into a genuine solution. The construction of $h_t^{\mathrm{app}}$ is similar to that in \cite[Sec.~3]{fredrickson2022asymptotic}.
 \subsection{Decoupled harmonic metric}\label{decoup_subsec}
 Since we are considering $\mathrm{SL}(2,\mathbb{C})$ Higgs bundles, the data on $\mathrm{det}\,E\simeq \mathcal{O}(-|S|)$ is fixed. In other words, we fix the holomorphic structure induced from $\bar{\partial}_E$ and the harmonic metric $h_{\mathrm{det}\,E}=\prod_{x\in S}|z-x|^2 $ adapted to the induced parabolic weight $\alpha_{x,1}+\alpha_{x,2}=1$ at $x$ (in $V$, $h_{\mathrm{det}\,E}=\prod_{x\in S}|1-xw|^2$). Let $(S_t,\mathcal{L}_t)$ be the spectral data associated with $(\bar{\partial}_E,\varphi_t)$. The spectral cover $\pi:S_t\to C$ ramifies at $Z_t$, the zero locus of $\tilde{\nu}_t(z)$. By \eqref{SpecLinedeg_eq}, $\mathrm{deg}\,\mathcal{L}_t=-|S|+N-2$. Now $\mathcal{L}_t$ has rank one, a parabolic structure at a point is equivalent to a parabolic weight since there is only a trivial filtration. For each of the two points in $\pi^{-1}(x)$, $x\in I_u\cup P_w$, let the weight be $\alpha_{x,1}$ if it corresponds to the eigenvalue $\rho_{x,m_{x}}$ of $\phi_{x,m_{x}}$ and otherwise be $\alpha_{x,2}$. For each point in $\pi^{-1}(Z_t\backslash(I_t\cup P_s))$, let the weight be $-1/2$. For each point $\pi^{-1}(x)$, $x\in I_t\cup P_s$, let the weight be $\alpha_{x,1}+\alpha_{x,2}-1/2=1/2$. Endowed with this parabolic structure we have $\mathrm{pdeg}\,\mathcal{L}_t=\mathrm{pdeg}\, E$. 
 By \cite{simpson_1990}, there is a unique (up to scaling) Hermitian-Einstein metric $h_{\mathcal{L}_t}$ on $\mathcal{L}_t$ adapted to the parabolic structure, which is flat in our setting. Finally, let $h_t^{\mathrm{dh}}$ be the metric on $E|_{C\backslash (Z_t\cup S)}$ defined as the orthogonal pushforward of $h_{\mathcal{L}_t}$. Then $h_t^{\mathrm{dh}}$ solves the decoupled Hitchin's equations:
 \[F_{h_t^{\mathrm{dh}}}=0,\quad [\varphi_t,\varphi_t^{\ast_{h_t^{\mathrm{dh}}}}]=0.\]
 $h_t^{\mathrm{dh}}$ is called a \emph{decoupled harmonic metric}. The metric on $\mathrm{det}\,E|_{C\backslash (Z_t\cup S)}$ induced by $h_t^{\mathrm{dh}}$ can be extended to $C\backslash\{S\}$ and is a harmonic metric. By uniqueness, we may assume that this metric is $h_{\mathrm{det}\,E}$ defined above, by rescaling $h_{\mathcal{L}_t}$ if necessary.

 \subsection{Normal forms}
 Near $x\in Z_t\cup S$, $\varphi_t$ and $h_t^{\mathrm{dh}}$ can be explicitly described using the following normal forms.
 \begin{proposition}\label{NormalForm_prop}
   There exist constants $1>\kappa,\kappa_1,\kappa_2>0$, such that for $t$ large enough, the followings hold.
   \begin{enumerate}[label=(\roman*)]
     \item There is a holomorphic coordinate $\zeta_{x,j,t}$ centered at $z_{x,j}(t)$ for $x\in I$, $j=0,\ldots,j(x)-1$ and a local holomorphic frame of $(E,\bar{\partial}_E)$ over $\widetilde{B}_{x,j,t}:=\{\,|\zeta_{x,j,t}|<\kappa t^{-1/j(x)} \,\}$ in which
     \begin{equation}\label{NormalFormZ_eq}
       \bar{\partial}_E=\bar{\partial},\quad \varphi_t=\lambda_{x,j}(t)\begin{pmatrix}
         0&1\\\zeta_{x,j,t}&0
       \end{pmatrix}\,\mathrm{d}\zeta_{x,j,t},\quad h_t^{\mathrm{dh}}=\begin{pmatrix}
         |\zeta_{x,j,t}|^{1/2}&0\\0&|\zeta_{x,j,t}|^{-1/2}
       \end{pmatrix}.
     \end{equation}
     Here the \emph{local mass} $\lambda_{x,j}(t)=|\nu_{x,j,t}(z_{x,j}(t))|^{1/2}$, where $\nu_{x,j,t}(z)=(z-z_{x,j}(t))^{-1}\nu_t(z)$.
 \item There is a holomorphic coordinate $\zeta_{x,t}$ centered at $x$ for $x\in P_w$, and a local holomorphic frame of $(E,\bar{\partial}_E)$ over $\widetilde{B}_{x,t}:=\{\,|\zeta_{x,t}|<\kappa t^{-2/3}\,\}$ in which
   \begin{equation}\label{NormalFormPw_eq}
     \bar{\partial}_E=\bar{\partial},\quad \varphi_t=\frac{1}{\zeta_{x,t}}\begin{pmatrix}
         0&-\mu_{x,2}(1-\epsilon_{x,0,t}^{-1}\zeta_{x,t})\\1&0
       \end{pmatrix}\,\mathrm{d}\zeta_{x,t},
   \end{equation}
   where $\epsilon_{x,0,t}=z_{x,0}(t)-x=O(1/t)$.
     \item There is a holomorphic coordinate $\zeta_{x,t}$ centered at $x$ for $x\in P_s$, and a local holomorphic frame of $(E,\bar{\partial}_E)$ over $\widetilde{B}_{x,t}:=\{\,|\zeta_{x,t}|<\kappa\,\}$ in which
   \begin{equation}\label{NormalFormPs_eq}
     \bar{\partial}_E=\bar{\partial},\quad \varphi_t=\lambda_{x}(t)\begin{pmatrix}
         0&1\\ \frac{1}{\zeta_{x,t}}&0
       \end{pmatrix}\,\mathrm{d}\zeta_{x,t},\quad h_t^{\mathrm{dh}}=\begin{pmatrix}
         |\zeta_{x,t}|^{1/2}&0\\0&|\zeta_{x,t}|^{3/2}
       \end{pmatrix}.
   \end{equation}
     Here $\lambda_{x}(t)=|\mathrm{Res}_x\nu_t(z)|^{1/2}$.
     \item There is a holomorphic coordinate $\zeta_{x,t}$ centered at $x$ for $x\in I_t$, and a local holomorphic frame of $(E,\bar{\partial}_E)$ over $\widetilde{B}_{x,t}:=\{\,|\zeta_{x,t}|<\kappa t^{-1/(m_x-1)}\,\}$ in which
 \begin{equation}\label{NormalFormIt_eq}
   \bar{\partial}_E=\bar{\partial},\quad \varphi_t=\lambda_{x}(t)\begin{pmatrix}
     0&\frac{1}{\zeta_{x,t}^{m_x-1}}\\ \frac{1}{\zeta_{x,t}^{m_x}}&0
   \end{pmatrix}\,\mathrm{d}\zeta_{x,t},\quad h_t^{\mathrm{dh}}=\begin{pmatrix}
     |\zeta_{x,t}|^{1/2}&0\\0&|\zeta_{x,t}|^{3/2}
   \end{pmatrix}.
 \end{equation}
     Here $\lambda_{x}(t)=|\nu_{x,t}(x)|^{1/2}$, and $\nu_{x,t}(z)=(z-x)^{2m_x-1}\nu_t(z)$.
     \item For $x\in I_u$, there is a local holomorphic frame of $(E,\bar{\partial}_E)$ over $\widetilde{B}_{x,t}:=\{\,|z-x|<\kappa t^{-1/j(x)}\,\}$ in which
     \begin{equation}\label{NormalFormIu_eq}
       \bar{\partial}_E=\bar{\partial},\quad \varphi_t=\frac{1}{z_x^{m_x}}(z_x^{2m_x}\nu_t(z))^{1/2}\sigma_3\,\mathrm{d}z,\quad h_t^{\mathrm{dh}}=\begin{pmatrix}
         |z_x|^{2\alpha_{x,1}}&0\\0&|z_x|^{2\alpha_{x,2}}
       \end{pmatrix}.
     \end{equation}
   \end{enumerate}
 The disks $\widetilde{B}_{x(,j),t}$ are disjoint, and $B_{x(,j),t,\kappa_1} \subset \widetilde{B}_{x(,j),t}\subset B_{x(,j),t,\kappa_2}$, where $B_\bullet$ are disks defined using the coordinate $z$, for example, $B_{x,j,t,\kappa_1}=\{\,|z-z_{x,j}(t)|<\kappa_1 t^{-1/j(x)}\,\}$, for $x\in I, j=0,\ldots,j(x)-1$.
 \end{proposition}
 \begin{proof}
   By \cite[Th.~6.1]{strebel_1984}, any quadratic differential having a critical point of odd order can be written in a standard form using some local holomorphic coordinate around that point. This can be applied to (\romannumeral1), (\romannumeral3), (\romannumeral4), where the quadratic differential $\det\varphi_t(z)$ has a critical point of order $1$, $-1$, $-(2m_x-1)$ respectively.

   (\romannumeral1) Let $z_{x,j,t}:=z-z_{x,j}(t)$. Consider the rescaling $z_{x,j,t}\mapsto (z_{x,j}(t)-x)^{-1}z_{x,j,t}:=\tilde{z}_{x,j,t}$. The quadratic differential $\det \varphi_t$ can be expressed in $\tilde{z}_{x,j,t}$ as \[\nu_{x,j,t}(z_{x,j}(t))\tilde{z}_{x,j,t}(z_{x,j}(t)-x)^3f_t(\tilde{z}_x)\,\mathrm{d}\tilde{z}_{x,j,t}^2.\] For $t$ large, $f_t(\tilde{z}_{x,j,t})$ is holomorphic on some disk $\{|\tilde{z}_{x,j,t}|<\kappa\}$, where it converges uniformly to a fixed holomorphic function $f(\tilde{z}_{x,j,t})$. By adapting the proof of \cite[Th.~6.1]{strebel_1984}, one can find a biholomorphic map $\tilde{\sigma}_{x,j,t}:\tilde{z}_{x,j,t}\mapsto\tilde{\zeta}_{x,j,t}$ from some disk $\{ |\tilde{z}_{x,j,t}|<\kappa_2 \}$ to its image, such that $\det \varphi_t=-\lambda_{x,j}(t)^2 (z_{x,j}(t)-x)^3 \tilde{\zeta}_{x,j,t}\,\mathrm{d}\tilde{\zeta}_{x,j,t}^2$, and $\tilde{\sigma}_{x,j,t}(0)=0$. $|\tilde{\zeta}_{x,j,t}/\tilde{z}_{x,j,t}|$ and its inverse are uniformly bounded. Let $\zeta_{x,j,t}:=(z_{x,j}(t)-x)\tilde{\zeta}_{x,j,t}$, and $\sigma_{x,j,t}$ be defined by $z_{x,j,t}\mapsto \zeta_{x,j,t}$. One can find constants $\kappa,\kappa_1$ so that $\widetilde{B}_{x,j,t}$ lies inside $\sigma_{x,j,t}(B_{x,j,t,\kappa_2})$ and its preimage $\sigma_{x,j,t}^{-1}(\widetilde{B}_{x,j,t})$ contains a smaller disk $B_{x,j,t,\kappa_1}$, for all sufficiently large $t$. In terms of the coordinate $\zeta_{x,j,t}$, we have $\det\varphi_t=-\lambda_{x,j}(t)^2\zeta_{x,j,t}\mathrm{d}\zeta_{x,j,t}^2$.

    As in \cite[Lem.~3.6]{mazzeo_swoboda_weiss_witt_2016}, we will write $\varphi_t$ in a standard form using local holomorphic gauge transformations, well defined over $\widetilde{B}_{x,j,t}$. Recall that $c(z)\neq 0$ in $\varphi_t$ is independent of $t$. For $t$ large, $c(z)$ is nonvanishing on $\widetilde{B}_{x,j,t}$, since otherwise $c(z_k)=0$ for a sequence $z_k\to x$, a contradiction. So we can write $\varphi_t$ in the coordinate $\zeta_{x,j,t}$ as $\varphi_t=\phi_{x,j,t}(\zeta_{x,j,t})\,\mathrm{d}\zeta_{x,j,t}$,   with $\phi_{x,j,t,21}(\zeta_{x,j,t})$ nonvanishing on $\widetilde{B}_{x,j,t}$, $\phi_{x,j,t,ij}(\zeta_{x,j,t})$ being the $(i,j)$-entry of $\phi_{x,j,t}(\zeta_{x,j,t})$. By a constant gauge transformation
   \[\begin{pmatrix}
     \phi_{x,j,t,11}(0)&1\\\phi_{x,j,t,21}(0)&0
   \end{pmatrix},\]
   we can make $\varphi_t=\phi^{(1)}_{x,j,t}(\zeta_{x,j,t})\,\mathrm{d}\zeta_{x,j,t}$, with
 \[\phi^{(1)}_{x,j,t}(0)=\begin{pmatrix}
   0&1\\0&0
 \end{pmatrix},\text{ and }\phi^{(1)}_{x,j,t,12}=\phi_{x,j,t,21}(\zeta_{x,j,t})/\phi_{x,j,t,21}(0).\]
 Now for $\kappa$ suitably small, $|\phi^{(1)}_{x,j,t,12}-1|<1/2$ and $\sqrt{\phi^{(1)}_{x,j,t,12}}$ is well defined on $\widetilde{B}_{x,j,t}$. Take the gauge transformation
 \[\frac{1}{\sqrt{\phi^{(1)}_{x,j,t,12}}}\begin{pmatrix}
   \phi^{(1)}_{x,j,t,12}&0\\-\phi^{(1)}_{x,j,t,11}&1
 \end{pmatrix}, \text{ such that }\varphi_t=\begin{pmatrix}
   0&1\\ \lambda_{x,j}(t)^2 \zeta_{x,j,t}&0
 \end{pmatrix}\,\mathrm{d}\zeta_{x,j,t}.\]
 Use the constant gauge transformation
 \[\begin{pmatrix}
   1&0\\0&\lambda_{x,j}(t)
 \end{pmatrix}, \text{ and then }\varphi_t=\lambda_{x,j}(t)\begin{pmatrix}
   0&1\\\zeta_{x,j,t}&0
 \end{pmatrix}\,\mathrm{d}\zeta_{x,j,t}.\]
 Finally, as in \cite[Prop.~3.5]{Fredrickson:2018fun}, $h_t^{\mathrm{dh}}$ has the form given above, after applying a holomorphic gauge transformation on $\widetilde{B}_{x,j,t}$ preserving $\varphi_t$.

 (\romannumeral2) We will find a holomorphic coordinate $\zeta_{x,t}$ centered at $x$ such that
 \begin{equation}\label{PwdetNormal_eq0}
 \det\varphi_t= \frac{\mathrm{d}\zeta_{x,t}^2}{\zeta_{x,t}^2}\mu_{x,2}(1-\epsilon_{x,0,t}^{-1}\zeta_{x,t})
 \end{equation}
 in $\widetilde{B}_{x,t}$. We can write \eqref{PwdetNormal_eq0} as
 \begin{equation}\label{PwdetNormal_eq}
 \frac{\mathrm{d}z_{x}^2}{z_{x}^2}(1-\epsilon_{x,0,t}^{-1}z_{x})f_t(z_{x}):=\mu_{x,2}^{-1}\det\varphi_t= \frac{\mathrm{d}\zeta_{x,t}^2}{\zeta_{x,t}^2}(1-\epsilon_{x,0,t}^{-1}\zeta_{x,t}).
 \end{equation}
 Let $\omega_{1,t}(z_x):=z_x^{-1}\sqrt{(1-\epsilon_{x,0,t}^{-1}z_{x})f_t(z_{x})}-z_x^{-1}$, which is holomorphic in $\{|z_x|<|\epsilon_{x,0,t}|\}$. Then the LHS of \eqref{PwdetNormal_eq} is $(\mathrm{d}N_t(z_x))^2$, where
 \[N_{1,t}(z_x)=\log z_x+\Omega_{1,t}(z_x),~ \Omega_{1,t}(z_x)=\int_0^{z_x}\omega_{1,t}(z)\,\mathrm{d}z. \]
 Similarly, the RHS of \eqref{PwdetNormal_eq} can be written as $(\mathrm{d}N_{2,t}(\zeta_{x,t}))^2$, where
 \begin{align*}
 N_{2,t}(\zeta_{x,t})&=\log \zeta_{x,t}+\Omega_{2,t}(\zeta_{x,t}),~\Omega_{2,t}(\zeta_{x,t})=\int_0^{\zeta_{x,t}} \omega_{2,t}(\zeta)\mathrm{d}\zeta,\\
 \omega_{2,t}(\zeta_{x,t})&=\zeta_{x,t}^{-1}\sqrt{1-\epsilon_{x,0,t}^{-1}\zeta_{x,t}}-\zeta_{x,t}^{-1}.
 \end{align*}
 Let $c_t=\Omega_{1,t}(\epsilon_{x,0,t})-\Omega_{2,t}(\epsilon_{x,0,t})$. Clearly, $\lim_{t\to\infty} c_t=0$. Exponentiating both sides of the equation $N_{1,t}(z_x)=c_t+N_{2,t}(\zeta_{x,t})$, we obtain
 \begin{equation}\label{PwdetNormalPot_eq}
 z_x\exp(\Omega_{1,t}(z_x)) =\zeta_{x,t}\exp(\Omega_{2,t}(\zeta_{x,t})+c_t).
 \end{equation}
 Consider the rescaled coordinates $\tilde{z}_{x,t}=\epsilon_{x,0,t}^{-1}z_x$ and $\tilde{\zeta}_{x,t}=\epsilon_{x,0,t}^{-1}\zeta_{x,t}$, and we can write \eqref{PwdetNormalPot_eq} in the form
 $F_t(\tilde{z}_{x,t})=G(\tilde{\zeta}_{x,t})$. $G'(0)\neq 0$, so there exists $\delta_1$ such that $G$ is invertible on $\{|\tilde{\zeta}_{x,t}|<2\delta_1\}$. Note that $\lim_{t\to\infty}F_t(z)=G(z)$ on $\{|z|<1\}$, then for $t$ large, we have
 \[F_t(\{|\tilde{z}_{x,t}|<\delta_1\})\subset  G(\{|\tilde{\zeta}_{x,t}|<2\delta_1\}).\]
 Let $\tilde{\sigma}_{x,t}:=G^{-1}\comp F_t, \tilde{z}_{x,t}\mapsto \tilde{\zeta}_{x,t}$, which is defined on $\{|\tilde{z}_{x,t}|<\delta_1\}$ and such that the following rescaled version of \eqref{PwdetNormal_eq} holds:
 \begin{equation}\label{PwdetNormalResc_eq}
 \frac{\mathrm{d}\tilde{z}_{x,t}^2}{\tilde{z}_{x,t}^2}(1-\tilde{z}_{x,t})\tilde{f}_t(\tilde{z}_{x,t})= \frac{\mathrm{d}\tilde{\zeta}_{x,t}^2}{\tilde{\zeta}_{x,t}^2}(1-\tilde{\zeta}_{x,t}).
 \end{equation}

 Next we show that $\tilde{\sigma}_{x,t}$ can be analytically extended to $\{|\tilde{\zeta}_{x,t}|<\delta_2\}$ for $t$ sufficiently large, where $\delta_2$ is some constant near $1$ to be specified later. We first determine $\tilde{\sigma}_{x,t}(\tilde{z}_{x,t})$ for $\tilde{z}_{x,t}\in U_{\delta_1,\delta_2}:=\{r\mathrm{e}^{\mathrm{i}\theta}\,|\,\delta_1\leqslant r<\delta_2,  -2\pi/3<\theta<2\pi/3\}$, using
 \begin{equation}\label{PwdetNormalPot_eq1}
 \int_{\delta_1/2}^{\tilde{z}_{x,t}}\frac{1}{z}\sqrt{(1-z)\tilde{f}_t(z)}\,\mathrm{d}z=\int_{\tilde{\sigma}_{x,t}(\delta_1/2)}^{\delta_1/2}+\int_{\delta_1/2}^{\tilde{z}_{x,t}}+\int_{\tilde{z}_{x,t}}^{\tilde{\sigma}_{x,t}(\tilde{z}_{x,t})}\frac{1}{\zeta}\sqrt{1-\zeta}\,\mathrm{d}\zeta.
 \end{equation}
 Here the integrals from $\delta_1/2$ to $\tilde{z}_{x,t}$ are taken along a curve in $U_{\delta_1/2,\delta_2}$, and are clearly independent of the choice of such a curve. Rewrite \eqref{PwdetNormalPot_eq1} as
 \begin{equation}\label{PwdetNormalPot_eq2}
 F_{1,t}(\tilde{z}_{x,t})=\int_{\tilde{z}_{x,t}}^{\tilde{\sigma}_{x,t}(\tilde{z}_{x,t})}\frac{1}{\zeta}\sqrt{1-\zeta}\,\mathrm{d}\zeta:=G_{\tilde{z}_{x,t}}(\tilde{\sigma}_{x,t}(\tilde{z}_{x,t})).
 \end{equation}
 $G_{\tilde{z}_{x,t}}'(\tilde{z}_{x,t})\neq 0$, so there exists $\delta_{x,t}$ such that $G_{\tilde{z}_{x,t}}$ is injective on $\{ |\zeta-\tilde{z}_{x,t}|<\delta_{x,t}\}$. For all $\tilde{z}_{x,t}\in \overline{U}_{\delta_1,\delta_2}$, we have $\delta_{x,t}\geqslant \underline{\delta}$ for some $\underline{\delta}>0$. Actually, suppose $G_{\tilde{z}_{x,t}}(\zeta_1)=G_{\tilde{z}_{x,t}}(\zeta_2)$ for $\zeta_1,\zeta_2\in \{ |\zeta-\tilde{z}_{x,t}|<\underline{\delta}\}$, then by Rolle's theorem for holomorphic functions, we have \[\mathrm{Re}(g(\zeta_3))=0,\quad \mathrm{Im}(g(\zeta_4))=0,\]
 for some $\zeta_3,\zeta_4\in\{ |\zeta-\tilde{z}_{x,t}|<\underline{\delta}\}$, where $g(\zeta)=\zeta^{-1}\sqrt{1-\zeta}$. The zero loci $Z_{\mathrm{Re}}$ and $Z_{\mathrm{Im}}$ of $\mathrm{Re}(g(\zeta))$ and  $\mathrm{Im}(g(\zeta))$ are independent of the choice of the branch of $\sqrt{1-\zeta}$. Explicitly,
 \[Z_{\mathrm{Re}}=\{ |\zeta-1|=1\}\cup [1,\infty),\quad  Z_{\mathrm{Im}}=(-\infty,1].\]
 Then one can find some $\underline{\delta}$ such that $\{ |\zeta-\tilde{z}_{x,t}|<\underline{\delta}\}$ does not simultaneously intersect with $Z_{\mathrm{Re}}$ and $Z_{\mathrm{Im}}$ for all $\tilde{z}_{x,t}\in \overline{U}_{\delta_1,\delta_2}$.  By the compactness of $\overline{U}_{\delta_1,\delta_2}$, there exists $\underline{\delta}_1$ such that $G_{\tilde{z}_{x,t}}(\{ |\zeta-\tilde{z}_{x,t}|<\underline{\delta}\})\supset \{ |\zeta|<\underline{\delta}_1\}$ for all $\tilde{z}_{x,t}\in \overline{U}_{\delta_1,\delta_2}$. Since $\lim_{t\to\infty}\tilde{f}_t(z)\to 1$ uniformly on $\{|z|<1\}$, for $t$ large, $|F_{1,t}(\tilde{z}_{x,t})|<\underline{\delta}_1$ on $\overline{U}_{\delta_1,\delta_2}$. Therefore we can define $\tilde{\sigma}_{x,t}(\tilde{z}_{x,t})$ on $U_{\delta_1,\delta_2}$ as $G_{\tilde{z}_{x,t}}^{-1}(F_{1,t}(\tilde{z}_{x,t}))$. Similarly, one can define $\tilde{\sigma}_{x,t}(\tilde{z}_{x,t})$ on $-U_{\delta_1,\delta_2}:=\{-\tilde{z}_{x,t}\in U_{\delta_1,\delta_2}\}$.

 Now we show that $\tilde{\sigma}_{x,t}$ can be extended across $1$. Let $\tau:z\mapsto z-1$ be the translation and $\hat{z}_{x,t}:=\tau(\tilde{z}_{x,t}), \hat{\zeta}_{x,t}:=\tau(\tilde{\zeta}_{x,t})$, then \eqref{PwdetNormalResc_eq} becomes
 \begin{equation}\label{PwdetNormal_eq3}
 \frac{\mathrm{d}\hat{z}_{x,t}^2}{(\hat{z}_{x,t}+1)^2}\hat{z}_{x,t}\hat{f}_t(\hat{z}_{x,t})= \frac{\mathrm{d}\hat{\zeta}_{x,t}^2}{(\hat{\zeta}_{x,t}+1)^2}\hat{\zeta}_{x,t}.
 \end{equation}
 By \cite[Th.~6.1]{strebel_1984}, there exists a biholomorphism $\hat{\sigma}_{2}: \hat{\zeta}_{x,t}\mapsto \xi_{x,t,2}$ from $\{
 |\hat{\zeta}_{x,t}|<2\delta_3\}$ to its image, unique up to a multiplicative factor $\exp(2k\pi\mathrm{i}/3)$ ($k=0,1,2$), such that the RHS of \eqref{PwdetNormal_eq3} equals $\xi_{x,t,2}\,\mathrm{d}\xi_{x,t,2}^2$. Similarly by the proof of \cite[Th.~6.1]{strebel_1984}, for $t$ large, there exists a biholomorphism $\hat{\sigma}_{t,1}: \hat{z}_{x,t}\mapsto \xi_{x,t,1}$ from a smaller disk $\{|\hat{z}_{x,t}|<\delta_3\}$ to its image, such that the LHS of \eqref{PwdetNormal_eq3} equals $\xi_{x,t,1}\,\mathrm{d}\xi_{x,t,1}^2$. $\hat{\sigma}_{t,1}$ is unique up to a multiplicative factor, which is chosen so that $-\pi/3\leqslant\mathrm{arg}(\hat{\sigma}_{t,1}(\delta_3/2)/\hat{\sigma}_2(\delta_3/2))<\pi/3$, and hence $\lim_{t\to\infty}\hat{\sigma}_{t,1}(z)=\hat{\sigma}_2(z)$ when $|z|<\delta_3$. We have $\hat{\sigma}_{t,1}(\{|\hat{z}_{x,t}|<\delta_3\})\subset \hat{\sigma}_2(\{|\hat{\zeta}_{x,t}|<2\delta_3\})$ for $t$ sufficiently large. Then we define $\hat{\sigma}_{x,t}:=\hat{\sigma}_2^{-1}\comp \hat{\sigma}_{t,1}$ and $\tilde{\sigma}_{x,t,1}:=\tau^{-1}\comp\hat{\sigma}_{x,t}\comp\tau$. Now $\tilde{\sigma}_{x,t,1}$ is defined on $\{|\tilde{z}_{x,t}-1|<\delta_3\}$, and \eqref{PwdetNormalResc_eq} holds for $\tilde{\zeta}_{x,t}=\tilde{\sigma}_{x,t,1}(\tilde{z}_{x,t})$. Recall that we have defined $\tilde{\sigma}_{x,t}$ on  $\{|\tilde{z}_{x,t}|<\delta_2\}$. Choose $\delta_2=1-\delta_3/2$, and we claim that $\tilde{\sigma}_{x,t,1}$ agrees with $\tilde{\sigma}_{x,t}$ on their common domain. It suffices to verify that $\tilde{\sigma}_{x,t,1}(1-3\delta_3/4)=\tilde{\sigma}_{x,t}(1-3\delta_3/4)$. For $t$ sufficiently large, $\tilde{\sigma}_{x,t}(1-3\delta_3/4),\tilde{\sigma}_{x,t,1}(1-3\delta_3/4)\in \{|z|<\delta_2\}$ and $|\tilde{\sigma}_{x,t}(1-3\delta_3/4)-\tilde{\sigma}_{x,t,1}(1-3\delta_3/4)|<\underline{\delta}$. By our constructions of $\tilde{\sigma}_{x,t,1}$ and $\tilde{\sigma}_{x,t}$, we have
 \[G_{\tilde{\sigma}_{x,t}(1-3\delta_3/4)}(\tilde{\sigma}_{x,t,1}(1-3\delta_3/4))=\int_{\tilde{\sigma}_{x,t}(1-3\delta_3/4)}^{\tilde{\sigma}_{x,t,1}(1-3\delta_3/4)} \frac{1}{\zeta}\sqrt{1-\zeta}=0. \]
 The claim follows from the injectivity of $G_{\tilde{\sigma}_{x,t}(1-3\delta_3/4)}$. We extend $\tilde{\sigma}_{x,t}$ to $\{|\tilde{z}_{x,t}|<\kappa t^{-2/3}|\epsilon_{x,0,t}|^{-1}\}$ for small $\kappa$ using the same method as extending it from $\{|\tilde{z}_{x,t}|<\delta_1\}$ to $\{|\tilde{z}_{x,t}|<\delta_2\}$. By rescaling we obtain $\sigma_{x,t}: z_x\mapsto \zeta_{x,t}$ as desired.

 We construct the inverse $\sigma_{x,t}^{-1}$ in a similar fashion. Recall that we expressed \eqref{PwdetNormalPot_eq} in the form  $F_t(\tilde{z}_{x,t})=G(\tilde{\zeta}_{x,t})$, and $G$ is invertible on $\{ |\tilde{\zeta}_{x,t}|<2\delta_1\}$. Let $2\epsilon_1=\min_{|\tilde{\zeta}_{x,t}|=2\delta_1}|G(\tilde{\zeta}_{x,t})|$. We have $F_t'(0)\neq 0$, $F_t(\tilde{z}_{x,t})\neq 0$ in $\{|\tilde{z}_{x,t}|<2\delta_1\}$ and $\min_{|\tilde{z}_{x,t}|=2\delta_1}|F_t(\tilde{z}_{x,t})|\geqslant \epsilon_1$ for $t$ large. By the inverse function theorem with an estimate on the range \cite[p.~234]{gamelin2001complex}, we have that $F_t(z)$ is invertible on $F_t^{-1}(\{|z|<\epsilon_1\})\cap \{|\tilde{z}_{x,t}|<2\delta_1\}$. Choose $\delta_1'$ such that $G(\{\tilde{\zeta}_{x,t}<\delta_1'\})\subset \{|z|<\epsilon_1\}$. Then we can define $\tilde{\sigma}_{x,t}^{-1}:\tilde{\zeta}_{x,t}\to\tilde{z}_{x,t}$ as $F_t^{-1}\comp G$ on $\{\tilde{\zeta}_{x,t}<\delta_1'\}$. As before, we extend $\tilde{\sigma}_{x,t}^{-1}$ to $\{\tilde{\zeta}_{x,t}<\delta_2'\}$ for some $\delta_2'$ close to $1$, using
 \begin{equation}\label{PwdetNormalPotInv_eq1}
 \int_{\delta_1'/2}^{\tilde{\zeta}_{x,t}}\frac{1}{\zeta}\sqrt{1-\zeta}\,\mathrm{d}\zeta=\int_{\tilde{\sigma}_{x,t}^{-1}(\delta_1'/2)}^{\delta_1'/2}+\int_{\delta_1'/2}^{\tilde{\zeta}_{x,t}}+\int_{\tilde{\zeta}_{x,t}}^{\tilde{\sigma}^{-1}_{x,t}(\tilde{\zeta}_{x,t})}\frac{1}{z}\sqrt{(1-z)\tilde{f}_t(z)}\,\mathrm{d}z.
 \end{equation}
 Rewrite the equation as
 \begin{equation}\label{PwdetNormalPotInv_eq2}
 G_{1,t}(\tilde{z}_{x,t})=\int_{\tilde{\zeta}_{x,t}}^{\tilde{\sigma}^{-1}_{x,t}(\tilde{\zeta}_{x,t})}\frac{1}{z}\sqrt{(1-z)\tilde{f}_t(z)}\,\mathrm{d}z:=F_{\tilde{\zeta}_{x,t}}(\tilde{\sigma}^{-1}_{x,t}(\tilde{\zeta}_{x,t})).
 \end{equation}
 The derivative of $F_{\tilde{\zeta}_{x,t}}$ is $g_{t,1}(z):=z^{-1}\sqrt{(1-z)\tilde{f}_t(z)}$. There exists $\underline{\delta}'$ such that the variation of the argument of $g_{t,1}$ over $\{|z-\tilde{\zeta}_{x,t}|<\underline{\delta}'\}$ is less than $\pi/2$ for any $\tilde{\zeta}_{x,t}\in \overline{U}_{\delta_1',\delta_2'}$, and then $F_{\tilde{\zeta}_{x,t}}$ is injective on this set. The remaining steps are the same.

 Write $\varphi_t$ in the coordinate $\zeta_{x,t}$ as $\varphi_t=\zeta_{x,t}^{-1}\phi_{x,t}\,\mathrm{d}\zeta_{x,t}$. As in (\romannumeral1), we can find a local holomorphic frame in which the $(2,1)$-entry of $\phi_{x,t}$ is $1$ and hence $\varphi_t$ has the form in \eqref{NormalFormPw_eq}.

 (\romannumeral3) is similar to \cite[Prop.~3.5]{fredrickson2022asymptotic}. (\romannumeral4) is similar to Lemma \ref{TIrregWeight_lem} and (\romannumeral3). (\romannumeral5) is similar to \cite[Prop.~3.7]{fredrickson2022asymptotic}, noting that $(z_x^{2m_x}\nu_t(z))^{1/2}$ is well defined on $\widetilde{B}_{x,t}$ for $\kappa$ suitably small. We can also choose $\kappa$ so that all the disks above are disjoint.
 \end{proof}

 \begin{remark}
   Intuitively, one may regard (\romannumeral3) as a limiting case of (\romannumeral2). When $t\to \infty$, the zero of $\varphi_t$ around $x\in P_w$ tends to $x$, yielding a nilpotent residue. More precisely, after a gauge transformation, we can express $\varphi_t$ around $x\in P_w$ as
   \[\epsilon_{x,0,t}^{1/2}\varphi_t=\frac{1}{\zeta_{x,t}}\begin{pmatrix}
     0&-\mu_{x,2}(\epsilon_{x,0,t}-\zeta_{x,t})\\ 1&0
   \end{pmatrix}\mathrm{d}\zeta_{x,t}\to \frac{1}{\zeta_{x,\infty}}\begin{pmatrix}
     0&\mu_{x,2}\zeta_{x,\infty}\\1&0
   \end{pmatrix}\mathrm{d}\zeta_{x,\infty},\]
   which has a strongly parabolic pole at $x$.
 \end{remark}

 \subsection{Approximate solutions}
 By (\romannumeral4) and (\romannumeral5) above, $h_t^{\mathrm{dh}}$ is already adapted to the parabolic structure at points in $I_u$ and $I_t$. In (\romannumeral1), $h_t^{\mathrm{dh}}$ is singular at $z_{x,j}(t)$, while in (\romannumeral3), $h_t^{\mathrm{dh}}$ is not adapted to the parabolic structure at $x\in P_s$. We desingularize $h_t^{\mathrm{dh}}$ near $Z_t\backslash(I_t\cup P_s)$ and $P_s$ using the following \emph{fiducial solutions} $h_{t}^{\mathrm{model}}$ \cite[Prop.s~3.8~\&~3.9]{fredrickson2022asymptotic}, and obtain the approximate metric $h_{t}^{\mathrm{app}}$. For $x\in P_w$, this kind of desingularization fails in a disk around $z_{x,0}(t)$ since the difference between the fiducial solution and $h_t^{\mathrm{dh}}$ does not decay to $0$ near the boundary when $t\to \infty$. Instead, we will use a model solution over $\widetilde{B}_{x,t}$ to be described later.
 \subsubsection{Fiducial Solutions}
 \begin{definition}\label{ApproxMet_def}
   \begin{enumerate}[label=(\roman*)]
     \item For $x\in I$, in the above holomorphic frame over $\widetilde{B}_{x,j,t}$, define
     \[h_{t}^{\mathrm{model}}=\begin{pmatrix}
       |\zeta_{x,j,t}|^{1/2}\mathrm{e}^{l_{x,j,t}(|\zeta_{x,j,t}|)}&0\\0&|\zeta_{x,j,t}|^{-1/2}\mathrm{e}^{-l_{x,j,t}(|\zeta_{x,j,t}|)}
     \end{pmatrix},\]
     where $l_{x,j,t}(r)=\psi_1(8\lambda_{x,j}(t)r^{3/2}/3)$, and $\psi_1$ satisfies the ODE
     \begin{equation}\label{PainleveODE_eq}
       \left(\frac{\mathrm{d}^2}{\mathrm{d}\rho^2}+\frac{1}{\rho}\frac{\mathrm{d}}{\mathrm{d}\rho}\right)\psi_1(\rho)=\frac{1}{2}\sinh(2\psi_1(\rho)),
     \end{equation}
     with asymptotics
     \[\psi_1(\rho)\sim \frac{1}{\pi}K_0(\rho)\text{ as }\rho\to\infty,\quad \psi_1(\rho)\sim -\frac{1}{3}\log\rho\text{ as }\rho\to 0.\]
   Here  $K_0$ is the modified Bessel function of the second kind, $K_0(\rho)\sim \rho^{-1/2}\mathrm{e}^{-\rho}$ as $\rho\to\infty$. Then we define
 \begin{equation}\label{approxMetricZ_eq}
   h_t^{\mathrm{app}}=\begin{pmatrix}
       |\zeta_{x,j,t}|^{1/2}\mathrm{e}^{\chi(|\zeta_{x,j,t}|\kappa^{-1}t^{1/j(x)})l_{x,j,t}(|\zeta_{x,j,t}|)}&0\\0&|\zeta_{x,j,t}|^{-1/2}\mathrm{e}^{-\chi(|\zeta_{x,j,t}|\kappa^{-1}t^{1/j(x)})l_{x,j,t}(|\zeta_{x,j,t}|)}
     \end{pmatrix},
 \end{equation}
   where $\chi$ is a nonincreasing smooth cutoff function with $\chi(r)=1$ for $r\leqslant 1/2$ and $\chi(r)=0$ for $r\geqslant 1$.
   \item For $x\in P_s$, in the above holomorphic frame over $\widetilde{B}_{x,t}$, define
   \[h_{t}^{\mathrm{model}}=\begin{pmatrix}
     |\zeta_{x,t}|^{1/2}\mathrm{e}^{m_{x,t}(|\zeta_{x,t}|)}&0\\0&|\zeta_{x,t}|^{3/2}\mathrm{e}^{-m_{x,t}(|\zeta_{x,t}|)}
   \end{pmatrix},\]
   where $m_{x,t}(r)=\psi_{2,x}(8\lambda_x(t)r^{1/2})$, and $\psi_{2,x}$ satisfies \eqref{PainleveODE_eq} with asymptotics
   \[\psi_{2,x}(\rho)\sim \frac{1}{\pi}K_0(\rho)\text{ as }\rho\to\infty,\quad \psi_{2,x}(\rho) \sim (1+2\alpha_{x,1}-2\alpha_{x,2})\log\rho\text{ as }\rho\to 0.\]
   Then we define
 \begin{equation}\label{approxMetricT_eq}
   h_t^{\mathrm{app}}=\begin{pmatrix}
       |\zeta_{x,t}|^{1/2}\mathrm{e}^{\chi(|\zeta_{x,t}|\kappa^{-1})m_{x,t}(|\zeta_{x,t}|)}&0\\0&|\zeta_{x,t}|^{3/2}\mathrm{e}^{-\chi(|\zeta_{x,t}|\kappa^{-1})m_{x,t}(|\zeta_{x,t}|)}
     \end{pmatrix}.
 \end{equation}
     \end{enumerate}
   According to \cite{mccoy1977painleve}, there exist unique functions $\psi_1$, $\psi_{2,x}$ solving \eqref{PainleveODE_eq} with the prescribed asymptotics.
 \end{definition}
 It is straightforward to check that $h_{t}^{\mathrm{model}}$ satisfies Hitchin's equation over $\widetilde{B}_{x(,j),t}$, smooth near $z_{x,j}(t)\in Z_t\backslash(I_t\cup P_s)$ and adapted to the parabolic structure at $x\in P_s$.
 \subsubsection{Model solution near $P_w$}\label{ModSol_subsubsec}
 Now we consider a Higgs bundle over  $\mathbb{C}P^1$ which is a model for $(\bar{\partial}_E,\varphi_t)|_{\widetilde{B}_{x,t}}$, where $x\in P_w$. Let $\hat{\mathcal{E}}=(\hat{E},\bar{\partial}_{\hat{E}})=\mathcal{O}(-1)\oplus\mathcal{O}(-1)$ be generated by $(e_1,e_2)$, and $\hat{\varphi}$ be the Higgs field defined in this frame by
 \[ \hat{\varphi}=\frac{\mathrm{d}z}{z}\begin{pmatrix} 0& -\mu_{x,2}(1-z)\\ 1&0\end{pmatrix}~\text{on}~ \mathbb{C}P^1\backslash{\{\infty\}},\quad  \hat{\varphi}=\frac{\mathrm{d}w}{w^2}\begin{pmatrix} 0& \mu_{x,2}(w-1)\\ -w&0\end{pmatrix}~\text{on}~ \mathbb{C}P^1\backslash{\{0\}}.\]
 The filtrations of $\hat{\mathcal{E}}$ at $0$ and $\infty$ are determined by $\varphi$ as in section \ref{irreg_subsec}. The parabolic weights are $\alpha_{x,1},\alpha_{x,2}$ at $0$ and $1/4,3/4$ at $\infty$. The Higgs bundle $(\bar{\partial}_{\hat{E}},\hat{\varphi})$ has a weakly parabolic pole at $0$ and an order $2$ twisted pole at $\infty$. This type of Higgs bundle is considered in
 \cite{harland22parabolic,mochizuki11wild,mochizuki21good}. By \cite[Th.~5.(\romannumeral1)]{harland22parabolic}, there exists a compatible harmonic metric $\hat{h}$. Define the decoupled harmonic $\hat{h}^{\mathrm{dh}}$ metric as in section \ref{decoup_subsec}, and let $\hat{\gamma}$ be defined by $\hat{h}(\cdot,\cdot)=\hat{h}^{\mathrm{dh}}(\cdot,\exp(-2\hat{\gamma})\cdot)$. Then $\hat{\gamma}$ is a smooth section of $\mathrm{i}\mathfrak{su}(E)$, where $\mathrm{i}\mathfrak{su}(E)$ is the bundle over $\mathbb{C}P^1\backslash\{0,1,\infty\}$ defined using $\hat{h}^{\mathrm{dh}}$. By \cite[Prop.~14.(a)]{harland22parabolic}, we have $|\hat{\gamma}|:=(\mathrm{tr}(\hat{\gamma}^2))^{1/2}$ is $O(|w|^\epsilon)$ near $\infty$ for some $\epsilon>0$. Moreover, by \cite[Th.~5.(\romannumeral2)]{harland22parabolic}, the curvature $F_{\hat{h}}$ of the Chern connection $D(\bar{\partial}_{\hat{E}},\hat{h})$ decays faster than any polynomial function in $|w|$ as $w\to 0$.

 \begin{lemma}\label{Fastdecay_lem}
   Let $\hat{\gamma}$ be defined as above. Then $|\hat{\nabla}^k\hat{\gamma}|=O(|w|^N)$ ($k=0,1,2$) as $w\to 0$ for any $N>0$, where $\hat{\nabla}$ is the Chern connection associated to $\bar{\partial}_{\hat{E}}$ and $\hat{h}^{\mathrm{dh}}$.
 \end{lemma}
 \begin{proof}
   For a section $\gamma$ of $\mathrm{i}\mathfrak{su}(E)$, define the following functionals
   \begin{align*}
     F_1(\gamma)&=F_{D(\mathrm{e}^{-\gamma}\bar{\partial}_E\mathrm{e}^{\gamma},\hat{h}^{\mathrm{dh}})}=\mathrm{e}^{-\gamma}(\bar{\partial}_{\hat{E}}(\mathrm{e}^{2\gamma}\partial_{\hat{h}^{\mathrm{dh}}}(\mathrm{e}^{-2\gamma})))\mathrm{e}^{\gamma},\\
     F_2(\gamma)&=[\mathrm{e}^{-\gamma}\hat{\varphi}\mathrm{e}^{\gamma},\mathrm{e}^{\gamma}\hat{\varphi}^{\ast_{\hat{h}^{\mathrm{dh}}}}\mathrm{e}^{-\gamma}].
   \end{align*}
   We have $F_1(0)=F_2(0)=0$ and $F_1(\hat{\gamma})+F_2(\hat{\gamma})=0$. Linearize $F_1,F_2$ and we obtain
   \begin{align*}
     L_1(\gamma)&=-2\bar{\partial}_{\hat{E}}\partial_{\hat{A}^{\mathrm{dh}}},\text{ where }\hat{A}^{\mathrm{dh}}=D(\bar{\partial}_{\hat{E}},\hat{h}^{\mathrm{dh}}),\\
     L_2(\gamma)&= M(\gamma),\text{ where }M(\gamma)=2[\hat{\varphi}^{\ast_{\hat{h}^{\mathrm{dh}}}},[\hat{\varphi},\gamma]].
   \end{align*}
   Similar to Proposition \ref{NormalForm_prop}.(\romannumeral4), in a neighborhood $W$ of $\infty$, there exists a holomorphic coordinate $\zeta=r\mathrm{e}^{\mathrm{i}\theta}$ centered at $\infty$ satisfying $c^{-1}|w|<|\zeta|<c|w|$ for some $c>0$ and a holomorphic trivialization $(e_1,e_2)$ of $\hat{\mathcal{E}}$ such that
   \[\hat{\varphi}=\begin{pmatrix}
     0& \frac{1}{\zeta}\\ \frac{1}{\zeta^2}&0
   \end{pmatrix}\mathrm{d}\zeta,\quad \hat{h}^{\mathrm{dh}}=\begin{pmatrix}
     r^{1/2}&0\\0&r^{3/2}
   \end{pmatrix}.\]
  In the unitary frame $(e_1/r^{1/4},e_2/r^{3/4})$, we have
   \[\hat{\varphi}=\begin{pmatrix}
     0& r^{-1/2} z^{-1}\\ r^{1/2}z^{-2}&0
   \end{pmatrix}\, \mathrm{d}\zeta,\quad \hat{h}^{\mathrm{dh}}=\mathrm{Id},\quad \mathrm{d}_{\hat{A}^{\mathrm{dh}}}=\mathrm{d}+\begin{pmatrix}
     \frac{1}{4}\mathrm{i}&0\\ 0&\frac{3}{4}\mathrm{i}
   \end{pmatrix}\mathrm{d}\theta.\]
   Since $\gamma=\gamma^{\ast}$, we have $\gamma=\left(\begin{smallmatrix}
     a&b\\\bar{b}&-a
   \end{smallmatrix}\right)$, where $a$ is real. Then
   \begin{align*}
     L_1(\gamma)=\frac{1}{2}\begin{pmatrix}
       \Delta a & \Delta_{-1/2} b\\  \overline{\Delta_{-1/2} b}& -\Delta a
   \end{pmatrix}\mathrm{d}\zeta \mathrm{d}\bar{\zeta},\quad L_2(\gamma)=\frac{-4}{r^3}\begin{pmatrix}
     2a& b-\mathrm{e}^{\mathrm{i}\theta}\bar{b}\\ \bar{b}-\mathrm{e}^{-\mathrm{i}\theta}b&-2a
   \end{pmatrix}\mathrm{d}\zeta \mathrm{d}\bar{\zeta},
   \end{align*}
   where $\Delta_{-1/2}:=\partial_r^2+r^{-1}\partial_r+r^{-2}(\partial_\theta-\mathrm{i}/2)^2$. The indicial roots of $L_1$ are integers for the diagonal terms, and constitute $\mathbb{Z}+1/2$ for the off-diagonal terms.

    By \cite[(6.1)]{fredrickson2022asymptotic}, we can write
   \begin{equation}\label{NonlinExpan_eq}
     F_1(\gamma)=B_1(\gamma)\bar{\partial}_{\hat{E}}\partial_{
     \hat{A}^{\mathrm{dh}}}\gamma+B_2(\gamma)\partial_{
     \hat{A}}\gamma \bar{\partial}_{\hat{E}}\gamma,\quad F_2(\gamma)=B_3(\gamma)(\hat{\varphi},\hat{\varphi}^{\ast}),
   \end{equation}
   where $B_1(\gamma)$ and $B_2(\gamma)$ are smooth functions of $\gamma$ and $B_3(\gamma)$ is a bilinear form with coefficients smooth in $\gamma$. We define the weighted spaces $r^\delta C_b^{k,\alpha}(W,\{\infty\};\mathrm{i}\mathfrak{su}(E))$ as in \cite[Sec.~5.1]{fredrickson2022asymptotic}. Then $\hat{\gamma}\in r^\epsilon C_b^{2,\alpha}$, where $0<\epsilon<1/2$ . Compose $L_1$ with $-\mathrm{i}\star$ to obtain $\tilde{L}_1$, which maps $\Gamma(\mathrm{i}\mathfrak{su}(E))$ to itself. By \eqref{NonlinExpan_eq}, $F_1(\hat{\gamma})-L_1(\hat{\gamma})\in r^{2\epsilon-2}C_b^{0,\alpha}$, and $F_1(\hat{\gamma})=O(r^N)$ for $N$ arbitrarily large, then $\tilde{L}_1(\hat{\gamma})\in r^{2\epsilon-2}C_b^{0,\alpha}$. Let $\hat{\gamma}=\left(\begin{smallmatrix}
     \hat{a}&\hat{b}\\\bar{\hat{b}}&-\hat{a}
   \end{smallmatrix}\right)$. By \cite[Prop.~5.1.(\romannumeral2)]{fredrickson2022asymptotic}, we have $\hat{a}=O(r^{2\epsilon})$, and $\hat{b}=O(r^{2\epsilon})$ if $2\epsilon<1/2$. Repeat this process if necessary, we may assume $2\epsilon>1/2$. Then we have $\hat{a}=O(r^{2\epsilon})$ and $\hat{b}=b_{-1}r^{1/2}+b_1r^{1/2}\mathrm{e}^{\mathrm{i}\theta}+O(r^{2\epsilon}) $, where $b_{\pm 1}\in \mathbb{C}$. Since $F_2(\hat{\gamma})=O(r^N)$, by considering its $(1,2)$-entry, we have
   \[b_{-1}r^{1/2}+b_1r^{1/2}\mathrm{e}^{\mathrm{i}\theta}-\overline{b_{-1}}r^{1/2}\mathrm{e}^{\mathrm{i}\theta}-\overline{b_1}r^{1/2}=O(r^{2\epsilon}),\]
 and hence $b_{-1}=\overline{b_1}$. Denote by $F_{j,n}$ ($j=1,2$) the order $n$ terms in the expansion of $F_j(\gamma)$ in $\gamma$. Then we have
 \begin{align*}
   F_{2,2}&=[[\hat{\varphi},\gamma],[\hat{\varphi},\gamma]^\ast]+\frac{1}{2}([[[\hat{\varphi},\gamma],\gamma],\hat{\varphi}^\ast]+[[[\hat{\varphi},\gamma],\gamma],\hat{\varphi}^\ast]^\ast)=0,\\
   F_{2,3}&=\frac{1}{6}([[[[\hat{\varphi},\gamma],\gamma],\gamma],\hat{\varphi}^\ast]+[[[[\hat{\varphi},\gamma],\gamma],\gamma],\hat{\varphi}^\ast]^\ast)\\&\quad+\frac{1}{2}([[[\hat{\varphi},\gamma],\gamma],[\hat{\varphi},\gamma]^\ast]+[[[\hat{\varphi},\gamma],\gamma],[\hat{\varphi},\gamma]^\ast]^\ast)\\&=\frac{4}{3r^3}
 \begin{pmatrix}
  F_{2,3,11} & F_{2,3,12} \\
  \widebar{F_{2,3,12}} & -F_{2,3,11}
 \end{pmatrix} \mathrm{d}\zeta \mathrm{d}\bar{\zeta},\\
 \text{ where }F_{2,3,11}&= a \big(-16 a^2+3 b^2 \mathrm{e}^{-\mathrm{i}\theta}+3 \bar{b}^2\mathrm{e}^{\mathrm{i}\theta} -10 b \bar{b}\big).
 \end{align*}
 From the $(1,1)$-entry of $F_2(\hat{\gamma})$, we see that
 \[-8\hat{a}+o(\hat{a})+O(r^2)=O(r^N),\]
 and then $\hat{a}=O(r^2)$. On the other hand,
 \begin{align*}
   F_{1,2}&=-2\gamma \bar{\partial}_{\hat{A}^{\mathrm{dh}}}\partial_{\hat{A}^{\mathrm{dh}}}\gamma-2(\bar{\partial}_{\hat{A}^{\mathrm{dh}}}\partial_{\hat{A}^{\mathrm{dh}}}\gamma)\gamma-4\bar{\partial}_{\hat{A}^{\mathrm{dh}}}\gamma\partial_{\hat{A}^{\mathrm{dh}}}\gamma\\&=\begin{pmatrix}
     F_{1,2,11}&F_{1,2,12}\\ F_{1,2,21}& F_{1,2,22}
 \end{pmatrix}\mathrm{d}\zeta \mathrm{d}\bar{\zeta},\\
 \text{ where }F_{1,2,11}&=a\Delta a+(b\overline{\Delta_{-1/2}b}+\bar{b}\Delta_{-1/2}b)/2+4(|\partial_{\bar{\zeta}}a|^2+|\partial_{\bar{\zeta}}b+(4\bar{\zeta})^{-1}b|^2).
 \end{align*}
 By considering the $(1,1)$-entry of $F_1(\hat{\gamma})$, we have
 \[|\partial_{\bar{\zeta}}\hat{b}+(4\bar{\zeta})^{-1}\hat{b}|^2+O(r^{2\epsilon-3/2})=O(r^N).\]
 Since $|\partial_{\bar{\zeta}}\hat{b}+(4\bar{\zeta})^{-1}\hat{b}|^2=|b_{-1}|^2r^{-1}/4+O(r^{2\epsilon-3/2})$, we have $b_{-1}=0$. Therefore $\hat{b}\in r^{\epsilon_1}C_b^{2,\alpha}$ for some $1/2<\epsilon_1<3/4$. As before, $\Delta_{-1/2}\hat{b}\in r^{2\epsilon_1-2}C_b^{0,\alpha}$ implies that $\hat{b}\in r^{2\epsilon_1}C_b^{2,\alpha}$. Repeat this process again, noting that $3/2$ is an indicial root of $\Delta_{-1/2}$, we have $\hat{b}=\overline{b_2}r^{3/2}\mathrm{e}^{-\mathrm{i}\theta}+b_2r^{3/2}\mathrm{e}^{2 \mathrm{i}\theta}+O(r^{\epsilon_2}) $ for some $3/2<\epsilon_2$. As before, we have $\hat{a}=O(r^6)$ and then $b_2=0$. By induction, $\hat{a},\hat{b}=O(r^N)$ for any $N$ and the lemma is proved.
 \end{proof}
 \begin{lemma}\label{PwmodelIso_lem}
   There exists an isomorphism $f_{x,t}:(\mathcal{E},\varphi_t,h_t^{\mathrm{dh}})|_{\widetilde{B}_{x,t}}\cong \hat{\sigma}_{x,t}^\ast(( \hat{\mathcal{E}},\hat{\varphi},\hat{h}^{\mathrm{dh}})|_{B'_{x,t}})$, where
   $B'_{x,t}:=\{|z|<\kappa t^{-2/3}|\epsilon_{x,0,t}|^{-1}\}$ and $\hat{\sigma}_{x,t}:\widetilde{B}_{x,t}\to \mathbb{C}P^1, \zeta_{x,t}\mapsto \epsilon_{x,0,t}^{-1}\zeta_{x,t}$.
 \end{lemma}
 \begin{proof}
   By Proposition \ref{NormalForm_prop} (\romannumeral2), there exists an isomorphism $(\mathcal{E},\varphi_t)|_{\widetilde{B}_{x,t}}\cong \hat{\sigma}_{x,t}^\ast(( \hat{\mathcal{E}},\hat{\varphi})|_{B'_{x,t}})$, denoted by $f_{x,t,1}$. Then $(\hat{\sigma}_{x,t}^{-1})^\ast h_t^{\mathrm{dh}}$ can be regarded as a decoupled harmonic metric on $(\hat{\mathcal{E}},\hat{\varphi})|_{B'_{x,t}}$ under the isomorphism $(\hat{\sigma}_{x,t}^{-1})^\ast\comp f_{x,t,1}$. $(\hat{\sigma}_{x,t}^{-1})^\ast h_t^{\mathrm{dh}}$ and $\hat{h}^{\mathrm{dh}}$ are orthogonal pushforwards under $\hat{\pi}$ of harmonic metrics $\hat{h}_1$, $\hat{h}_2$ on $\mathcal{L}$, where $\hat{\pi}:\mathcal{\hat{L}}\to\mathcal{\hat{E}}$ is the spectral cover. Denote $\hat{B}_{t,1}:=\hat{\pi}^{-1}(B'_{x,t})$, which is a ramified double covering of $B'_{x,t}$ and hence a disk. Let $s$ be a holomorphic trivialization of $\hat{\mathcal{L}}$, $\hat{\pi}^{-1}(1)=:\{p_0\}$ and $\hat{\pi}^{-1}(0)=:\{p_1,p_2\}$, then $\hat{h}_1(s,s)=\hat{u}\hat{h}_2(s,s)$, where $\log\hat{u}$ is a harmonic function on $\hat{B}_{t,1}\backslash\{p_0,p_1,p_2\}$. Let $\tilde{z}$ be a holomorphic coordinate on $\hat{B}_{t,1}$ and $\tilde{z}_i$ be the coordinates of $p_i$. By Bôcher's Theorem, there exist unique real numbers $a_i$ such that $\log\hat{u}(\tilde{z})=\mathrm{Re}(\hat{v}(\tilde{z}))+\sum_{i=1}^3 a_i\log|\tilde{z}-\tilde{z}_i|$, $\hat{v}$ is a holomorphic function on $\hat{B}_{t,1}\backslash\{p_0,p_1,p_2\}$. Since $\hat{h}_1,\hat{h}_2$ are compatible with the same parabolic data at $p_i$, we have $a_i=0$, and $\hat{v}$ is bounded, hence can be extended to $\hat{B}_{t,1}$. Therefore $(\hat{\sigma}_{x,t}^{-1})^\ast h_t^{\mathrm{dh}}=\exp(\mathrm{Re}(v(z)))\hat{h}^{\mathrm{dh}}$, where $v(z)$ is a holomorphic function on $B'_{x,t}$. Let $(e_1,e_2)$ be a holomorphic frame of $\hat{\mathcal{E}}=\mathcal{O}(-1)\oplus  \mathcal{O}(-1)$ over $B'_{x,t}$, and $(e_1,e_2)\mapsto (\exp(-v/2)e_1,\exp(-v/2)e_2)$ be the automorphism of $(\hat{\mathcal{E}},\hat{\varphi})|_{B'_{x,t}}$, whose composition with $f_{x,t,1}$ gives the desired $f_{x,t}$.
 \end{proof}
 Define $\hat{h}_t^{\mathrm{app}}(\cdot,\cdot)=\hat{h}(\cdot,\mathrm{exp}(2\chi(|z|\kappa^{-1}t^{2/3}|\epsilon_{x,0,t}|) \hat{\gamma})\cdot)$ on $B'_{x,t}$, and $h_t^{\mathrm{app}}=f_{x,t}^\ast\hat{\sigma}_{x,t}^\ast \hat{h}_t^{\mathrm{app}}$. Finally, set $h_{t}^{\mathrm{app}}=h_t^{\mathrm{dh}}$ on the complement of the above disks. $h_t^{\mathrm{app}}$ approximately solves Hitchin's equation in the following sense.
 \begin{proposition}
   For any $N>0$, we have
   \begin{equation}\label{ErrorEst_eq}
   \big\lVert F_{h_t^\mathrm{app}}+\big[\varphi_t,\varphi_t^{\ast_{h_t^\mathrm{app}}}\big] \big\rVert_{L^2}=O(t^{-N}).
 \end{equation}
 If $P_w=\varnothing$, then there exists positive constants $t_0$, $c,c'$, such that for any $t\geqslant t_0$,
   \begin{equation}\label{ErrorEst_eq1}
   \big\lVert F_{h_t^\mathrm{app}}+\big[\varphi_t,\varphi_t^{\ast_{h_t^\mathrm{app}}}\big] \big\rVert_{L^2}\leqslant c\mathrm{e}^{-c' t^{\sigma}},
 \end{equation}
 where the $L^2$-norm is taken with respect to $h_t^{\mathrm{app}}$ and a fixed Riemannian metric on $C$, and
 \[\sigma=\min\left(\frac{m_0-1}{m_0+N-4},\min_{x\in I_u\backslash\{0\}}\frac{m_x-1}{m_x},\min_{x\in I_t}\frac{2m_x-3}{2(m_x-1)},\frac{1}{2}\delta_{P_s\neq \varnothing}\right),\]
 $\delta_{P_s\neq \varnothing}=1$ if $P_s\neq \varnothing$, otherwise it equals $2$.
 \end{proposition}
 \begin{proof}
   Since $h_t^{\mathrm{app}}$ coincides with $h_{t}^{\mathrm{model}}$ and $h_t^{\mathrm{dh}}$ away from some annuli, solving Hitchin's equation, we only need to estimate the error in these regions. Around the zeros of $\varphi_t$, let $\mathrm{Ann}_{x,j}=\{\,\kappa t^{-1/j(x)}/2<|\zeta_{x,j,t}|<\kappa t^{-1/j(x)}\,\}$. Denote $|\zeta_{x,j,t}|$ by $r$ for simplicity.
   \begin{align*}
   F_{h_t^{\mathrm{app}}}&=-\frac{1}{4}\left(\Delta l_{x,j,t}(r)\chi\left(r\kappa^{-1}t^{1/j(x)}\right)\right)\sigma_3\, \mathrm{d}\zeta_{x,j,t}\mathrm{d}\bar{\zeta}_{x,j,t}\\
   &=-\frac{1}{4}\Big(l_{x,j,t}''(r)\chi\left(r\kappa^{-1}t^{1/j(x)}\right)+2\kappa^{-1}t^{1/j(x)}l_{x,j,t}'(r)\chi'\left(r\kappa^{-1}t^{1/j(x)}\right)\\&\quad+\kappa^{-2}t^{2/j(x)}l_{x,j,t}(r)\chi''\left(r\kappa^{-1}t^{1/j(x)}\right)\\&\quad+r^{-1}l_{x,j,t}'(r)\chi\left(r\kappa^{-1}t^{1/j(x)}\right)+r^{-1}\kappa^{-1}t^{1/j(x)}l_{x,j,t}(r)\chi'\left(r\kappa^{-1}t^{1/j(x)}\right)\Big)\sigma_3\,\mathrm{d}\zeta_{x,j,t}\mathrm{d}\bar\zeta_{x,j,t}.
 \end{align*}
 On the other hand,
 \[  \big[\varphi_t,\varphi_t^{\ast_{h_t^{\mathrm{app}}}}\big]=2\lambda_{x,j}(t)^2r\sinh\left(2l_{x,j,t}(r)\chi\left(r\kappa^{-1}t^{1/j(x)}\right)\right)\sigma_3\,\mathrm{d}\zeta_{x,j,t}\mathrm{d}\bar\zeta_{x,j,t}.\]
 By \eqref{PainleveODE_eq} and the definition of $l_{x,j,t}$, we have
 \[l''_{x,j,t}(r)+r^{-1}l_{x,j,t}'(r)=8\lambda_{x,j}(t)^2 r\sinh(2l_{x,j,t}),\quad l_{x,j,t}=O(\mathrm{e}^{-c'\lambda_{x,j}(t)t^{-3/2j(x)}}).\]
 We can write $F_{h_t^\mathrm{app}}+\big[\varphi_t,\varphi_t^{\ast_{h_t^\mathrm{app}}}\big]=E_{x,j,t}\sigma_3\,\mathrm{d}\zeta_{x,j,t}\mathrm{d}\bar\zeta_{x,j,t}$, the $L^2$ norm of which is bounded by $ct^{-1/j(x)}\sup|E_{x,j,t}|$. Here and below, $c,c'>0$ are generic constants which may vary, but does not depend on $t$. By Lemma \ref{HiggsDetRoot_lem}, as $t\to\infty$,
 \[\lambda_{x,j}(t)\sim t^{1/2}\prod_{y\in S}|z_{x,j}(t)-y|^{-m_x}\prod_{(y,k)\neq(x,j)}|z_{x,j}(t)-z_{y,k}(t)|^{1/2}\sim \begin{cases}
   c' t^{\frac{2m_x+1}{2j(x)}}&\text{ if }x\in I_u,\\
   c' t^{\frac{m_x}{j(x)}}&\text{ if }x\in I_t.
 \end{cases}\]
 Then for $x\in I_u$,
 \begin{align*}
   |E_{x,j,t}|&\leqslant c\left(\lambda_{x,j}(t)^2r(\left|\sinh(2l_{x,j,t})-\sinh(2l_{x,j,t}\chi)\right|+|\chi-1||\sinh(2l_{x,j,t})|)\right.\\&\left.\quad+t^{1/j(x)}|l'_{x,j,t}|+t^{2/j(x)}|l_{x,j,t}|\right)\\
   &\leqslant c\left(t^{2m_x/j(x)}\cosh(2l_{x,j,t})|l_{x,j,t}|+t^{1/j(x)}|l'_{x,j,t}|+t^{2/j(x)}|l_{x,j,t}|\right)\\&\leqslant c\left(t^{2m_x/j(x)}|l_{x,j,t}|+t^{2/j(x)}|l'_{x,j,t}|\right)\leqslant c \mathrm{e}^{-c' t^{(m_x-1)/j(x)}}.
 \end{align*}
 For $x\in I_t$, we have $|E_{x,j,t}|\leqslant c\mathrm{e}^{-c' t^{(2m_x-3)/(2m_x-2)}}$.

 For $x\in P_s$, let $\mathrm{Ann}_x=\{\,\kappa /2<|\zeta_{x,t}|<\kappa \,\}$, and denote $r=|\zeta_{x,t}|$. Similar as before,
 \begin{align*}
   F_{h_t^{\mathrm{app}}}  &=-\frac{1}{4}\Big(m_{x,t}''(r)\chi\left(r\kappa^{-1}\right)+2\kappa^{-1}m_{x,t}'(r)\chi'\left(r\kappa^{-1}\right)+\kappa^{-2}m_{x,t}(r)\chi''\left(r\kappa^{-1}\right)\\&\quad+r^{-1}m_{x,t}'(r)\chi\left(r\kappa^{-1}\right)+r^{-1}\kappa^{-1}m_{x,t}(r)\chi'\left(r\kappa^{-1}\right)\Big)\sigma_3\,\mathrm{d}\zeta_{x,t}\mathrm{d}\bar{\zeta}_{x,t},\\
   \big[\varphi_t,\varphi_t^{\ast_{h_t^{\mathrm{app}}}}\big]
  &=2\lambda_{x}(t)^2r^{-1}\sinh\left(2m_{x,t}(r)\chi\left(r\kappa^{-1}\right)\right)\sigma_3\,\mathrm{d}\zeta_{x,t}\mathrm{d}\bar{\zeta}_{x,t}.
 \end{align*}
 Now we have
 \[m_{x,t}''(r)+r^{-1}m_{x,t}'(r)=8\lambda_{x}(t)^2 r^{-1}\sinh(2m_{x,t}),\quad m_{x,j,t}=O(\mathrm{e}^{-c'\lambda_{x}(t)}),\]
 where $\lambda_{x}(t)\sim c' t^{1/2}$. For $E_{x,t}$ defined similarly as above, we have the estimate \[|E_{x,t}|\leqslant c\mathrm{e}^{-c' t^{1/2}}.\]

 Finally, for $x\in P_w$, let $\mathrm{Ann}_{x}=\{\,\kappa t^{-2/3}/2<|\zeta_{x,t}| <\kappa t^{-2/3} \,\}$, where $|E_{x,t}|=O(t^{-N})$ for any $N>0$ by \eqref{NonlinExpan_eq} and Lemmas \ref{Fastdecay_lem}, \ref{PwmodelIso_lem}.

 Combining all the estimates above, the result follows.
 \end{proof}
 \begin{remark}\label{WeakParabolic_rmk}
   If we use $h_t^{\mathrm{app}}$ in Definition \ref{ApproxMet_def}.(\romannumeral1) to desingularize $h_t^{\mathrm{dh}}$ around $z_{x,0}(t)$ for $x\in P_w$, then $\lambda_{x,j}(t)\sim c' t^{3/2}$, and $l_{x,j,t}=O(\mathrm{e}^{-c' t^{3/2}t^{-3/2}})=O(1)$. Hence the error does not decay to $0$ as $t\to\infty$.
 \end{remark}

 \section{Linear Analysis}
 We wish to obtain the solution $h_t$ to Hitchin's equation by perturbing $h_t^{\mathrm{app}}$. This is equivalent to finding an $\mathrm{SL}(E)$-valued $h_t^{\mathrm{app}}$-Hermitian section $H_t$ such that $h_t(v,w)=h_t^{\mathrm{app}}(H_tv,w)$. Let $\mathrm{e}^{\gamma_t}=H_t^{-1/2}$, then $\mathrm{e}^{\gamma_t}\cdot h_t(v,w)=h_t(H_t^{-1/2}v,H_t^{-1/2}w)=h_t^{\mathrm{app}}(v,w)$. By the gauge invariance of Hitchin's equations, $(\bar{\partial}_E,\varphi_t,h_t)$ is a solution is equivalent to that \[\mathrm{e}^{\gamma_t}\cdot (\bar{\partial}_E,\varphi_t,h_t)=(\mathrm{e}^{-\gamma_t}\bar{\partial}_E\mathrm{e}^{\gamma_t},\mathrm{e}^{-\gamma_t}\varphi_t \mathrm{e}^{\gamma_t},h_t^{\mathrm{app}})\]
 is a solution. Therefore, we need to find $\gamma_t\in \Omega^0(\mathrm{i}\mathfrak{su}(E))$ (the bundle is defined over $C\backslash S$ using $E$ and $h_t^{\mathrm{app}}$), such that $F_t(\gamma_t)=0$, where
 \[F_t(\gamma):=F_{A_t^{\exp(\gamma)}}+\big[\mathrm{e}^{-\gamma}\varphi_t\mathrm{e}^{\gamma},\mathrm{e}^{\gamma}\varphi_t^{\ast_{h_t^{\mathrm{app}}}}\mathrm{e}^{-\gamma}\big],\]
 $A_t$ is the Chern connection $D(\bar{\partial}_E,h_t^{\mathrm{app}})$ and $A_t^{\exp(\gamma)}=D(\mathrm{e}^{-\gamma}\bar{\partial}_E\mathrm{e}^{\gamma},h_t^{\mathrm{app}})$. Let $\mathrm{d}_{A_t}=\partial_{A_t}+\bar{\partial}_E$ be the exterior derivative associated with $A_t$, $\nabla_{A_t}$ be the covariant derivative. We have
 \[ F_{A_t^{\mathrm{exp}(\gamma)}}=\mathrm{e}^{-\gamma}\left(F_{A_t}+\bar{\partial}_{A_t}\left(\mathrm{e}^{2\gamma} \partial_{A_t}\left(\mathrm{e}^{-2\gamma}\right)\right)\right)\mathrm{e}^\gamma.\]
 Linearize $F_t$ at $\gamma=0$ we get
 \begin{align*}
   DF_t(\gamma)&=[F_{A_t},\gamma]-2\bar{\partial}_{A_t}\partial_{A_t}\gamma+[[\varphi_t,\gamma],\varphi_t^{\ast_{h_t^{\mathrm{app}}}}]+[\varphi_t,[\gamma,\varphi_t^{\ast_{h_t^{\mathrm{app}}}}]]\\&=\mathrm{i}\star\Delta_{A_t}\gamma+M_{\varphi_t}\gamma,
 \end{align*}
 where $\Delta_{A_t}=\mathrm{d}_{A_t}^\ast \mathrm{d}_{A_t}$ on $\Omega^0(\mathrm{i}\mathfrak{su}(E))$, $M_{\varphi_t}\gamma=[\varphi_t^{\ast_{h_t^{\mathrm{app}}}},[\varphi_t,\gamma]]-[\varphi_t,[\varphi_t^{\ast_{h_t^{\mathrm{app}}}},\gamma]]$, and we used the identity $2\bar{\partial}_A\partial_A=F_A-\mathrm{i}\star\Delta_A$. By composing $-\mathrm{i}\star: \Omega^2(\mathfrak{su}(E))\to\Omega^0(\mathrm{i}\mathfrak{su}(E))$ we get the operator
 \begin{equation}\label{LinOp_eq}
   L_t(\gamma)=\Delta_{A_t}\gamma-\mathrm{i}\star M_{\varphi_t}\gamma.
 \end{equation}
 Consider the irregular connection defined by $D_t^\gamma=\mathrm{d}_{A_t^{\mathrm{exp}(\gamma)}}+\mathrm{e}^{-\gamma}\varphi_t\mathrm{e}^{\gamma}+\mathrm{e}^{\gamma}\varphi_t^{\ast}\mathrm{e}^{-\gamma}$, here and below we drop the subscript of $\ast_{h_t^\mathrm{app}}$. Then its curvature $\left(D_t^\gamma\right)^2=F_t(\gamma)$, since $\bar{\partial}_E \varphi_t=\partial_{A_t}\varphi_t^\ast=0$. So we can also regard the problem as finding $\gamma$ so that this connection is flat. When $\gamma=0$, write the operator as $D_t$, which is the sum of the unitary part $\mathrm{d}_{A_t}$ and the self-adjoint part $\Psi_t:=\varphi_t+\varphi_t^\ast$.
 \begin{lemma}(Weitzenböck formulas)\label{weitz_lemma}
   Suppose $\gamma\in\Omega^0(\mathrm{i}\mathfrak{su}(E))$, $w\in \Omega^p(\mathrm{End}\,E)$, then \begin{enumerate}[label=(\roman*)]
     \item $(D_t^\ast D_t+D_tD_t^\ast)w=\nabla_{A_t}^\ast\nabla_{A_t} w+(\Psi_t\otimes )^\ast\Psi_t\otimes w+\mathscr{F}(w)+\mathscr{R}(w)$, where $\mathscr{F},\mathscr{R}$ are curvature operators, $\Psi_t\otimes\cdot$ is a combination of tensor product on the manifold part and Lie bracket on the bundle part. More explicitly, let $e_1,e_2$ be a local orthonormal frame of $TC$ and $e^1,e^2$ be the dual frame, write the curvature $D_t^2$ as $F_{t,12}e^1\wedge e^2$, then $\mathscr{F}(w)=[F_{t,12},(e^1\wedge i_{e_2}-e^2\wedge  i_{e_1})w]$, and
   \begin{align*}
     \mathscr{R}(w)(X_1,\ldots,X_p)&=\sum_{j=1}^pw(X_1,\ldots,\mathrm{Ric}(X_j),\ldots,X_p)\\&+\sum_{a=1}^2\sum_{j<k}(-1)^{j+k+1}w(e_a,R(X_j,X_k)e_a,X_1,\ldots,\hat{X}_j,\ldots,\hat{X}_k,\ldots,X_p).\hspace{-1.5cm}
   \end{align*}
     \item $D_t^\ast D_t\gamma=\mathrm{d}_{A_t}^\ast \mathrm{d}_{A_t}\gamma+([\Psi_t,\cdot])^\ast[\Psi_t,\gamma]=L_t\gamma$.
     \item $\mathrm{Re}\left(\langle \mathrm{d}_{A_t}\gamma,[\Psi_t,\gamma]\rangle\right)=0$, so the following pointwise equality holds
    \[\langle D_t\gamma,D_t\gamma\rangle= |D_t\gamma|^2=|\mathrm{d}_{A_t}\gamma|^2+|[\Psi_t,\gamma]|^2.\]
   \end{enumerate}
 \end{lemma}
 \begin{proof}
   (\romannumeral1) This follows from \cite[Theorem 5.4]{biquard_1997}, since $\mathrm{d}_{A_t}^\ast \Psi_t=0$.\\
   (\romannumeral2) The first equality is a special case of (\romannumeral1). The adjoint of $[\varphi_t,\cdot]$ is $-\bar{\ast}[\varphi_t,\bar{\ast}\cdot]$, where $\bar{\ast}$ is the conjugate linear operator combining $\bar{\star}$ and the Hermitian adjoint. Actually, for a fixed holomorphic coordinate $z$ where $\varphi=\phi\,\mathrm{d}z$ (subscript $t$ omitted),
   \begin{align*}
     \langle [\phi,a],b\,\mathrm{d}z\rangle&=\langle [\phi,a]\,\mathrm{d}z,b\,\mathrm{d}z\rangle=2\langle [\phi,a],b\rangle,\\
     \langle a,-\bar{\ast}[\phi,\bar{\ast} b\,\mathrm{d}z]\rangle&=\langle a,-\bar{\ast} [\phi,ib^\ast \mathrm{d}\bar{z}]\rangle=\langle a,-[b,\phi^\ast]\bar{\star} (\mathrm{i}\,\mathrm{d}z \mathrm{d}\bar{z})\rangle=2\langle a,[\phi^\ast,b]\rangle=2\langle [\phi,a],b\rangle.
   \end{align*}
 Similarly, the adjoint of $[\varphi_t^\ast,\cdot]$ is $-\bar{\ast}[\varphi_t^\ast,\bar{\ast}\cdot]$. Then (recall that $\bar{\ast}\mathrm{d}z=\mathrm{i}\,\mathrm{d}\bar{z}
 $, $\bar{\ast}\mathrm{d}\bar{z}=-\mathrm{i}\,\mathrm{d}z$)
 \begin{align*}
   ([\Psi_t,\cdot])^\ast[\Psi_t,\gamma]&=-\bar{\ast}[\varphi_t,\bar{\ast}[\varphi_t,\gamma]]-\bar{\ast}[\varphi_t^\ast,\bar{\ast}[\varphi_t^\ast,\gamma]]\\
   &=-\bar{\ast}[\varphi_t,\mathrm{i}[\varphi_t^\ast,\gamma]]+\bar{\ast}[\varphi_t^\ast,\mathrm{i}[\varphi_t,\gamma]]\\
   &=-\mathrm{i}\star [\varphi_t^\ast,[\varphi_t,\gamma]]+\mathrm{i}\star [\varphi_t,[\varphi_t^\ast,\gamma]]=-\mathrm{i}\star M_{\varphi_t}\gamma.
 \end{align*}
 (\romannumeral3) Write $\mathrm{d}_{A_t}\gamma=\gamma_z\,\mathrm{d}z+\gamma_{\bar{z}}\,\mathrm{d}\bar{z}$. $\gamma_{\bar{z}}=\gamma_z^\ast$, since $\gamma=\gamma^\ast$ and $\mathrm{d}_{A_t}$ is a Hermitian connection. As above, write $\varphi=\varphi\, \mathrm{d}z$, then \begin{align*}
   \mathrm{Re}\left(\langle \mathrm{d}_{A_t}\gamma,[\varphi+\varphi^\ast,\gamma]\rangle\right)&=2\mathrm{Re}\left(\mathrm{tr}(\gamma_z[\varphi,\gamma]^\ast)+\mathrm{tr}(\gamma_{\bar{z}}[\varphi^\ast,\gamma]^\ast)\right)\\
   &=2\mathrm{Re}\left(\mathrm{tr}(\gamma_z[\varphi,\gamma]^\ast)+\mathrm{tr}(\gamma_z^\ast[\gamma,\varphi])\right)\\
   &=2\mathrm{Re}\big(\mathrm{tr}(\gamma_z[\varphi,\gamma]^\ast)-\overline{\mathrm{tr}(\gamma_z[\varphi,\gamma]^\ast)}\big)=0.
   \end{align*}
 \end{proof}
 Choose a local holomorphic trivialization as in Proposition \ref{NormalForm_prop} (\romannumeral3). Consider the map $\sigma_{x,t}:D_{x,t}\to \widetilde{B}_{x,t},~\xi_{x,t}\to \xi_{x,t}^2=\zeta_{x,t}$, where $D_{x,t}=\{\,|\xi_{x,t}|<\kappa^{1/2}t^{-1/2(m_x-1)}\,\}$. Then for the pullbacks $\widetilde{E}=\sigma_{x,t}^\ast E$, $\widetilde{\varphi}_t=\sigma_{x,t}^\ast \varphi_t$, $\tilde{h}_t^{\mathrm{app}}=\sigma_{x,t}^\ast h_t^{\mathrm{app}}$, we have
 \[ \tilde{\bar{\partial}}_E=\bar{\partial},\quad \widetilde{\varphi}_t=2\lambda_x(t)\begin{pmatrix}
     0&\frac{1}{\xi_{x,t}^{2m_x-3}}\\\frac{1}{\xi_{x,t}^{2m_x-1}}&0
 \end{pmatrix}\,\mathrm{d}\xi_{x,t},\quad \tilde{h}_t^{\mathrm{app}}=\begin{pmatrix}
   |\xi_{x,t}|&0\\0&|\xi_{x,t}|^3
 \end{pmatrix}. \]
 After applying the gauge transformation $\left(\begin{smallmatrix}
   \xi_{x,t}&-\xi_{x,t}\\1&1
 \end{smallmatrix}\right)$ (noninvertible at $0$),
 \begin{equation}\label{DiagHiggs2Cover_eq}
   \widetilde{\varphi}_t=2\lambda_x(t)\begin{pmatrix}
       \frac{1}{\xi_{x,t}^{2(m_x-1)}}&0\\0&\frac{1}{\xi_{x,t}^{2(m_x-1)}}
   \end{pmatrix}\,\mathrm{d}\xi_{x,t},\quad \tilde{h}_t^{\mathrm{app}}=2\begin{pmatrix}
     |\xi_{x,t}|^3&0\\0&|\xi_{x,t}|^3
   \end{pmatrix}.
 \end{equation}
 Let $\hat{B}_{x,t}:=\{\,|\zeta_{x,t}|<\lambda_x(t)^{-2}\,\}$ for $x\in P_s$, and $\hat{B}_{x,t}:=\{\,|\zeta_{x,t}|<2|\epsilon_{x,0,t}|\,\}$ for $x\in P_w$. Choose $t_0\geqslant 1$ so that when $t\geqslant t_0$, we have $\lambda_x(t)^{-2}<\kappa/2$, and then $h_t^{\mathrm{app}}=h_t^{\mathrm{model}}$ over $\hat{B}_{x,t}$. From now on, we assume that $t\geqslant t_0$.
 \begin{definition}
   Let $r_t$ be a smooth weight function on $C$ such that
 \begin{enumerate}[label=(\roman*)]
   \item  $r_t=1$ in $C\backslash (\bigcup_{x\in I}\widetilde{B}_{x,t}\cup\bigcup_{x\in P}\hat{B}_{x,t})$.
   \item  In $\widetilde{B}_{x,t}$, $x\in S$, $r_t$ is nondecreasing in $|\zeta_{x,t}|$ (or $|z_x|$ for $x\in I_u$).
   \item \[r_t=\begin{cases}
       |t^{1/j(x)}z_x|\quad &\text{if }x\in I_u,~|z_x|\leqslant \kappa t^{-1/j(x)}/2,  \\
       |t^{1/(m_x-1)}\zeta_{x,t}|^{1/2}\quad &\text{if }x\in I_t,~|\zeta_{x,t}|\leqslant \kappa t^{-1/(m_x-1)}/2,\\
       |\lambda_x(t)^2\zeta_{x,t}|\quad &\text{if }x\in P_s,~|\zeta_{x,t}|\leqslant \lambda_x(t)^{-2}/2,\\
       |\epsilon_{x,0,t}^{-1}\zeta_{x,t}|\quad &\text{if }x\in P_w,~|\zeta_{x,t}|\leqslant |\epsilon_{x,0,t}|/2.
     \end{cases}\]
 \end{enumerate}
 \end{definition}
   Denote $B_x=\{\,|z_x|<\kappa_2\,\}$, where $\kappa_2$ is from Proposition \ref{NormalForm_prop}, and we can assume that $\kappa_2$ is small enough so that the disks $B_{x}$ are disjoint.
 \begin{definition}\label{Metricgt_def}
   Let $g_t=f_t g_C$ be the metric on $C$ conformal to the round metric $g_C$ (Hitchin's equations are conformally invariant, we consider $L_t$ above with $\star$ and $\Delta_{A_t}$ defined using $g_t$), which satisfies
 \begin{enumerate}[label=(\roman*)]
   \item $f_t=1$ in $C\backslash(\bigcup_{x\in S}B_x)$, independent of $t$ in $C\backslash(\bigcup_{x\in I_t\cup P}B_x)$.
   \item $g_t$ induces the volume form $(i/8|\zeta_{x,t}|)\,\mathrm{d}\zeta_{x,t} \mathrm{d}\bar{\zeta}_{x,t}$ when restricted to $\widetilde{B}_{x,t}$ for $x\in I_t$, its pullback under $\sigma_{x,t}$ is $(i/2)\,\mathrm{d}\xi_{x,t} \mathrm{d}\bar{\xi}_{x,t}$, the Euclidean volume form on $D_{x,t}$.
   \item $g_t$ is Euclidean in $B_x$ for $x\in I_u$.
   \item $g_t$ is Euclidean in $\widetilde{B}_{x,t}$ for $x\in P$.
   \item $|\nabla_{g_C}^j \log\,f_t|_{g_C}\leqslant c$ in $B_x$ for $x\in I_u\cup P$, $j=0,1,2$, and $c^{-1}|z_x|^{-1}\leqslant f_t\leqslant c|z_x|^{-1}$, $|\nabla_{g_C}^j \log\,f_t|_{g_C}\leqslant ct^{j/(m_x-1)}$ in $B_x\backslash \widetilde{B}_{x,t}$ for $x\in I_t$, $j=1,2$.
 \end{enumerate}
 \end{definition}

 \begin{definition}
   For smooth sections $w\in\Omega^p(\mathrm{End} E)$, define
 \begin{equation}\label{SobolevSpace_eq}
 \lVert w\rVert_{L_\delta^p}=\lVert r_t^{-\delta-2/p} w\rVert_{L^p},~  \lVert w\rVert_{L_\delta^{k,p}}^p=\lVert w\rVert_{L_\delta^{k-1,p}}^p+\lVert \nabla_{A_t}w\rVert_{L_\delta^{k-1,p}}^p+\lVert \Psi_t\otimes w\rVert_{L_\delta^{k-1,p}}^p
 \end{equation}
 inductively for $k=1,2\ldots$ Here the $L^p$ norms are taken with respect to $g_t$ and $h_t^{\mathrm{app}}$. Let the spaces $L_\delta^p$ and $L_\delta^{k,p}$ be the completions of smooth sections over $C\backslash S$ with respect to these norms. We also introduce spaces with additional weighted conditions,
 \[\hat{L}_{\delta}^{k,p}=\left\{\,f\in L_{\delta}^{k,p},r_t^{j-k}\nabla_{A_t}^j f\in L_\delta^p\text{ for }0\leqslant j\leqslant k\right\}.\]
 By using cutoff functions of the form $1-\chi(nr_t)$, one can see that $C_0^\infty(C\backslash S)$ is dense in $\hat{L}_\delta^{k,p}$.
 \end{definition}

 In $\widetilde{B}_{x,t}$ for $x\in I_u$, we can decompose $\mathrm{End}\, E$ as $\mathrm{End}_D\, E\oplus \mathrm{End}_T\, E$ in the holomorphic frame above, where $\mathrm{End}_D$ is the diagonal part. Then $\mathrm{End}_D\, E=\ker \mathrm{ad}(\varphi_t), \mathrm{End}_T\, E=\left(\ker \mathrm{ad}(\varphi_t)\right)^\perp$. $\varphi_t$ acts on $u\in\mathrm{End}_T\, E$ as
 \[[\varphi_t,u]= 2 \frac{\sqrt{-z_x^{2m_x}\nu_t(z)}}{z_x^{m_x}}\begin{pmatrix}
 0&u_{12}\\-u_{21}&0
 \end{pmatrix}\mathrm{d}z\Rightarrow |[\varphi_t,u]|=2\frac{|-z_x^{2m_x}\nu_t(z)|^{1/2}}{|z_x|^{m_x}}|u|,\]
 and then
 \begin{equation}\label{HiggsActBd_eq}
 c^{-1}t^{m_x/j(x)} |u|r_t^{-m_x}  \leqslant |[\varphi_t,u]|\leqslant ct^{m_x/j(x)} |u|r_t^{-m_x}\text{ in }\widetilde{B}_{x,t}.
 \end{equation}
  Therefore, over $\widetilde{B}_{x,t}$, the spaces defined above are equivalent (may not be uniformly in $t$) to the spaces defined in \cite{biquard_boalch_2004}, using the irregular Hermitian connection $A_t+\varphi_t-\varphi_t^\ast$ and the decomposition of $\mathrm{End}\, E$ into components $\mathrm{End}(E)_0=\mathrm{End}_D\, E$, $\mathrm{End}(E)_{m_x}=\mathrm{End}_T\, E$. We collect the following facts over $\widetilde{B}_{x,t}$ from \cite[Sec.~3]{biquard_boalch_2004}.
 \begin{lemma}\label{SobEmbed_lem}
 \begin{enumerate}[label=(\roman*)]
 \item For functions, we have $\hat{L}_{\delta-1}^{1,p}\hookrightarrow L_\delta^q$ if $1/2\geqslant 1/p-1/q$ (with strict inequality for $q=\infty$). In particular, for $p>2$ we have $\hat{L}_{\delta-1}^{1,p}\hookrightarrow C_\delta^0$.
 \item If $\delta<0$ and the function $f$ vanishes on $\partial \widetilde{B}_{x,t}$, or if $\delta>0$ and $f$ vanishes near $x$, then
 \[\lVert \partial f/\partial r\rVert_{L_{\delta-1}^p}\geqslant c^{-1}t^{1/j(x)}\lVert f\rVert_{L_{\delta}^p},\]
 where $r=|z_x|$. This implies that $\hat{L}_{\delta-1}^{1,p}=L_{\delta-1}^{1,p}$ for $\delta<0$ and $p\in[1,\infty)$.
 \item If $u\in L_{-2+\delta}^{2,2}$ is a section of $\mathrm{End}_D\, E$, then $u$ is continuous and $u-u(x)\in\hat{L}_{-2+\delta}^{2,2}\subset C_\delta^0\cap L_\delta^p$.
 \item For $\delta'<\delta$, $p>2$, there are compact Sobolev embeddings for sections of $\mathrm{End}_T\, E$: \[
 L_{\delta}^{1,2}\hookrightarrow L_{\delta'+1+(2m_x-2)/p}^p,\quad L_{\delta}^{1,p}\hookrightarrow C_{\delta'+m_x-2(m_x-1)/p}^0,\quad
 L_{\delta}^{2,2}\hookrightarrow C_{\delta'+m_x+1}^0.
 \]
 \end{enumerate}
 \end{lemma}

 Now we consider the Dirichlet problem on $\widetilde{B}_{x,t}$ and derive some a priori estimates. The proofs of the following two results are adapted from \cite[Cor.~4.2,~Lem.~4.4]{biquard_boalch_2004}, with the $t$-dependence of the operator $L_t$ taken into account.

 \begin{lemma}\label{aprioriEst_lem}
   Suppose $w\in L_{-2+\delta}^{1,2}(\Omega^1(\mathrm{End}_T\, E))$, which vanishes on $\partial \widetilde{B}_{x,t}$, then
   \[\lVert D_t w\rVert_{L_{-2+\delta}^2}+\lVert D_t^\ast w\rVert_{L_{-2+\delta}^2}\geqslant c^{-1}\big(\lVert\nabla_{A_t} w\rVert_{L_{-2+\delta}^2}+\lVert \Psi_t\otimes w\rVert_{L_{-2+\delta}^2}\big),\]
   for all weights $\delta$, when $t$ is sufficiently large.
 \end{lemma}
 \begin{proof}
    By Lemma \ref{weitz_lemma} (\romannumeral1), noting that the curvature terms vanish,
    \[\lVert (D_t+D_t^\ast)(r_t^{1-\delta}w)\rVert_{L^2}^2=\lVert \nabla_{A_t}(r_t^{1-\delta}w)\rVert_{L^2}^2+\lVert r_t^{1-\delta}\Psi_t\otimes w\rVert_{L^2}^2.\]
    The commutator can be controlled as
  \[|[D_t+D_t^\ast,r_t^{1-\delta}]w|+|[\nabla_{A_t},r_t^{1-\delta}]w|\leqslant ct^{1/j(x)}|r_t^{-\delta}w|.\]
  Then we have \[\lVert r_t^{1-\delta}(D_t+D_t^\ast)w\rVert_{L^2}\geqslant c^{-1}\left(\lVert r_t^{1-\delta}\nabla_{A_t} w\rVert_{L^2}+\lVert r_t^{1-\delta}\Psi_t\otimes w\rVert_{L^2}\right)-ct^{1/j(x)}\lVert r_t^{-\delta}w\rVert_{L^2}.\]
 On the other hand, \[\lVert r_t^{-\delta}w\rVert_{L^2}\leqslant \lVert r_t^{1-\delta}r_t^{-m_x}w\rVert_{L^2}\leqslant ct^{-m_x/j(x)} \lVert r_t^{1-\delta}\Psi_t\otimes w\rVert_{L^2}.\]
  Thus \begin{align*}
    \lVert r_t^{1-\delta}(D_t+D_t^\ast)w\rVert_{L^2}&\geqslant (c^{-1}-ct^{-(m_x-1)/j(x)})\left(\lVert r_t^{1-\delta}\nabla_{A_t} w\rVert_{L^2}+\lVert r_t^{1-\delta}\Psi_t\otimes w\rVert_{L^2}\right)\\
    &\geqslant (c^{-1}/2)\left(\lVert r_t^{1-\delta}\nabla_{A_t} w\rVert_{L^2}+\lVert r_t^{1-\delta}\Psi_t\otimes w\rVert_{L^2}\right),
  \end{align*}
  when $t$ is large enough.
 \end{proof}
 \begin{proposition}\label{analysisIu_prop}
   On $\widetilde{B}_{x,t}$ for $x\in I_u$, $L_t:L_{-2+\delta}^{2,2}(\mathrm{i}\mathfrak{su}(E))\to L_{-2+\delta}^{2}(\mathrm{i}\mathfrak{su}(E))$ is an isomorphism with zero Dirichlet boundary condition for small weights $\delta>0$. When restricted to the off-diagonal part, the same holds for any $\delta$. Moreover, the norm of $L_t^{-1}$ is uniformly bounded in $t$.
 \end{proposition}
 \begin{proof}
 Observe that the decomposition $\mathrm{End}\, E=\mathrm{End}_D\, E\oplus \mathrm{End}_T\, E$ is preserved by $L_t$, we can consider the two components separately. For the off-diagonal part, we find a solution of $L_tu=v$ by minimizing the functional (recall Lemma \ref{weitz_lemma} (\romannumeral3))
 \[S(u)=\frac{1}{2}\int_{\widetilde{B}_{x,t}} |D_t u|^2\,\mathrm{dvol}-\langle u,v\rangle_{L^2}=\frac{1}{2}\lVert \mathrm{d}_{A_t}u\rVert_{L^2}^2+\lVert[\varphi_t,u]\rVert_{L^2}^2-\langle u,v\rangle_{L^2}\]
 among $u\in L_{-1}^{1,2}(\mathrm{i}\mathfrak{su}_T(E))$ vanishing on the boundary, where $v\in L_{-m_x-1}^2(\mathrm{i}\mathfrak{su}_T(E))$. $S(u)$ is continuous: \[S(u)\leqslant \lVert \nabla_{A_t} u\rVert_{L^2}^2+\lVert \Psi_t\otimes u\rVert_{L^2}^2+\lVert r_t^{-m_x}u \rVert_{L^2}\lVert r_t^{m_x} v\rVert_{L^2}\leqslant \lVert u\rVert_{L_{-1}^{1,2}}^2+c\lVert u\rVert_{L_{-1}^{1,2}}\lVert v\rVert_{L_{-m_x-1}^2},\]
 where we used \eqref{HiggsActBd_eq} in the last inequality. On the other hand,
 \begin{align}
   S(u)&=\frac{1}{2}\big(\lVert \mathrm{d}_{A_t}u\rVert_{L^2}^2+2\lVert [\varphi_t,u]\rVert_{L^2}^2\big)-\langle u,v\rangle_{L^2}\notag\\&\geqslant c^{-1}\left(\lVert \nabla_{A_t} u\rVert_{L^2}^2+\lVert \Psi_t\otimes u\rVert_{L^2}^2+t^{2m_x/j(x)} \big\lVert r_t^{-m_x}u\right\rVert_{L^2}^2\big)- \lVert r_t^{-m_x}u\rVert_{L^2}\lVert v\rVert_{L_{-m_x-1}^2}\notag\\
   &\geqslant c^{-1}\lVert u\rVert_{L_{-1}^{1,2}}^2-c\lVert u\rVert_{L_{-1}^{1,2}}\lVert v\rVert_{L_{-m_x-1}^2}.\label{lowerbdS_eq}
 \end{align}
 Therefore $S(u)$ is coercive. If $\mathrm{d}_{A_t}u=0$, then $|u|^2$ is constant and $u=0$ since it vanishes on the boundary. This implies that $S$ is strictly convex, and a unique minimum of $S$ can be found, being a weak solution of the equation of the form \[\Delta_{A_t}\left( r_t^{m_x} u\right)=r_t^{m_x}v+P_0\left(r_t^{-m_x}u\right)+P_1\left(\nabla_{A_t} u\right),\]
 where $P_0,P_1:L^2\to L^2$ are bounded. By elliptic regularity, $r_t^{m_x} \nabla_{A_t}^2 u\in L^2$, and then $u\in L_{-m_x-1}^{2,2}$. With the Dirichlet boundary condition, we have obtained an isomorphism $L_t: L_{-m_x-1}^{2,2}\to L_{-m_x-1}^2$. The solution satisfies $\langle D_tu,D_tf\rangle_{L^2}=\langle v,f\rangle_{L^2}$ for any $f\in L_{-1}^{1,2}$ vanishing on the boundary. Choose $f=u$ (the first inequality follows by setting $v=0$ in \eqref{lowerbdS_eq}),
 \begin{align}
   \lVert u\rVert_{L_{-1}^{1,2}}^2\leqslant c\lVert D_tu\rVert_{L^2}^2= c\langle v,u\rangle_{L^2}&\leqslant ct^{-m_x/j(x)}\lVert \Psi_t\otimes u\rVert_{L^2}\lVert r_t^{m_x}v\rVert_{L^2}\notag\\&\leqslant ct^{-m_x/j(x)}\lVert u\rVert_{L_{-1}^{1,2}}\lVert v\rVert_{L_{-m_x-1}^2}\notag\\\Rightarrow~\lVert u\rVert_{L_{-1}^{1,2}}&\leqslant ct^{-m_x/j(x)}\lVert v\rVert_{L_{-m_x-1}^2}.\label{estimate1_eq}
 \end{align}
 For a one form $w$ with coefficient in $\mathrm{i}\mathfrak{su}(E)$, by Lemma \ref{aprioriEst_lem} we have \[\lVert\nabla_{A_t} w\rVert_{L_{-m_x-1}^2}\leqslant c\big(\lVert D_t w\rVert_{L_{-m_x-1}^2}+\lVert D_t^\ast w\rVert_{L_{-m_x-1}^2}\big).\]
 Then for $w=\nabla_{A_t}u$,
 \begin{align*}
   D_tw&=(\mathrm{d}_{A_t}+[\Psi_t,\cdot])(\mathrm{d}_{A_t}u)=[\Psi_t,\mathrm{d}_{A_t}u],\\
   D_t^\ast w&=D_t^\ast D_tu-D_t^\ast[\Psi_t,u]=D_t^\ast D_tu+\bar{\ast}[\bar{\ast}\Psi_t,\mathrm{d}_{A_t}u]-[\Psi_t,\cdot]^\ast[\Psi_t,u],\\
   \lVert r_t^{m_x} (\Psi_t\otimes\cdot)^2 u\rVert_{L^2}&\leqslant t^{m_x/j(x)}\lVert \Psi_t\otimes u\rVert_{L^2}\leqslant t^{m_x/j(x)}\lVert u\rVert_{L_{-1}^{1,2}}, \\
   \lVert r_t^{m_x}\Psi_t\otimes \nabla_{A_t}u\rVert_{L^2}&\leqslant t^{m_x/j(x)}\lVert \nabla_{A_t}u\rVert_{L^2}\leqslant t^{m_x/j(x)}\lVert u\rVert_{L_{-1}^{1,2}},\\
   \lVert r_t^{m_x}\nabla_{A_t}^2 u\rVert_{L^2}&\leqslant c\big(t^{m_x/j(x)}\lVert u\rVert_{L_{-1}^{1,2}}+\lVert D_t^\ast D_tu\rVert_{L_{-m_x-1}^2}\big)\leqslant c\lVert v\rVert_{L_{-m_x-1}^2},
 \end{align*}
 by \eqref{estimate1_eq}. Similarly, with $w=[\Psi_t,u]$ one deduces that $\lVert r_t^{m_x} \nabla_{A_t}[\Psi_t,u]\rVert_{L^2}\leqslant c\lVert v\rVert_{L_{-m_x-1}^2}$. Hence
 \begin{align*}
   \lVert u\rVert_{L_{-m_x-1}^{2,2}}&\leqslant c\big(\lVert u\rVert_{L_{-1}^{1,2}}+\lVert r_t^{m_x}\nabla_{A_t}^2u\rVert_{L^2}+\lVert r_t^{m_x}\nabla_{A_t}[\Psi_t,u]\rVert_{L^2}\\&\quad+\lVert r_t^{m_x}(\Psi_t\otimes\cdot)^2u\rVert_{L^2}+\lVert r_t^{m_x}\Psi_t\otimes \nabla_{A_t}u\rVert_{L^2}\big)\leqslant c\lVert v\rVert_{L_{-m_x-1}^2}.
 \end{align*}

 Next we prove that the isomorphism extends to all weights, by showing that the inverse is continuous in other weighted spaces. For any weight $\delta$, it suffices to prove
 \[\lVert r_t^\delta D_t^\ast D_t u\rVert_{L_{-m_x-1}^2}=\lVert D_t^\ast D_t u\rVert_{L_{-m_x-1-\delta}^2}\geqslant c^{-1}\lVert u\rVert_{L_{-m_x-1-\delta}^{2,2}},\]
 using the established a priori estimates
 \begin{equation}
   \lVert r_t^\delta u\rVert_{L_{-m_x-1}^{2,2}}\leqslant c\lVert D_t^\ast D_t r_t^\delta u\rVert_{L_{-m_x-1}^2}.\label{wtApriori_eq}
 \end{equation}
 We have
  \begin{align*}
   \lVert [D_t^\ast D_t,r_t^\delta]u\rVert_{L_{-m_x-1}^2}&\leqslant ct^{2/j(x)}\lVert r_t^{\delta-2} u\rVert_{L_{-m_x-1}^2}\\&\quad+ct^{1/j(x)}\big(\lVert r_t^{\delta-1} \nabla_{A_t} u\rVert_{L_{-m_x-1}^2}+\lVert r_t^{\delta-1} [\varphi_t,u]\rVert_{L_{-m_x-1}^2}\big)\\
   &\leqslant ct^{2/j(x)} \lVert r_t^\delta r_t^{-2m_x} u\rVert_{L_{-m_x-1}^2}\\&\quad+ct^{1/j(x)}\big(\lVert r_t^\delta r_t^{-m_x}\nabla_{A_t} u\rVert_{L_{-m_x-1}^2}+\lVert r_t^\delta r_t^{-m_x}[\varphi_t, u] \rVert_{L_{-m_x-1}^2}\big)\\
   &\leqslant c t^{(1-m_x)/j(x)} \big(\lVert r_t^\delta \Psi_t\otimes \nabla_{A_t}u\rVert_{L_{-m_x-1}^2}+\lVert r_t^\delta (\Psi_t\otimes \cdot)^2u\rVert_{L_{-m_x-1}^2}\big)\\
   &\leqslant ct^{(1-m_x)/j(x)} \lVert u\rVert_{L_{-m_x-1-\delta}^{2,2}}.
 \end{align*}
  Similarly, noting that $|[\nabla_{A_t},r_t^\delta]w|\leqslant c t^{(1-m_x)/j(x)}|r_t^\delta \Psi_t\otimes w|$,
 \begin{align*}
   \lVert u\rVert_{L_{-m_x-1-\delta}^{2,2}}&\leqslant \lVert r_t^\delta u\rVert_{L_{-m_x-1}^{2,2}}+\lVert  [\nabla_{A_t}^2,r_t^\delta]u\rVert_{L_{-m_x-1}^2}+\lVert [\nabla_{A_t},r_t^\delta]u\rVert_{L_{-m_x-1}^2}\\&\hspace{2.75cm}+\lVert \Psi_t\otimes [\nabla_{A_t},r_t^\delta]u\rVert_{L_{-m_x-1}^2}+\lVert [\nabla_{A_t},r_t^\delta]\Psi_t\otimes u\rVert_{L_{-m_x-1}^2} \\
   &\leqslant \lVert r_t^\delta u\rVert_{L_{-m_x-1}^{2,2}}+ct^{(1-m_x)/j(x)} \lVert u\rVert_{L_{-m_x-1-\delta}^{2,2}}.
 \end{align*}
 Finally by \eqref{wtApriori_eq},
 \begin{align*}
   \lVert r_t^\delta D_t^\ast D_t u\rVert_{L_{-m_x-1}^2}&\geqslant \lVert D_t^\ast D_tr_t^\delta u \rVert_{L_{-m_x-1}^2}-\lVert [D_t^\ast D_t,r_t^\delta]u\rVert_{L_{-m_x-1}^2}\\&\geqslant c^{-1}\lVert r_t^\delta u\rVert_{L_{-m_x-1}^{2,2}}-\lVert [D_t^\ast D_t,r_t^\delta]u\rVert_{L_{-m_x-1}^2}\\
   &\geqslant (c^{-1}-ct^{(1-m_x)/j(x)})\lVert u\rVert_{L_{-m_x-1-\delta}^{2,2}}\geqslant (c^{-1}/2)\lVert u\rVert_{L_{-m_x-1-\delta}^{2,2}},
 \end{align*}
 for $t$ large. We have proved that for the off-diagonal part, $L_t: L_{-2+\delta}^{2,2}\to L_{-2+\delta}^2$ is an isomorphism for all weights $\delta$ with the norm of the inverse uniformly bounded in $t$.

 For the diagonal part, the equation becomes the Laplace equation $\Delta u=v$ with the Dirichlet boundary condition. Let $\tau: \widetilde{B}_{x,t}\to B_{0}(\kappa):=\{\,|\tilde{z}|<\kappa\,\},~ z_x\mapsto t^{1/j(x)}z_x:=\tilde{z}$ be the rescaling map, and $\tilde{u}={(\tau^{-1})}^\ast u$, $\tilde{v}= {(\tau^{-1})}^\ast v$. Then $\Delta u=v$ is equivalent to $t^{2/j(x)}\Delta \tilde{u}=\tilde{v}$ on $B_{0}(\kappa)$, where $\Delta: L_{-2+\delta}^{2,2}\to L_{-2+\delta}^2$ with the Dirichlet boundary condition is an isomorphism by classical elliptic theory on cylinders (see \cite{biquard_boalch_2004}), and the weighted spaces on $B_{0}(\kappa)$ are defined similarly as in \eqref{SobolevSpace_eq}, using the Euclidean distance from $0$ as the weight function, which is equivalent to ${(\tau^{-1})}^\ast r_t$ uniformly in $t$. Then $\lVert  \tilde{v}\rVert_{L_{-2+\delta}^2}\geqslant c^{-1}\lVert t^{2/j(x)} \tilde{u}\rVert_{L_{-2+\delta}^{2,2}}$, which implies that $\lVert v\rVert_{L_{-2+\delta}^2}\geqslant c^{-1}\lVert u\rVert_{L_{-2+\delta}^{2,2}}$.
 \end{proof}

 Next we study the behavior of $L_t$ around $x\in I_t$. Define $L_\delta^{k,p}(D_{x,t}),\hat{L}_\delta^{k,p}(D_{x,t})$ as in \eqref{SobolevSpace_eq} using $\tilde{r}_t=\sigma_{x,t}^\ast r_t$, $\widetilde{\varphi}_t$, $\tilde{h}_t^{\mathrm{app}}$, the Euclidean metric on $D_{x,t}$ (which is $\sigma_{x,t}^{\ast}g_t$), and $\widetilde{A}_t=D(\tilde{\bar{\partial}}_E,\tilde{h}_t^{\mathrm{app}})$. For $w\in \Omega^q(\mathrm{End}\,E)$ we have \begin{equation}\label{NormRelation_eq}
      \lVert \sigma_{x,t}^\ast w\rVert_{L_\delta^{k,p}(D_{x,t})}=2\lVert w\rVert_{L_\delta^{k,p}(\widetilde{B}_{x,t})}.
    \end{equation}
 Let $\widetilde{L}_t(\gamma)=\Delta_{\widetilde{A}_t}\gamma-\mathrm{i}\star M_{\widetilde{\varphi}_t}\gamma$, then $\widetilde{L}_t(\sigma^\ast\gamma)=\sigma^\ast (L_t\gamma)$. The previous analysis near $x\in I_u$ yields the following.

 \begin{corollary}\label{analysisIt_cor}
   On $\widetilde{B}_{x,t}$ for $x\in I_t$, $L_t:L_{-2+\delta}^{2,2}(\mathrm{i}\mathfrak{su}(E))\to L_{-2+\delta}^{2}(\mathrm{i}\mathfrak{su}(E))$ is an isomorphism with zero Dirichlet condition on the boundary for small weights $\delta>0$. Moreover, the norm of $L_t^{-1}$ is uniformly bounded in $t$.
 \end{corollary}
 \begin{proof}
  Let $v\in L_{-2+\delta}^{2}(\mathrm{i}\mathfrak{su}(E))(\widetilde{B}_{x,t})$ then $\tilde{v}=\sigma_{x,t}^\ast v\in L_{-2+\delta}^2(\mathrm{i}\mathfrak{su}(E))(D_{x,t})$, and by Proposition \ref{analysisIu_prop} with $m_x$ and $j(x)$ both replaced by $2(m_x-1)$, there exists $\tilde{u}\in L_{-2+\delta}^{2,2}(\mathrm{i}\mathfrak{su}(E))(D_{x,t})$ such that $\widetilde{L}_t \tilde{u}=\tilde{v}$ and $\tilde{u}=0$ on $\partial D_{x,t}$. Note that $\widetilde{L}_t$ is $\mathbb{Z}_2$-equivariant, meaning that $\iota^\ast(\widetilde{L}_t (\tilde{u}))=\widetilde{L}_t (\iota^\ast \tilde{u})$ for the involution $\iota:\xi_{x,t}\mapsto -\xi_{x,t}$. Then $\widetilde{L}_t(\iota^\ast \tilde{u}-\tilde{u})=\iota^\ast \tilde{v}-\tilde{v}=0$, and $\iota^\ast \tilde{u}=\tilde{u}$ since $\widetilde{L}_t$ is injective. $\tilde{u}$ descends to $\widetilde{B}_{x,t}$ to give $u\in L_{-2+\delta}^{2,2}(\mathrm{i}\mathfrak{su}(E))(\widetilde{B}_{x,t})$ with $L_t u=v$ and $u=0$ on $\partial \widetilde{B}_{x,t}$. If $L_t u=0$, then $\widetilde{L}_t(\sigma_{x,t}^\ast u)=0$ and $\sigma_{x,t}^\ast u=0$ which implies that $u=0$. Therefore $L_t$ is an isomorphism. The uniform boundedness of $L_t^{-1}$ follows from that of $\widetilde{L}_t^{-1}$ and \eqref{NormRelation_eq}.
 \end{proof}
  The local behavior of $L_t$ around $x\in P_s$ has been studied in \cite{fredrickson2022asymptotic}, and we have the following result.
 \begin{proposition}\label{analysisPs_prop}
   On $\hat{B}_{x,t}$, for $x\in P_s$, $L_t:L_{-2+\delta}^{2,2}(\mathrm{i}\mathfrak{su}(E))\to L_{-2+\delta}^{2}(\mathrm{i}\mathfrak{su}(E))$ is an isomorphism with zero Dirichlet condition on the boundary for small weights $\delta>0$. Moreover, the norm of $L_t^{-1}$ is uniformly bounded in $t$.
 \end{proposition}
 \begin{proof}
 Consider $\tau:\hat{B}_{x,t}\to B_0(1):=\{\,|\tilde{z}|<1\,\},~\zeta_{x,t}\to \lambda_x(t)^2\zeta_{x,t}:=\tilde{z}$. Let $\widetilde{\varphi}_t=(\tau^{-1})^\ast\varphi_t$, and $\tilde{h}_t^{\mathrm{app}}=(\tau^{-1})^\ast \tilde{h}_t^{\mathrm{app}}$. Then by \eqref{NormalFormPs_eq} and \eqref{approxMetricT_eq},
 \[\widetilde{\varphi}_t=\begin{pmatrix}
   0&\frac{1}{\lambda_x(t)}\\ \frac{\lambda_x(t)}{\tilde{z}}&0
 \end{pmatrix}\,\mathrm{d}\tilde{z},\quad \tilde{h}_t^{\mathrm{app}}=\begin{pmatrix}
   \frac{1}{\lambda_x(t)}\tilde{r}^{1/2}\mathrm{e}^{\psi_{2,x} (8\tilde{r}^{1/2})}&0\\0&\frac{1}{\lambda_x(t)^3}\tilde{r}^{3/2}\mathrm{e}^{-\psi_{2,x}(8\tilde{r}^{1/2})}
 \end{pmatrix},\]
 where $\tilde{r}=|\tilde{z}|$. By the gauge transformation $g=\left(\begin{smallmatrix}
   \lambda_x(t)^{-1/2}&0\\0&\lambda_x(t)^{1/2}
 \end{smallmatrix}\right)$, we have
 \[\widetilde{\varphi}_t=\begin{pmatrix}
   0&1\\ \frac{1}{\tilde{z}}&0
 \end{pmatrix}\,\mathrm{d}\tilde{z},\quad \tilde{h}_t^{\mathrm{app}}=\frac{1}{\lambda_x(t)^2}\begin{pmatrix}
   \tilde{r}^{1/2}\mathrm{e}^{\psi_{2,x} (8\tilde{r}^{1/2})}&0\\0&\tilde{r}^{3/2}\mathrm{e}^{-\psi_{2,x}(8\tilde{r}^{1/2})}
 \end{pmatrix}.\]
 Then in this trivialization, $\widetilde{\varphi}_t$ and $\lambda_x(t)^2 \tilde{h}_t^{\mathrm{app}}$ are independent of $t$, which together with the Euclidean metric $g_e$ on $B_0(1)$ define an operator $\tilde{L}$ as in \eqref{LinOp_eq}. The indicial root analysis as in \cite[Sec.~4.2,~Lem. 5.8]{fredrickson2022asymptotic} implies that $\tilde{L}: L_{-2+\delta}^{2,2}(\mathrm{i}\mathfrak{su}(E))\to L_{-2+\delta}^{2}(\mathrm{i}\mathfrak{su}(E))$ (these spaces are defined on $B_0(1)$ as in \eqref{SobolevSpace_eq} using $g_e$, $\widetilde{\varphi}_t$, $\lambda_x(t)^2 \tilde{h}_t^{\mathrm{app}}$) is an isomorphism for $\delta>0$ small. Therefore $L_t$ is also an isomorphism, and the uniform boundedness of $L_t^{-1}$ follows in the same way as in Proposition \ref{analysisIu_prop} for the diagonal part.
 \end{proof}

 For $x\in P_w$, similar to the previous proposition, the pullback of $(\mathcal{E},\varphi_t,h_t^{\mathrm{app}})|_{\hat{B}_{x,t}}$ under $\hat{\sigma}_{x,t}^{-1}|_{B'}$ to $B':=\{|z|<2\}$ becomes a standard one. Using \cite[Sec.~4.3]{fredrickson2022asymptotic} or \cite[Lem.~4.4]{biquard_boalch_2004}, we obtain the following result.
 \begin{proposition}\label{analysisPw_prop}
   On $\hat{B}_{x,t}$, for $x\in P_w$, $L_t:L_{-2+\delta}^{2,2}(\mathrm{i}\mathfrak{su}(E))\to L_{-2+\delta}^{2}(\mathrm{i}\mathfrak{su}(E))$ is an isomorphism with zero Dirichlet condition on the boundary for small weights $\delta>0$. Moreover, the norm of $L_t^{-1}$ is uniformly bounded in $t$.
 \end{proposition}

 Let $U_{\mathrm{ext}}:=C\backslash (\bigcup_{x\in I} \frac{1}{2}\widetilde{B}_{x,t}\cup\bigcup_{x\in P}\frac{1}{2}\hat{B}_{x,t})$, where $\frac{1}{2}\widetilde{B}_{x,t}=\{\,|\zeta_{x,t}|< \kappa t^{-1/(m_x-1)}/2\,\}$ for $x\in I_t$, $\frac{1}{2}\widetilde{B}_{x,t}=\{\,|z_x|< \kappa t^{-1/j(x)}/2\,\}$ for $x\in I_u$, and $\frac{1}{2}\hat{B}_{x,t}$ for $x\in P$ are defined similarly. In $U_{\mathrm{ext}}$, $c^{-1}\leqslant r_t\leqslant 1$, which is immaterial in the analysis, so we omit the subscript of $L^{k,2}_{-2+\delta}$.

 \begin{proposition}\label{analysisUext_prop}
   On $U_{\mathrm{ext}}$, $L_t:L^{2,2}(\mathrm{i}\mathfrak{su}(E))\to L^2(\mathrm{i}\mathfrak{su}(E))$ is an isomorphism when imposing the zero Dirchlet boundary condition. The norm of $L_t^{-1}$ is bounded by $ct^4$.
 \end{proposition}
 \begin{proof}
 The invertibility is standard (see \cite[Sec.~5.2]{mazzeo_swoboda_weiss_witt_2016}). Let $u\in L^{2,2}(\mathrm{i}\mathfrak{su}(E))$ be vanishing on $\partial U_{\mathrm{ext}}$, then \begin{align*}
     \lVert \nabla_{A_t} u\rVert_{L^2}^2+\lVert [\Psi_t,u]\rVert_{L^2}^2 =\langle D_t^\ast D_t u,u\rangle_{L^2}\leqslant \lVert L_t u\rVert_{L^2}\lVert u\rVert_{L^2}.
   \end{align*}
 Let $\check{U}_{\mathrm{ext}}=\{|w|<2\kappa^{-1}t^{1/(m_0+N-4)}\}$, where $w=1/z$. Then $\check{U}_{\mathrm{ext}}$ is a disk centered at $\infty$ and contains $U_{\mathrm{ext}}$. Let $L_\mathrm{e}^2$ be the $L^2$ norm on $\check{U}_{\mathrm{ext}}$ defined using the Euclidean metric $g_e$ on this disk. By the Poincaré inequality and Kato's inequality, $\lVert u\rVert_{L_\mathrm{e}^2}\leqslant ct^{1/(m_0+N-4)}\lVert \nabla_{A_t} u\rVert_{L_\mathrm{e}^2}\leqslant ct^{1/2}\lVert \nabla_{A_t} u\rVert_{L_\mathrm{e}^2}$. By the Definition \ref{Metricgt_def} of $g_t$, we have $c^{-1}t^{-1}g_t\leqslant g_e\leqslant ct^{4/(m_0+N-4)} g_t\leqslant ct^2g_t$ on $U_{\mathrm{ext}}$. Therefore,
 \[\lVert u\rVert_{L^2}^2\leqslant ct\lVert u\rVert_{L_\mathrm{e}^2}^2\leqslant ct^2\lVert \nabla_{A_t}u\rVert_{L_\mathrm{e}^2}^2= ct^2\lVert \nabla_{A_t}u\rVert_{L^2}^2,\]
 where the last equality follows from the conformal invariance of the $L^2$ norm of a $1$-form. Then we have
   \begin{align*}
     \lVert u\rVert_{L^2}^2&\leqslant ct^2\lVert L_t u\rVert_{L^2}\lVert u\rVert_{L^2}\,\Rightarrow\, \lVert u\rVert_{L^2}\leqslant ct^2\lVert L_t u\rVert_{L^2},\\
     \lVert \nabla_{A_t}u\rVert_{L^2}&\leqslant ct\lVert L_tu\rVert_{L^2},~\lVert \Psi_t\otimes u\rVert_{L^2}\leqslant ct\lvert L_t u\rVert_{L^2}.
   \end{align*}
 We deduce that \begin{align*}
     \lVert (\Psi_t\otimes\cdot)^2 u\rVert_{L^2}&\leqslant c\sup |\varphi_t|\lVert \Psi_t\otimes u\rVert_{L^2}\leqslant ct\sup |\varphi_t|\lVert L_tu\rVert_{L^2},\\
     \lVert \Psi_t\otimes \nabla_{A_t}u\rVert_{L^2}&\leqslant c\sup |\varphi_t|\lVert  \nabla_{A_t}u\rVert_{L^2}\leqslant ct\sup |\varphi_t|\lVert L_tu\rVert_{L^2}.
     \end{align*}
   By Lemma \ref{weitz_lemma} (\romannumeral1), and the conformal change formula for $\mathscr{R}$,
   \begin{align*}
     \lVert \nabla_{A_t}^2 u\rVert_{L^2}&\leqslant c\left(\lVert D_t\nabla_{A_t} u\rVert_{L^2}+\lVert D_t^\ast \nabla_{A_t} u\rVert_{L^2}+t^3\lVert\nabla_{A_t} u\rVert_{L^2}+\lVert \mathscr{F}(\nabla_{A_t} u)\rVert_{L^2}\right)\\
     &\leqslant c\left(t\sup|\varphi_t|+t^4+t\max\sup |E_{x(,j),t}|\right)\lVert L_t u\rVert_{L^2},
   \end{align*}
   where the error terms $E_{x(,j),t}$ decay faster than any polynomial in $t$. In $U_{\mathrm{ext}}\backslash \big(\bigcup_{x\in I} \widetilde{B}_{x,j,t}\cup\bigcup_{x\in P}\hat{B}_{x,t}\big)$, $\varphi_t$ can be diagonalized as $\sqrt{-\nu_t(z)}\sigma_3\,\mathrm{d}z$, and $h_t^{\mathrm{app}}$ is also diagonal so that (recall that the norm $|\mathrm{d}z|$ is taken with respect to $g_t$) \[|\varphi_t|^2=2|\nu_t(z)||\mathrm{d}z|^2\leqslant ct^2.\]
   In $\widetilde{B}_{x,j,t}$, $\varphi_t$ and $h_t^\mathrm{app}$ are given by \eqref{NormalFormZ_eq} and \eqref{approxMetricZ_eq} in some holomorphic frame, so we have
   \begin{align*}
     |\varphi_t|^2&=2\lambda_{x,j}(t)^2r\cosh\left(2l_{x,j,t}(r)\chi\left(r\kappa^{-1}t^{1/j(x)}\right)\right)|\mathrm{d}\zeta_{x,j,t}|^2\\
     &\leqslant c\lambda_{x,j}(t)^{4/3} \rho^{2/3}\cosh(2\psi_1(\rho)\chi)|\mathrm{d}\zeta_{x,j,t}|^2\leqslant ct^2,
   \end{align*}
 where $r=|\zeta_{x,j,t}|, \rho=8\lambda_{x,j}(t)r^{3/2}/3$. In $\widetilde{B}_{x,t}\cap U_{\mathrm{ext}}$ for $x\in P_s$,  $\varphi_t$ and $h_t^\mathrm{app}$ are given by \eqref{NormalFormPs_eq} and \eqref{approxMetricT_eq} in some holomorphic frame, then
 \begin{align*}
   |\varphi_t|^2&=2\lambda_{x}(t)^2r\cosh\left(2m_{x,t}(r)\chi\left(r\kappa^{-1}\right)\right)|\mathrm{d}\zeta_{x,t}|^2\\
   &\leqslant c \rho^2\cosh(2\psi_{2,x}(\rho)\chi)|\mathrm{d}\zeta_{x,t}|^2\leqslant c,
 \end{align*}
 where $r=|\zeta_{x,t}|, \rho=8\lambda_{x}(t)r^{1/2}<8$. In $\hat{B}_{x,t}\cap U_{\mathrm{ext}}$ for $x\in P_w$, by considering the pullback from the model Higgs bundle in Section \ref{ModSol_subsubsec}, we have $|\varphi_t|^2\leqslant c$. Therefore, $\lVert \nabla_{A_t}^2 u\rVert_{L^2}\leqslant ct^4\lVert L_t u\rVert_{L^2}$. Similarly, applying the Weitzenböck formula for $w=[\Psi_t,u]$ gives $\lVert \nabla_{A_t} [\Psi_t,u] \rVert_{L^2}\leqslant ct^4\lVert L_t u\rVert_{L^2}$. Combining all the estimates above, we have \[\lVert u\rVert_{L^{2,2}(U_{\mathrm{ext}})}\leqslant ct^4\lVert L_t u\rVert_{L^2}.\]
 \end{proof}
 Finally, we consider the global behavior of $L_t$.
 \begin{lemma}\label{GlobIso_lem}
   $L_t:L_{-2+\delta}^{2,2}(\mathrm{i}\mathfrak{su}(E))\to L_{-2+\delta}^2(\mathrm{i}\mathfrak{su}(E))$ is an isomorphism for small $\delta>0$.
 \end{lemma}
 \begin{proof}
   By \cite[Lemma 5.1]{biquard_boalch_2004}, $L_t$ is Fredholm of index zero (the proof still work for $x\in I_t$ by considering the local lifted problem as in Corollary \ref{analysisIt_cor} and for $x\in P_s$ by \cite[Lem. 5.8]{fredrickson2022asymptotic}). For $u\in\ker L_t$, near $x\in I_u\cup P$ we have $|D_t u|=O(r^{\delta-1})$ ($r=|\zeta_{x,t}|$ for $x\in P$, and $r=|z_x|$ for $x\in I_u$) by the proof of \cite[Lem.~4.6]{biquard_boalch_2004}. The boundary term in the integration by parts is $\lim_{\epsilon\to 0^+}\int_{r=\epsilon} \langle \partial_r u,u\rangle r\,\mathrm{d}\theta$ \cite[Lem.~5.5]{fredrickson2022asymptotic} which vanishes for $x\in I_u\cup P$ since $\langle[\varphi_t,u],u\rangle=0$ and $\langle \partial_r u,u\rangle=O(r^{\delta-1})$, and similarly vanishes for $x\in I_t$ by considering the double lifting. Therefore \begin{equation}\label{intbyParts_eq}
   \langle L_t u,u\rangle_{L^2}=\lVert \mathrm{d}_{A_t}u\rVert_{L^2}^2+\lVert [\Psi_t,u]\rVert_{L^2}^2.
 \end{equation}
 Now we have $\mathrm{d}_{A_t}u=[\varphi_t,u]=0$, which implies that $u\equiv 0$ as in the proof of \cite[Lem.~5.6]{fredrickson2022asymptotic} since $Z_t\backslash (I_t\cup P_s)\neq \varnothing$. $L_t$ has trivial kernel, and then it is an isomorphism.
 \end{proof}

 \begin{lemma}\label{PoincareIneq_lem}Let $u\in L_{-2+\delta}^{2,2}(\mathrm{i}\mathfrak{su}(E))$, then for $t$ large enough,
   \begin{equation}\label{PoincareIneq_eq}
     \lVert u\rVert_{L_{-\delta}^2}^2\leqslant ct^{8/3}\left(\lVert \nabla_{A_t}u\rVert_{L^2}^2+\lVert [\varphi_t,u]\rVert_{L^2}^2\right).
   \end{equation}
 \end{lemma}
 \begin{proof}
 The idea is to use $[\varphi_t,u]$ to control $u$ around a zero of $\tilde{\nu}_t(z)$ near each $x\in I$ (we choose $z_{x,0}(t)$, which has asymptotics given in Lemma \ref{HiggsDetRoot_lem}), and use $\nabla_{A_t}u$ to control $u$ elsewhere. The inequality is obtained by gluing these estimates together.

 For $t$ large, $x\in I_u$, define $W_x:=\{\,|\zeta_{x,0,t}|<t^{-(2m_x+1)/3j(x)}\}\subset \frac{1}{2}\widetilde{B}_{x,0,t}$, a disk centered at $z_{x,0}(t)$. Then $\chi(r\kappa^{-1}t^{1/j(x)})=1$ on $W_x$, where $r=|\zeta_{x,0,t}|$, and \[\rho=8\lambda_{x,0}(t)r^{3/2}/3\leqslant c\text{ on }W_x.\]  Recall that $\varphi_t$ is given by \eqref{NormalFormZ_eq} in some holomorphic frame over $W_x$. Since $u$ is a section of $\mathrm{i}\mathfrak{su}(E)$, it has form
 \[u=\begin{pmatrix}
   p& \bar{s} r^{-1}\mathrm{e}^{-2l_{x,0,t}(r)}\\ s&-p
 \end{pmatrix},\]
 where $p$ is real valued. $\varphi_t$ acts on $u$ as
 \[[\varphi_t,u]=\lambda_{x,0}(t)\begin{pmatrix}
   s-\bar{s}\zeta_{x,0,t} r^{-1}\mathrm{e}^{-2l_{x,0,t}(r)}&-2p\\
   2\zeta_{x,0,t} p&-s+\bar{s}\zeta_{x,0,t} r^{-1}\mathrm{e}^{-2l_{x,0,t}(r)}
 \end{pmatrix}\,\mathrm{d}\zeta_{x,0,t}.\]
 Then on $W_x$ we have \begin{align}
   |u|^2&=2p^2+2|s|^2 r^{-1}\mathrm{e}^{-2l_{x,0,t}(r)}\leqslant 2p^2+c|s|^2t^{(2m_x+1)/3j(x)}\rho^{-2/3}\mathrm{e}^{-2\psi_1(\rho)}\notag\\&\leqslant 2p^2+ct^{(2m_x+1)/3j(x)}|s|^2,\notag  \\
   |[\varphi_t,u]|^2&= \lambda_{x,0}(t)^2\big(8p^2 r\cosh\left(2l_{x,0,t}(r)\right)+2\left|s-\bar{s}\zeta_{x,0,t} r^{-1}\mathrm{e}^{-2l_{x,0,t}(r)}\right|^2\big)\notag\\
   &\geqslant c^{-1}t^{(2m_x+1)/j(x)}\left(p^2t^{-(2m_x+1)/3j(x)}\rho^{2/3}\cosh\left(2\psi_1(\rho)\right)+|s|^2\left(1-\mathrm{e}^{-2\psi_1(\rho)}\right) \right)\notag\\&\geqslant c^{-1}t^{(4m_x+2)/3j(x)}|u|^2.\label{HiggsIuEst_eq}
 \end{align}

 Consider $\chi_x(r):=\chi(t^{(2m_x+1)/3j(x)}r)$, then $\chi_x=0$ in $C\backslash W_x$ and $\chi_x=1$ in $\hat{W}_x:=\frac{1}{2}W_x$. Let $r_{t,x}$ be the function obtained by extending $r_t|_{\widetilde{B}_{x,t}}$ by $1$, and $\hat{g}_t=\hat{f}_tg_C$ where $\hat{f}_t$ is the extension of $f_t|_{C\backslash \bigcup_{y\in I_t}B_y}$ by $1$. We claim that for all $v\in L_{-2+\delta}^{2,2}$ with $v$ vanishing in $\hat{W}_x$,
 \begin{equation}\label{hardyineqIu_eq}
   \lVert v\rVert_{L_{-\delta,x}^2}^2\leqslant ct^{(4m_x+2)/3j(x)}\lVert \nabla_{A_t}v\rVert_{L_{-\delta-1,x}^2}^2,
 \end{equation}
 where the subscripts indicate the norms are defined using $r_{t,x},\hat{g}_t$ instead of $r_t,g_t$. By Kato's inequality, we may assume that $v$ is a function, and prove the inequality with $\nabla_{A_t}v$ replaced by $\mathrm{d}v$. First consider the biholomorphism $\tau_x:C\to C_1=\mathbb{CP}^1$, defined by $z_x\mapsto c_xt^{1/j(x)}z_x$, where $c_x$ is the constant satisfying $\lim_{t\to\infty}\tau_x(z_{x,0}(t))=x+1$. Define weighted spaces on $C_1$ using the round metric on $C_1$ and the weight function $r_x$ given by the geodesic distance from $x$, then $r_x$ is uniformly equivalent to $(\tau_x^{-1})^\ast r_{t,x}$, and $c^{-1}t^{-2/j(x)}(\tau_x^{-1})^\ast \hat{g}_t\leqslant g_{C_1}$. Then for $v_1=(\tau_x^{-1})^\ast v$ we have
 \begin{align*}
   \lVert v\rVert_{L_{-\delta,x}^2}^2&=\int_{C} r_{t,x}^{2\delta-2}|v|^2\,\mathrm{dvol}_{\hat{g}_t}\leqslant c \int_{C_1} r_x^{2\delta-2}|v_1|^2\,(\tau_x^{-1})^\ast(\mathrm{dvol}_{\hat{g}_t})\\&\leqslant ct^{2/j(x)}\int_{C_1} r_x^{2\delta-2}|v_1|^2\,\mathrm{dvol}_{C_1}=ct^{2/j(x)}\lVert v_1\rVert_{L_{-\delta}^2(C_1)}^2,\\
   \lVert \mathrm{d}v_1\rVert_{L_{-\delta-1}^2(C_1)}^2&=\int_{C_1}r_x^{2\delta}|\mathrm{d}v_1|_{g_{C_1}}^2\,\mathrm{dvol}_{C_1}\leqslant c\int_C r_{t,x}^{2\delta}|\mathrm{d}v|_{\hat{g}_t}^2\,\mathrm{dvol}_{\hat{g}_t}=c\lVert \mathrm{d}v\rVert_{L_{-\delta-1,x}^2}^2.
 \end{align*}
 \eqref{hardyineqIu_eq} will follow from
 \begin{equation}\label{hardyineqIu_eq1}
   \lVert v_1\rVert_{L_{-\delta}^2(C_1)}^2\leqslant ct^{(4m_x-4)/3j(x)}\lVert \mathrm{d}v_1\rVert_{L_{-\delta-1}^2(C_1)}^2.
 \end{equation}
 Suppose $C_1$ is embedded in $\mathbb{R}^3$ so that $x=(0,1,0)$ and $x+1=(0,0,1)$. Consider the stereographic projection $\hat{\tau}_x:C_1\backslash \{x+1\}\to \mathbb{R}^2$, which sends $x+1$ to $\infty$ and $x$ to $(0,1)$. Note that as $t\to\infty$, $\tau_x(z_{x,0}(t))\to x+1$ and the diameter of $\tau_x(\hat{W}_x)$ is $O(t^{-(2m_x-2)/3j(x)})$, then for some constant $\kappa_x$, we have $\hat{\tau}_x^{-1}\left(\mathbb{R}^2\backslash C_2\right)\subset \tau_x(\hat{W}_x)$, where $C_2:=B_{(0,1)}(\kappa_xt^{(2m_x-2)/3j(x)})$. Define weighted spaces on $C_2$ using the Euclidean metric and the weight function $\hat{r}$ given by the Euclidean distance from $(0,1)$. Then $v_2=(\hat{\tau}_x^{-1})^\ast v_1$ vanishes on $\partial C_2$, $c^{-1}(\hat{\tau}_x^{-1})^\ast r_x\leqslant \hat{r}\leqslant ct^{(2m_x-2)/3j(x)}(\hat{\tau}_x^{-1})^\ast r_x$, and $c^{-1}(\hat{\tau}_x^{-1})^\ast g_{C_1}\leqslant g_{C_2}$ on $C_2$. We deduce that
 \begin{align*}
 \lVert v_1\rVert_{L_{-\delta}^2(C_1)}^2&=\int_{C_1}r_x^{2\delta-2}|v_1|^2\,\mathrm{dvol}_{C_1}\leqslant ct^{(2-2\delta)(2m_x-2)/3j(x)}\lVert v_2\rVert_{L_{-\delta}^2(C_2)}^2,\\
 \lVert \mathrm{d}v_2\rVert_{L_{-\delta-1}^2(C_2)}^2&=\int_{C_2}\hat{r}^{2\delta}|\mathrm{d}v_2|_{g_{C_2}}^2\,\mathrm{dvol}_{C_2}\leqslant ct^{(4m_x-4)\delta/3j(x)}\lVert \mathrm{d}v_1\rVert_{L_{-\delta-1}^2(C_1)}^2.
 \end{align*}
 An analogue of Lemma \ref{SobEmbed_lem} (\romannumeral2) implies that \[\lVert v_2\rVert_{L_{-\delta}^2(C_2)}^2\leqslant c\lVert \mathrm{d}v_2\rVert_{L_{-\delta-1}^2(C_2)}^2,\]
 since $v_2$ vanishes on $\partial C_2$. Therefore
 \begin{align*}
   \lVert v_1\rVert_{L_{-\delta}^2(C_1)}^2\leqslant c\lVert v_2\rVert_{L_{-\delta}^2(C_2)}^2\leqslant c\lVert \mathrm{d}v_2\rVert_{L_{-\delta-1}^2(C_2)}^2\leqslant ct^{(4m_x-4)\delta/3j(x)}\lVert \mathrm{d}v_1\rVert_{L_{-\delta-1}^2(C_1)}^2,
 \end{align*}
 verifying \eqref{hardyineqIu_eq1}, and \eqref{hardyineqIu_eq} holds.

 For $t$ large, $x\in I_t$, $W_x:=\{\,|\zeta_{x,0,t}|<t^{-2m_x/3j(x)}\,\}\subset \frac{1}{2}\widetilde{B}_{x,0,t}$, so that $\chi(r\kappa^{-1}t^{1/j(x)})=1$ and $\rho\leqslant c$ on $W_x$, where $r=|\zeta_{x,0,t}|$ and $\rho=8\lambda_{x,0}(t)r^{3/2}/3$. Similar to \eqref{HiggsIuEst_eq},
 \begin{equation}
   |[\varphi_t,u]|^2\geqslant c^{-1}t^{4m_x/3j(x)} |u|^2.\label{HiggsItEst_eq}
 \end{equation}
 Consider $\chi_x(r):=\chi(t^{2m_x/3j(x)}r)$, then $\chi_x=0$ in $C\backslash W_x$ and $\hat{\chi}=1$ in $\hat{W}_x:=\frac{1}{2}W_x$. Let $r_{t,x}$ be the function on $C$ which is the extension of $r_t|_{\widetilde{B}_{x,t}}$ by $1$, and $g_{t,x}=f_{t,x}g_C$ where $f_{t,x}$ is the extension of $f_t|_{C\backslash \bigcup_{y\in I_t\backslash\{x\}}B_y}$ by $1$. We claim that for all $v\in L_{-2+\delta}^{2,2}$ with $v$ vanishing in $\hat{W}_x$,
 \begin{equation}\label{hardyineqIt_eq}
   \lVert v\rVert_{L_{-\delta,x}^2}^2\leqslant ct^{4m_x/3j(x)}\lVert \nabla_{A_t}v\rVert_{L_{-\delta-1,x}^2}^2,
 \end{equation}
 where the subscripts indicate the norms are defined using $r_{t,x},g_{t,x}$ instead of $r_t$, $g_t$. As above we prove the inequality for $v$ being a function. Consider the biholomorphism $\tau_x: C\to C_1$, defined by $z_x\mapsto c_xt^{1/j(x)}z_x$, where $c_x$ is the constant chosen so that $\lim_{t\to\infty}\tau_x(z_{x,0}(t))=x+1$.  Define weighted spaces on $C_1$ using the singular metric $g_{C_1}=r_x^{-1}g_{\mathrm{FS}(C_1)}$ where $g_{\mathrm{FS}(C_1)}$ is the Fubini-Study metric on $C_1$, and the weight function $r_x^{1/2}$ where $r_x$ is the geodesic distance from $x$. Note that $r_x^{1/2}$ is uniformly equivalent to $(\tau_x^{-1})^\ast r_{t,x}$, and $c^{-1}t^{-2/j(x)}(\tau_x^{-1})^\ast g_{t,x}\leqslant g_{C_1}$. Then for $v_1=(\tau^{-1})^\ast v$ we have
 \[\lVert v\rVert_{L_{-\delta,x}^2}^2\leqslant ct^{2/j(x)}\lVert v_1\rVert_{L_{-\delta}^2(C_1)}^2,\quad \lVert \mathrm{d}v_1\rVert_{L_{-\delta-1}^2(C_1)}^2\leqslant c\lVert \mathrm{d}v\rVert_{L_{-\delta-1,x}^2}^2.\]
 \eqref{hardyineqIt_eq} follows once we prove
 \begin{equation}\label{hardyineqIt_eq1}
   \lVert v_1\rVert_{L_{-\delta}^2(C_1)}^2\leqslant ct^{(4m_x-6)/3j(x)}\lVert \mathrm{d}v_1\rVert_{L_{-\delta-1}^2(C_1)}^2.
 \end{equation}
 Suppose $C_1$ is embedded in $\mathbb{R}^3$ as above and consider the stereographic projection $\hat{\tau}_x:C_1\backslash \{x+1\}\to \mathbb{R}^2$, which sends $x+1$ to $\infty$ and $x$ to $(0,1)$. Note that the diameter of $\tau_x(\hat{W}_x)$ is $O(t^{-(2m_x-3)/3j(x)})$, then for some constant $\kappa_x$, we have $\hat{\tau}_x^{-1}\left(\mathbb{R}^2\backslash C_2\right)\subset \tau_x(W_2)$, where $C_2:=B_{(0,1)}(\kappa_xt^{(2m_x-3)/3j(x)})$. Define weighted spaces on $C_2$ using the singular metric $g_{C_2}=\hat{r}^{-1}g_{e(C_2)}$ where $g_{e(C_2)}$ is the Euclidean metric, and the weight function $\hat{r}^{1/2}$ where $\hat{r}$ is the Euclidean distance from $(0,1)$. Then $v_2=(\hat{\tau}_x^{-1})^\ast v_1$ vanishes on $\partial C_2$, $c^{-1}(\hat{\tau}_x^{-1})^\ast r_x\leqslant \hat{r}\leqslant ct^{(2m_x-3)/3j(x)}(\hat{\tau}_x^{-1})^\ast r_x$, and $c^{-1}(\hat{\tau}_x^{-1})^\ast g_{C_1}\leqslant g_{C_2}$ on $C_2$. We have
 \[\lVert v_1\rVert_{L_{-\delta}^2(C_1)}^2\leqslant ct^{(2-2\delta)(2m_x-3)/3j(x)}\lVert v_2\rVert_{L_{-\delta}^2(C_2)}^2,\quad \lVert \mathrm{d}v_2\rVert_{L_{-\delta-1}^2(C_2)}^2\leqslant ct^{(4m_x-6)\delta/3j(x)}\lVert \mathrm{d}v_1\rVert_{L_{-\delta-1}^2(C_1)}^2.\]
 Let $\check{L}_{\delta}^{k,p}(C_2)$ be defined using the Euclidean metric $g_{e(C_2)}$ and the weight function $\hat{r}$. An analogue of Lemma \ref{SobEmbed_lem} (\romannumeral2) implies that \[\lVert v_2\rVert_{\check{L}_{-\delta/2}^2(C_2)}^2\leqslant c\lVert \mathrm{d}v_2\rVert_{\check{L}_{-\delta/2-1}^2(C_2)}^2,\]
 since $v_2$ vanishes on $\partial C_2$. This is equivalent to $\lVert v_2\rVert_{L_{-\delta}^2(C_2)}^2\leqslant c\lVert \mathrm{d}v_2\rVert_{L_{-\delta-1}^2(C_2)}^2$. Therefore \eqref{hardyineqIt_eq1} and \eqref{hardyineqIt_eq} are verified.

 For $x\in P_w$, let $L_{-\delta,x}$ be defined as that for $x\in  I_u$. For all $v\in L_{-2+\delta}^{2,2}$ with $v$ vanishing in $\hat{W}_0$, using the stereographic projection sending $0$ to $\infty$ as above, we have
 \begin{equation}\label{hardyineqPw_eq}
 \lVert v\rVert_{L_{-\delta,x}^2}^2\leqslant ct^{(4m_0+2)/3j(0)}\lVert \nabla_{A_t}v\rVert_{L_{-\delta-1,x}^2}^2.
 \end{equation}

 For $x\in P_s$, let $\hat{\chi}_x(r):=\chi(\lambda_x(t)^2 r)$, where $r=|\zeta_{x,t}|$. In $\hat{B}_{x,t}$, similar to \eqref{HiggsIuEst_eq}, we have \begin{equation}\label{HiggsTEst_eq}
 |[\varphi_t,u]|^2\geqslant c^{-1}t^2|u|^2.
 \end{equation}
 Let $v$ be vanishing on $\partial\hat{B}_{x,t}$, then as before,
 \begin{equation}\label{hardyineqPs_eq}
   \lVert v\rVert_{L_{-\delta}^2(\hat{B}_{x,t})}^2\leqslant ct^{-2}\lVert \nabla_{A_t}u\rVert_{L_{-\delta-1}^2(\hat{B}_{x,t})}^2.
 \end{equation}

 Applying \eqref{hardyineqIu_eq}, \eqref{hardyineqIt_eq} to $v=(1-\chi_x)u$, \eqref{hardyineqPw_eq} to $v=(1-\chi_0)u$, and \eqref{hardyineqPs_eq} to $v=\hat{\chi}_x u$, and using \eqref{HiggsIuEst_eq}, \eqref{HiggsItEst_eq}, \eqref{HiggsTEst_eq},
 \begin{align*}
   \lVert u\rVert_{L_{-\delta}^2}^2&\leqslant c\sum_{x\in I}\lVert u\rVert_{L_{-\delta,x}^2}^2+c\sum_{x\in P_w}\lVert (1-\chi_0)u\rVert_{L_{-\delta,x}^2}^2 + c\sum_{x\in P_s} \lVert\hat{\chi}_x u\rVert_{L_{-\delta}^2(\hat{B}_{x,t})}^2\\& \leqslant c\sum_{x\in I}\lVert \chi_xu\rVert_{L^2(W_x)}^2+c\sum_{x\in I}\lVert (1-\chi_x)u\rVert_{L_{-\delta,x}^2}^2+c\sum_{x\in P_w}\lVert (1-\chi_0)u\rVert_{L_{-\delta,x}^2}+ c\sum_{x\in P_s} \lVert\hat{\chi}_x u\rVert_{L_{-\delta}^2(\hat{B}_{x,t})}^2\\&\leqslant c\Big(\sum_{x\in I_u}t^{-(4m_x+2)/3j(x)}+\sum_{x\in I_t}t^{-4m_x/3j(x)}\Big)\lVert [\varphi_t,u]\rVert_{L^2}^2\\&\quad+c\sum_{x\in I_u} t^{(4m_x+2)/3j(x)}\lVert \nabla_{A_t}(1-\chi_x)u\rVert_{L_{-\delta-1,x}^2}^2+c\sum_{x\in I_t} t^{4m_x/3j(x)}\lVert \nabla_{A_t}(1-\chi_x)u\rVert_{L_{-\delta-1,x}^2}^2\\&\quad+ct^{(4m_0+2)/3j(0)}\sum_{x\in P_w}\lVert \nabla_{A_t}(1-\chi_0)u\rVert_{L_{-\delta-1,x}^2}^2+ct^{-2}\sum_{x\in P_s}\lVert \nabla_{A_t}\hat{\chi}_x u\rVert_{L_{-\delta-1,x}^2}^2\\
   &\leqslant ct^{8/3}\left(\lVert \nabla_{A_t}u\rVert_{L^2}^2+\lVert [\varphi_t,u]\rVert_{L^2}^2\right),
 \end{align*}
 where in the last inequality we used $(4m_x+2)/3j(x)\leqslant 5/3$ for $x\in I_u$, $4m_x/3j(x)\leqslant 8/3$ for $x\in I_t$, and the $L^2$ norm is uniformly stronger than various $L_{-\delta-1}^2$ norms.
 \end{proof}

 \begin{proposition}\label{InvNormBd_prop}
   The norm of $L_t^{-1}:L_{-2+\delta}^2(\mathrm{i}\mathfrak{su}(E))\to L_{-2+\delta}^{2,2}(\mathrm{i}\mathfrak{su}(E))$ is bounded by $ct^{26/3}$.
 \end{proposition}
 \begin{proof}
 Let $v\in C_0^\infty(C\backslash S)$, then $u=L_t^{-1}(v)\in L_{-2+\delta}^{2,2}$ and the boundary terms in the integration by parts formula vanish as above, so \eqref{intbyParts_eq} holds. By the density of $C_0^\infty(C\backslash S)$ in $L_{-2+\delta}^{2}$ and the continuity of $L_t^{-1}, \mathrm{d}_{A_t}, [\Psi_t,\cdot]$, \eqref{intbyParts_eq} holds for all $u\in L_{-2+\delta}^{2,2}$. Using \eqref{intbyParts_eq} and \eqref{PoincareIneq_eq},
   \begin{align*}
    \lVert u\rVert_{L_{-\delta}^2}^2&\leqslant ct^{8/3}\left(\lVert \nabla_{A_t} u\rVert_{L^2}^2+\lVert [\varphi_t,u]\rVert_{L^2}^2\right)\leqslant ct^{8/3}\langle L_t u,u\rangle_{L^2}\\&\leqslant ct^{8/3}\lVert u\rVert_{L_{-\delta}^2}\lVert L_t u\rVert_{L_{-2+\delta}^2},\\
    \lVert u\rVert_{L_{-2+\delta}^2}&\leqslant \lVert u\rVert_{L_{-\delta}^2}\leqslant ct^{8/3}\lVert L_t u\rVert_{L_{-2+\delta}^2},\\
   \lVert \nabla_{A_t}u\rVert_{L_{-2+\delta}^2}&\leqslant \lVert\nabla_{A_t}u\rVert_{L^2}\leqslant ct^{4/3}\lVert L_t u\rVert_{L_{-2+\delta}^2},\\\lVert [\varphi_t,u]\rVert_{L^2}&\leqslant ct^{4/3}\lVert L_t u\rVert_{L_{-2+\delta}^2}.
  \end{align*}
 For $x\in I$, let $\hat{\chi}_x$ be the function supported in $\widetilde{B}_{x,t}$ where it satisfies $\hat{\chi}_x(r)=\chi(\kappa^{-1} t^{1/j(x)} r)$ ($r=|z_x|$ for $x\in I_u$ and $r=|\zeta_{x,t}|$ for $x\in I_t$). Recall that we have defined $\hat{\chi}_x$ for $x\in P_s$ in the proof of the previous proposition, and we define $\hat{\chi}_x$ for $x\in P_w$ similarly. Using Propositions \ref{analysisIu_prop}, \ref{analysisIt_cor}, \ref{analysisPs_prop}, \ref{analysisPw_prop}, \ref{analysisUext_prop} we have
  \begin{align}
    \lVert u\rVert_{L_{-2+\delta}^{2,2}}&\leqslant \sum_{x\in S}\lVert \hat{\chi}_x u\rVert_{L_{-2+\delta}^{2,2}}+\lVert (1-\sum_{x\in S} \hat{\chi}_x)u \rVert_{L^{2,2}}\\&\leqslant c\sum_{x\in S}\lVert L_t\hat{\chi}_x u\rVert_{L_{-2+\delta}^2}+ct^4\lVert L_t(1-\sum_{x\in S}\hat{\chi}_x) u\rVert_{L^2}\notag\\
    &\leqslant ct^4\big(\lVert  L_t u\rVert_{L_{-2+\delta}^2}+t^2\lVert u\rVert_{L_{-\delta}^2}+t\lVert  \nabla_{A_t} u\rVert_{L^2}+t\lVert [\varphi_t,u]\rVert_{L^2}\big)\notag\\
    &\leqslant ct^{26/3}\lVert L_t u\rVert_{L_{-2+\delta}^2}.\label{InverseBound_eq}
  \end{align}
 \end{proof}

 \section{Deforming to a Solution}
 Now we find a solution to Hitchin's equation by deforming the approximate solution $(A_t,\varphi_t)$ via a complex gauge transformation $\exp(\gamma)$ for $\gamma\in\Omega^0(\mathrm{i}\mathfrak{su}(E))$, which satisfies
 \[F_t(\gamma)=F_{A_t^{\exp(\gamma)}}+\left[\mathrm{e}^{-\gamma}\varphi_t\mathrm{e}^{\gamma},\mathrm{e}^{\gamma}\varphi_t^\ast \mathrm{e}^{-\gamma}\right]=:F_{A_t}+[\varphi_t,\varphi_t^\ast]+L_t\gamma +Q_t\gamma=0,\]
 where by abuse of notation we still denote by $L_t$ the linear operator in the previous section composed with the isomorphism $\mathrm{i}\star$. We will analyze the nonlinear term $Q_t$ in a similar way as that in \cite{mazzeo_swoboda_weiss_witt_2016}. Denote
 \begin{align*}
   A_t^{\exp \gamma}&=A_t+(\bar{\partial}_{A_t}-\partial_{A_t})\gamma+R_{A_t}(\gamma),\\ \varphi_t^{\exp \gamma}&=\varphi_t+[\varphi_t,\gamma]+R_{\varphi_t}(\gamma),\\
   R_{A_t}(\gamma)&=\mathrm{e}^{-\gamma}(\bar{\partial}_{A_t}\mathrm{e}^{\gamma})-(\partial_{A_t}\mathrm{e}^{\gamma})\mathrm{e}^{-\gamma}-(\bar{\partial}_{A_t}-\partial_{A_t})\gamma,\\
   R_{\varphi_t}(\gamma)&=\mathrm{e}^{-\gamma}\varphi_t \mathrm{e}^{\gamma}-[\varphi_t,\gamma]-\varphi_t.
 \end{align*}
 Then we have \begin{align}
   Q_t(\gamma)&=\mathrm{d}_{A_t}(R_{A_t}(\gamma))+[R_{\varphi_t}(\gamma),\varphi_t^\ast]+[\varphi_t,R_{\varphi_t}(\gamma)^\ast]\notag\\&\quad+\frac{1}{2}[(\bar{\partial}_{A_t}-\partial_{A_t})\gamma+R_{A_t}(\gamma),(\bar{\partial}_{A_t}-\partial_{A_t})\gamma+R_{A_t}(\gamma)]\notag\\
   &\quad+[[\varphi_t,\gamma]+R_{\varphi_t}(\gamma),([\varphi_t,\gamma]+R_{\varphi_t}(\gamma))^\ast].\label{NonlinTerm_eq}
 \end{align}
 For $u\in \Omega^p(\mathrm{End}\, E)$, define \[\lVert u\rVert_{\widetilde{L}_{\delta}^{k,q}}:=\bigg(\sum_{j=0}^k\lVert \nabla_{A_t}^j u\rVert_{L_{\delta}^q}^q\bigg)^{1/q}.\]
 Clearly $\lVert u\rVert_{\widetilde{L}_{\delta}^{k,q}}\leqslant \lVert u\rVert_{L_{\delta}^{k,q}}$.

 \begin{lemma}
   $\widetilde{L}_{-2+\delta}^{2,2}(\mathrm{End}\,E)$ is an algebra, and $\widetilde{L}_{-2+\delta}^{1,2}(\mathrm{End}\,E)$ is a module over this algebra. Moreover,
   \begin{align}
     \lVert uv\rVert_{\widetilde{L}_{-2+\delta}^{2,2}}&\leqslant ct^5 \lVert u\rVert_{\widetilde{L}_{-2+\delta}^{2,2}}\lVert v\rVert_{\widetilde{L}_{-2+\delta}^{2,2}},\\
     \lVert uv\rVert_{\widetilde{L}_{-2+\delta}^{1,2}}&\leqslant ct^{9/2} \lVert u\rVert_{\widetilde{L}_{-2+\delta}^{2,2}}\lVert v\rVert_{\widetilde{L}_{-2+\delta}^{1,2}},\\
     \lVert u\rVert_{L^\infty}&\leqslant ct^3 \lVert u\rVert_{\widetilde{L}_{-2+\delta}^{2,2}}.
   \end{align}
 \end{lemma}
 \begin{proof}
   Let $\tau_x$ and $\hat{\chi}_x$ be defined as in the previous section, and denote $u_x:=\hat{\chi}_xu$, $u_{\mathrm{ext}}:=(1-\sum_{x\in S}\hat{\chi}_x)u$. For $x\in I_u$, $\tau_x$ maps $\widetilde{B}_{x,t}$ to $B_0(\kappa |c_x|):=\{\,|z|<\kappa |c_x|\,\}$. Denote the Euclidean metric on this disk by $g_e$, then $c^{-1}t^{2/j(x)}(\tau_x^{-1})^\ast g_t\leqslant g_e\leqslant ct^{2/j(x)}(\tau_x^{-1})^\ast g_t$. Define the Sobolev spaces on $B_0(\kappa |c_x|)$ using the weight function $r=|z|$ and the metric $g_e$. For $w$ a $p$-form with values in $\mathrm{End} E$, using the Sobolev embedding $L_{-2+\delta}^{1,2}\hookrightarrow L_{-1+\delta}^q$ on $B_0(\kappa |c_x|)$ for any fixed $q\in[2,\infty)$,
   \begin{align*}
       \lVert w_x\rVert_{L_{-1+\delta}^q}&\leqslant ct^{-2/qj(x)}\lVert (\tau_x^{-1})^{\ast}|w_x|\rVert_{L_{-1+\delta}^q(B_0(\kappa |c_x|))}\\&\leqslant ct^{-2/qj(x)}\big(\lVert (\tau_x^{-1})^\ast |w_x|\rVert_{L_{-2+\delta}^2(B_0(\kappa |c_x|))}+\lVert (\tau_x^{-1})^\ast \mathrm{d}|w_x|\rVert_{L_{-2+\delta}^2(B_0(\kappa |c_x|))}\big)\\
       &\leqslant ct^{-2/qj(x)}\big(t^{1/j(x)}\lVert w_x\rVert_{L_{-2+\delta}^2}+\lVert \nabla_{A_t}w_x\rVert_{L_{-2+\delta}^2}\big)\leqslant ct^{(q-2)/qj(x)}\lVert w_x\rVert_{\widetilde{L}_{-2+\delta}^{1,2}}.
   \end{align*}
   Similarly for $x\in P$, we have $\lVert w_x\rVert_{L_{-1+\delta}^q}\leqslant ct^{(q-2)/q}\lVert w_x\rVert_{\widetilde{L}_{-2+\delta}^{1,2}}$. For $x\in I_t$, let $\check{L}_{\delta}^{k,m}(B_0(\kappa |c_x|))$ be defined using $r^{-1}g_e$ and the weight function $r^{1/2}$, then $\check{L}_{\delta}^{q}$ is equivalent to $L_{\delta/2}$ for functions on $B_0(\kappa |c_x|)$, and we have
   \begin{align*}
     \lVert w_x\rVert_{L_{-1+\delta}^q}&\leqslant ct^{-2/qj(x)}\lVert (\tau_x^{-1})^{\ast}|w_x|\rVert_{\check{L}_{-1+\delta}^q(B_0(\kappa |c_x|))}\leqslant ct^{-2/qj(x)}\lVert (\tau_x^{-1})^{\ast}|w_x|\rVert_{L_{(-1+\delta)/2}^q(B_0(\kappa |c_x|))}
     \\&\leqslant ct^{-2/qj(x)}\big(\lVert (\tau_x^{-1})^\ast |w_x|\rVert_{L_{(-1+\delta)/2-1}^2(B_0(\kappa |c_x|))}+\lVert (\tau_x^{-1})^\ast \mathrm{d}|w_x|\rVert_{L_{(-1+\delta)/2-1}^2(B_0(\kappa |c_x|))}\big)
     \\&\leqslant ct^{-2/qj(x)}\big(\lVert (\tau_x^{-1})^\ast |w_x|\rVert_{\check{L}_{-2+\delta}^2(B_0(\kappa |c_x|))}+\lVert (\tau_x^{-1})^\ast \mathrm{d}|w_x|\rVert_{\check{L}_{-2+\delta}^2(B_0(\kappa |c_x|))}\big)\\
     &\leqslant ct^{-2/qj(x)}\big(t^{1/j(x)}\lVert w_x\rVert_{L_{-2+\delta}^2}+\lVert \nabla_{A_t}w_x\rVert_{L_{-2+\delta}^2}\big)\leqslant ct^{(q-2)/qj(x)}\lVert w_x\rVert_{\widetilde{L}_{-2+\delta}^{1,2}}.
 \end{align*}
 On $U_{\mathrm{ext}}$, $c^{-1}g_C\leqslant g_t\leqslant ctg_C$, and $r_t$ is uniformly bounded, then
 \begin{align}
   \lVert w_{\mathrm{ext}}\rVert_{L^q}&\leqslant ct^{1/q}\lVert w_{\mathrm{ext}}\rVert_{L^q(g_C)}\leqslant ct^{1/q}\left(\lVert w_{\mathrm{ext}}\rVert_{L^2(g_C)}+\lVert \nabla_{A_t}w_{\mathrm{ext}}\rVert_{L^2(g_C)}\right)\notag\\
   &\leqslant ct^{1/q}\left(t^{\max(p-1,0)/2}\lVert w_{\mathrm{ext}}\rVert_{L^2}+t^{p/2}\lVert \nabla_{A_t}w_{\mathrm{ext}}\rVert_{L^2}\right).\label{SobineqExt_eq}
 \end{align}
 Combining the above inequalities,
   \[\lVert w\rVert_{L_{-1+\delta}^q}\leqslant \sum_{x\in S}\lVert w_x\rVert_{L_{-1+\delta}^q}+\lVert w_{\mathrm{ext}}\rVert_{L_{-1+\delta}^q}\leqslant ct^{2+p/2}\lVert w\rVert_{\widetilde{L}_{-2+\delta}^{1,2}}.\]
 For $w_i\in \Omega^{p_i}(\mathrm{End} E)$ ($i=1,2$),
   \begin{align}
     \lVert w_1w_2\rVert_{L_{-2+\delta}^2}&=\lVert r_t^{1-\delta}w_1w_2\rVert_{L^2}\leqslant \lVert r_t^{1-2\delta}w_1w_2\rVert_{L^2}\leqslant \lVert r_t^{1/2-\delta}w_1\rVert_{L^4}\lVert r_t^{1/2-\delta}w_2\rVert_{L^4}\notag\\&=\lVert w_1\rVert_{L_{-1+\delta}^4}\lVert w_2\rVert_{L_{-1+\delta}^4}\leqslant ct^{4+(p_1+p_2)/2}\lVert w_1\rVert_{\widetilde{L}_{-2+\delta}^{1,2}}\lVert w_2\rVert_{\widetilde{L}_{-2+\delta}^{1,2}}.\label{SobMultIneq_eq}
   \end{align}
 For $2<q<2/(1-\delta)$, $x\in I_u$, using $L_{-1+\delta}^q\subset L^q$, $ L^{1,q}\subset L^\infty$ on $B_0(\kappa|c_x|)$,
   \begin{align*}
     \lVert u_x\rVert_{L^\infty}&=\lVert (\tau_x^{-1})^\ast |u_x|\rVert_{L^\infty(B_0(\kappa|c_x|))}\\&\leqslant c\big(\lVert (\tau_x^{-1})^{\ast}|u_x|\rVert_{L_{-1+\delta}^q(B_0(\kappa|c_x|)}+\lVert (\tau_x^{-1})^\ast \mathrm{d}|u_x|\rVert_{L_{-1+\delta}^q(B_0(\kappa|c_x|)}\big)\\
     &\leqslant c\big(t^{2/qj(x)}\lVert u_x\rVert_{L_{-1+\delta}^q}+t^{(2-q)/qj(x)}\lVert\nabla_{A_t}u_x\rVert_{L_{-1+\delta}^q}\big)\\
     &\leqslant c\big(t^{2/qj(x)}t^{(q-2)/qj(x)}\lVert u_x\rVert_{\widetilde{L}_{-2+\delta}^{1,2}}+\lVert \nabla_{A_t}u_x\rVert_{\widetilde{L}_{-2+\delta}^{1,2}}\big)\leqslant ct^{1/j(x)}\lVert u\rVert_{\widetilde{L}_{-2+\delta}^{2,2}}.
   \end{align*}
   Similarly, we have $\lVert u_x\rVert_{L^\infty}\leqslant ct\lVert u_x\rVert_{\widetilde{L}_{-2+\delta}^{2,2}}$ for $x\in P$, and $\lVert u_x\rVert_{L^\infty}\leqslant ct^{1/j(x)}\lVert u_x\rVert_{\widetilde{L}_{-2+\delta}^{2,2}}$ for $x\in I_t$. On $U_{\mathrm{ext}}$, by \eqref{SobineqExt_eq},
   \[  \lVert u_{\mathrm{ext}}\rVert_{L^\infty}\leqslant c\left(\lVert u_{\mathrm{ext}}\rVert_{L^q(g_C)}+\lVert \nabla_{A_t}u_{\mathrm{ext}}\rVert_{L^q(g_C)}\right)\leqslant ct^{1/2}\lVert u_{\mathrm{ext}}\rVert_{\widetilde{L}_{-2+\delta}^{2,2}}.\]
 Therefore,
   \[\lVert u\rVert_{L^{\infty}}\leqslant \sum_{x\in S}\lVert u_x\rVert_{L^{\infty}}+\lVert u_{\mathrm{ext}}\rVert_{L^{\infty}}\leqslant ct^3\lVert u\rVert_{\widetilde{L}_{-2+\delta}^{2,2}}.\]
   Then by \eqref{SobMultIneq_eq},
   \begin{align*}
     \lVert \nabla_{A_t}(uv)\rVert_{L_{-2+\delta}^2}&\leqslant \lVert (\nabla_{A_t}u)v\rVert_{L_{-2+\delta}^2}+\lVert u\nabla_{A_t}v\rVert_{L_{-2+\delta}^2}\\&\leqslant ct^{9/2}\lVert u\rVert_{\widetilde{L}_{-2+\delta}^{2,2}}\lVert v\rVert_{\widetilde{L}_{-2+\delta}^{1,2}}+\lVert u\rVert_{L^\infty}\lVert v\rVert_{\widetilde{L}_{-2+\delta}^{1,2}}
   \leqslant ct^{9/2}\lVert u\rVert_{\widetilde{L}_{-2+\delta}^{2,2}}\lVert v\rVert_{\widetilde{L}_{-2+\delta}^{1,2}}.
   \end{align*}
   Finally we estimate
   \begin{align*}
     \lVert\nabla_{A_t}^2 uv\rVert_{L_{-2+\delta}^2}&\leqslant \lVert u\rVert_{L^\infty}\lVert \nabla_{A_t}^2 v\rVert_{L_{-2+\delta}^2}+ct^5\lVert \nabla_{A_t}u\rVert_{\widetilde{L}_{-2+\delta}^{1,2}}\lVert \nabla_{A_t}v\rVert_{\widetilde{L}_{-2+\delta}^{1,2}}+\lVert v\rVert_{L^\infty}\lVert \nabla_{A_t}^2 u\rVert_{L_{-2+\delta}^2}\\
     &\leqslant ct^5\lVert u\rVert_{\widetilde{L}_{-2+\delta}^{2,2}}\lVert v\rVert_{\widetilde{L}_{-2+\delta}^{2,2}}.
   \end{align*}
 \end{proof}
 \begin{theorem}\label{DefSol_thm}
 There exists $t_N$ such that for any $t\geqslant t_N$, there is a unique $h_t^{\mathrm{app}}$-Hermitian $\gamma_t\in L_{-2+\delta}^{2,2}(\mathrm{i}\mathfrak{su}(E))$ with $\lVert \gamma_t\rVert_{L_{-2+\delta}^{2,2}}\leqslant c_Nt^{-N}$, such that $F_t(\gamma_t)=0$. If $P_w=\varnothing$, then the decay rate of $\gamma_t$ can be improved to $\lVert \gamma_t\rVert_{L_{-2+\delta}^{2,2}}\leqslant c\exp(-c't^\sigma)$ for some constants $c,c'$. Hence $h_t$ defined by $h_t(v,w)=h_t^{\mathrm{app}}(\mathrm{e}^{-\gamma_t}v,\mathrm{e}^{-\gamma_t}w)$ is a harmonic metric associated to $(\bar{\partial}_E,\varphi_t)$ and is compatible with the parabolic structure at each $x\in S$.
 \end{theorem}
 \begin{proof}
 Let $0<\rho\leqslant 1$ and $\gamma_1,\gamma_2\in B_\rho:=\{\gamma\in L_{-2+\delta}^{2,2}(\mathrm{i}\mathfrak{su}(E)):\lVert \gamma\rVert_{L_{-2+\delta}^{2,2}}<\rho\}$. We estimate each of the five terms in the difference $Q_t(\gamma_1)-Q_t(\gamma_2)$, with $Q_t$ given by \eqref{NonlinTerm_eq}. First, \begin{align}
   \lVert R_{A_t}(\gamma_1)-R_{A_t}(\gamma_2))\rVert_{\widetilde{L}_{-2+\delta}^{1,2}}&\leqslant\lVert \left(\mathrm{e}^{-\gamma_1}(\bar{\partial}_{A_t}\mathrm{e}^{\gamma_1})-\mathrm{e}^{-\gamma_2}(\bar{\partial}_{A_t}\mathrm{e}^{\gamma_2})-\bar{\partial}_{A_t}(\gamma_1-\gamma_2)\right)\rVert_{\widetilde{L}_{-2+\delta}^{1,2}}\notag\\
   &\quad +\lVert \left((\partial_{A_t}\mathrm{e}^{\gamma_1})\mathrm{e}^{-\gamma_1}-(\partial_{A_t}\mathrm{e}^{\gamma_2})\mathrm{e}^{-\gamma_2}+\partial_{A_t}(\gamma_1-\gamma_2)\right)\rVert_{\widetilde{L}_{-2+\delta}^{1,2}}.\label{NonlinearAaDiff_eq}
 \end{align}
 The first term is bounded by
 \begin{align*}
   \lVert (\mathrm{e}^{-\gamma_1}-\mathrm{e}^{-\gamma_2})\bar{\partial}_{A_t} \mathrm{e}^{\gamma_1}\rVert_{\widetilde{L}_{-2+\delta}^{1,2}}+\lVert \mathrm{e}^{-\gamma_2}(\bar{\partial}_{A_t}(\mathrm{e}^{\gamma_1}-\mathrm{e}^{\gamma_2}))-\bar{\partial}_{A_t}(\gamma_1-\gamma_2) \rVert_{\widetilde{L}_{-2+\delta}^{1,2}}:=\uppercase\expandafter{\romannumeral1}+\uppercase\expandafter{\romannumeral2}.
 \end{align*}
 Write $\mathrm{e}^\gamma=1+\gamma+S(\gamma)$, then
 \begin{align*}
   \uppercase\expandafter{\romannumeral1}&\leqslant c t^{9/2}\lVert \mathrm{e}^{-\gamma_1}-\mathrm{e}^{-\gamma_2}\rVert_{\widetilde{L}_{-2+\delta}^{2,2}}\lVert \bar{\partial}_{A_t}\mathrm{e}^{\gamma_1}\rVert_{\widetilde{L}_{-2+\delta}^{1,2}}\leqslant c t^{9/2}\lVert \mathrm{e}^{-\gamma_1}- \mathrm{e}^{-\gamma_2}\rVert_{\widetilde{L}_{-2+\delta}^{2,2}}\lVert \gamma_1+S(\gamma_1)\rVert_{\widetilde{L}_{-2+\delta}^{2,2}}
   \\&\leqslant c t^{9/2}\lVert \gamma_1-\gamma_2\rVert_{\widetilde{L}_{-2+\delta}^{2,2}}\bigg(\sum_{j=0}^\infty \frac{(\rho ct^5)^j}{j!}\bigg)\rho\bigg(\sum_{j=0}^\infty \frac{(\rho ct^5)^j}{(j+1)!}\bigg)\\
   &\leqslant ct^{9/2}\rho\lVert \gamma_1-\gamma_2\rVert_{\widetilde{L}_{-2+\delta}^{2,2}},
 \end{align*}
 for $\rho<t^{-5}$. Similarly,
 \begin{align*}
   \uppercase\expandafter{\romannumeral2}&=\lVert (1-\gamma_2+S(-\gamma_2))\left(\bar{\partial}_{A_t}(\gamma_1-\gamma_2+S(\gamma_1)-S(\gamma_2))\right)-\bar{\partial}_{A_t}(\gamma_1-\gamma_2)\rVert_{\widetilde{L}_{-2+\delta}^{1,2}}\\
   &\leqslant \lVert \bar{\partial}_{A_t}(S(\gamma_1)-S(\gamma_2))\rVert_{\widetilde{L}_{-2+\delta}^{1,2}}+\lVert (-\gamma_2+S(-\gamma_2))\bar{\partial}_{A_t}(\gamma_1-\gamma_2+S(\gamma_1)-S(\gamma_2))\rVert_{\widetilde{L}_{-2+\delta}^{1,2}}\\
   &\leqslant \lVert S(\gamma_1)-S(\gamma_2)\rVert_{\widetilde{L}_{-2+\delta}^{2,2}}+  ct^{9/2}\lVert -\gamma_2+S(-\gamma_2)\rVert_{\widetilde{L}_{-2+\delta}^{2,2}}\lVert \gamma_1-\gamma_2+S(\gamma_1)-S(\gamma_2)\rVert_{\widetilde{L}_{-2+\delta}^{2,2}}\\
   &\leqslant ct^{9/2}\rho \lVert \gamma_1-\gamma_2\rVert_{\widetilde{L}_{-2+\delta}^{2,2}}.
 \end{align*}
 The estimate for the other term in \eqref{NonlinearAaDiff_eq} is analogous, and we obtain \begin{align}
     \lVert R_{A_t}(\gamma_1)-R_{A_t}(\gamma_2))\rVert_{\widetilde{L}_{-2+\delta}^{1,2}}&\leqslant ct^{9/2}\rho\lVert \gamma_1-\gamma_2\rVert_{\widetilde{L}_{-2+\delta}^{2,2}}\leqslant ct^{9/2}\rho\lVert \gamma_1-\gamma_2\rVert_{L_{-2+\delta}^{2,2}},\label{RAtIneq_eq}\\
     \lVert \mathrm{d}_{A_t}(R_{A_t}(\gamma_1)-R_{A_t}(\gamma_2))\rVert_{L_{-2+\delta}^2}&\leqslant ct^{9/2}\rho\lVert \gamma_1-\gamma_2\rVert_{L_{-2+\delta}^{2,2}}.\notag
 \end{align}
 Let $\beta_1,\ldots,\beta_n\in L_{-2+\delta}^{2,2}$, where $n\geqslant 2$. By the Jacobi identity, $\left|[\varphi_t^\ast,[\beta_1,\cdots,[\beta_n,\varphi_t]\cdots]]\right|$ is bounded by
 \begin{align*}
   & \left|[[\varphi_t^\ast,\beta_1],[\beta_2,\cdots,[\varphi_t,\beta_n]\cdots]]\right|+\left|[\beta_1,[[\varphi_t^\ast,\beta_2],[\beta_3,\cdots[\varphi_t,\beta_n]\cdots]]]\right|\\
   +\cdots&+\left|[\beta_1,\cdots[[\varphi_t^\ast,\beta_{n-1}],[\varphi_t,\beta_n]]\cdots]\right|+\left|[\beta_1,\cdots[\beta_{n-1},[\varphi_t^\ast,[\varphi_t,\beta_n]]]\cdots]\right|\\
   &\leqslant 2^n\bigg(\sum_{j=1}^{n-1}\big(|[\varphi_t^\ast,\beta_j]||[\varphi_t,\beta_n]|\prod_{m\neq j,n}\lVert \beta_m\rVert_{L^\infty}\big)+|[\varphi_t^\ast,[\varphi_t,\beta_n]]|\prod_{m=1}^{n-1}\lVert \beta_m\rVert_{L^\infty}\bigg).
 \end{align*}
 Integrate both sides, then $\left\lVert[\varphi_t^\ast,[\beta_1,\cdots,[\beta_n,\varphi_t]\cdots]]\right\rVert_{L_{-2+\delta}^2}$ is bounded by
 \begin{align*}
     & 2^n\bigg(\sum_{j=1}^{n-1}ct^5\lVert [\varphi_t^\ast,\beta_j]\rVert_{\widetilde{L}_{-2+\delta}^{1,2}}\lVert [\varphi_t,\beta_n]\rVert_{\widetilde{L}_{-2+\delta}^{1,2}}(ct^3)^{n-2}\prod_{m\neq j,n}\lVert \beta_m\rVert_{\widetilde{L}_{-2+\delta}^{2,2}} \\
     &\qquad+ \lVert[\varphi_t^\ast,[\varphi_t,\beta_n]]\rVert_{L_{-2+\delta}^2}(ct^3)^{n-1}\prod_{m=1}^{n-1}\lVert \beta_m\rVert_{\widetilde{L}_{-2+\delta}^{2,2}} \bigg)\\
     &\leqslant (ct^{9/2})^{n-1}t^{1/2}\prod_{m=1}^n\lVert \beta_m\rVert_{L_{-2+\delta}^{2,2}}.
 \end{align*}
 Then for $\gamma_1,\gamma_2\in B_\rho$,
 \begin{align*}
 \left\lVert [\varphi_t^\ast,(\mathrm{ad}_{\gamma_1})^n\varphi_t-(\mathrm{ad}_{\gamma_2})^n\varphi_t]\right\rVert_{L_{-2+\delta}^2}&\leqslant \left\lVert[\varphi_t^\ast,[\gamma_1-\gamma_2,[\gamma_1,\cdots,[\gamma_1,\varphi_t]\cdots]]]\right\rVert_{L_{-2+\delta}^2}\\&\quad+\left\lVert[\varphi_t^\ast,[\gamma_2,[\gamma_1-\gamma_2,[\gamma_1,\cdots,[\gamma_1,\varphi_t]\cdots]]]]\right\rVert_{L_{-2+\delta}^2}\\&\quad+\cdots+\left\lVert[\varphi_t^\ast,[\gamma_2,[\gamma_2,\cdots,[\gamma_1-\gamma_2,\varphi_t]\cdots]]]\right\rVert_{L_{-2+\delta}^2}\\
   &\leqslant  n(ct^{9/2}\rho)^{n-1}t^{1/2}\lVert \gamma_1-\gamma_2\rVert_{L_{-2+\delta}^{2,2}}.
 \end{align*}
 Then we compute that
 \begin{align*}
   \lVert [R_{\varphi_t}(\gamma_1),\varphi_t^\ast]-[R_{\varphi_t}(\gamma_2),\varphi_t^\ast]\rVert_{L_{-2+\delta}^2}&=\lVert[\varphi_t^\ast,\sum_{j=2}^\infty \frac{1}{j!}(\mathrm{ad}_{-\gamma_1})^j\varphi_t-(\mathrm{ad}_{-\gamma_2})^j\varphi_t]\rVert_{L_{-2+\delta}^2}\\
   &\leqslant \sum_{j=2}^\infty \frac{1}{j!}j(c t^{9/2}\rho)^{j-1}t^{1/2}\lVert \gamma_1-\gamma_2\rVert_{L_{-2+\delta}^{2,2}}\\
   &\leqslant ct^{5}\rho\lVert \gamma_1-\gamma_2\rVert_{L_{-2+\delta}^{2,2}}.
 \end{align*}
 Similarly,
 \[\lVert [\varphi_t,R_{\varphi_t}(\gamma_1)^\ast]-[\varphi_t,R_{\varphi_t}(\gamma_2)^\ast]\rVert_{L_{-2+\delta}^2}\leqslant ct^{5}\rho\lVert \gamma_1-\gamma_2\rVert_{L_{-2+\delta}^{2,2}}.\]
 The $L_{-2+\delta}^2$ norm of the fourth term in $Q_t(\gamma_1)-Q_t(\gamma_2)$ is bounded by
 \begin{align*}
 &\,\quad \lVert[(\bar{\partial}_{A_t}-\partial_{A_t})(\gamma_1-\gamma_2)+R_{A_t}(\gamma_1)-R_{A_t}(\gamma_2),(\bar{\partial}_{A_t}-\partial_{A_t})\gamma_1+R_{A_t}(\gamma_1)]\rVert_{L_{-2+\delta}^2}\\
   &+\lVert[(\bar{\partial}_{A_t}-\partial_{A_t})\gamma_2+R_{A_t}(\gamma_2),(\bar{\partial}_{A_t}-\partial_{A_t})(\gamma_1-\gamma_2)+R_{A_t}(\gamma_1)-R_{A_t}(\gamma_2)]\rVert_{L_{-2+\delta}^2}:=\uppercase\expandafter{\romannumeral1}'+\uppercase\expandafter{\romannumeral2}'.
 \end{align*}
 By \eqref{SobMultIneq_eq}, \eqref{RAtIneq_eq},
 \begin{align*}
   \uppercase\expandafter{\romannumeral1}'&\leqslant ct^5\lVert (\bar{\partial}_{A_t}-\partial_{A_t})(\gamma_1-\gamma_2)+R_{A_t}(\gamma_1)-R_{A_t}(\gamma_2)\rVert_{\widetilde{L}_{-2+\delta}^{1,2}}\lVert (\bar{\partial}_{A_t}-\partial_{A_t})\gamma_1+R_{A_t}(\gamma_1)\rVert_{\widetilde{L}_{-2+\delta}^{1,2}}\\
   &\leqslant ct^5\lVert \gamma_1-\gamma_2\rVert_{L_{-2+\delta}^{2,2}}(1+ct^{9/2}\rho)^2\rho\leqslant ct^5\rho\lVert \gamma_1-\gamma_2\rVert_{L_{-2+\delta}^{2,2}}.
 \end{align*}
 Similarly, $\uppercase\expandafter{\romannumeral2}'$ is also bounded by $ct^5\rho\lVert \gamma_1-\gamma_2\rVert_{L_{-2+\delta}^{2,2}}$. Consider $\beta_1,\ldots,\beta_n\in L_{-2+\delta}^{2,2}$, then
 \begin{align*}
   \lVert [\beta_1,[\beta_2,\cdots,[\beta_n,\varphi_t]\cdots]]\rVert_{\widetilde{L}_{-2+\delta}^{1,2}}&\leqslant (  ct^{9/2})^{n-1}\prod_{m=1}^{n-1}\lVert\beta_m\rVert_{L_{-2+\delta}^{2,2}}\lVert [\varphi_t,\beta_n]\rVert_{\widetilde{L}_{-2+\delta}^{1,2}}\\
   &\leqslant (  ct^{9/2})^{n-1}\prod_{m=1}^{n}\lVert\beta_m\rVert_{L_{-2+\delta}^{2,2}}.
 \end{align*}
 This implies that
 \[
   \lVert R_{\varphi_t}(\gamma_1)-R_{\varphi_t}(\gamma_2)\rVert_{\widetilde{L}_{-2+\delta}^{1,2}}\leqslant \sum_{j=2}^\infty \frac{1}{j!}j(  ct^{9/2}\rho)^{j-1}\lVert \gamma_1-\gamma_2\rVert_{L_{-2+\delta}^{2,2}}\leqslant ct^{9/2}\rho\lVert \gamma_1-\gamma_2\rVert_{L_{-2+\delta}^{2,2}}.
 \]
 The $L_{-2+\delta}^2$ norm of the last term in $Q_t(\gamma_1)-Q_t(\gamma_2)$ is bounded by
 \begin{align*}
 &\,\quad\lVert [[\varphi_t,\gamma_1-\gamma_2]+R_{\varphi_t}(\gamma_1)-R_{\varphi_t}(\gamma_2),([\varphi_t,\gamma_1]+R_{\varphi_t}(\gamma_1))^\ast]\rVert_{L_{-2+\delta}^2}\\&+ \lVert [[\varphi_t,\gamma_2]+R_{\varphi_t}(\gamma_2),([\varphi_t,\gamma_1-\gamma_2]+R_{\varphi_t}(\gamma_1)-R_{\varphi_t}(\gamma_2))^\ast]\rVert_{L_{-2+\delta}^2}:=\uppercase\expandafter{\romannumeral1}''+\uppercase\expandafter{\romannumeral2}'',
 \end{align*}
 where $\uppercase\expandafter{\romannumeral1}''$ can be estimated as
 \begin{align*}
   \uppercase\expandafter{\romannumeral1}''&\leqslant ct^5\lVert[\varphi_t,\gamma_1-\gamma_2]+R_{\varphi_t}(\gamma_1)-R_{\varphi_t}(\gamma_2) \rVert_{\widetilde{L}_{-2+\delta}^{1,2}}\lVert [\varphi_t,\gamma_1]+R_{\varphi_t}(\gamma_1) \rVert_{\widetilde{L}_{-2+\delta}^{1,2}}\\
   &\leqslant ct^5\lVert \gamma_1-\gamma_2\rVert_{L_{-2+\delta}^{2,2}}(1+ct^{9/2}\rho)^2\rho\leqslant ct^5\rho\lVert \gamma_1-\gamma_2\rVert_{L_{-2+\delta}^{2,2}},
 \end{align*}
 similar for $\uppercase\expandafter{\romannumeral2}''$. Therefore
 \begin{equation}\label{NonlinearEst_eq}
   \lVert Q_t(\gamma_1)-Q_t(\gamma_2)\rVert_{L_{-2+\delta}^2}\leqslant ct^{5}\rho\lVert \gamma_1-\gamma_2\rVert_{L_{-2+\delta}^{2,2}}.
 \end{equation}
 Now $F_t(\gamma)=0$ is equivalent to \[\gamma=-L_t^{-1}(F_{A_t}+[\varphi_t,\varphi_t^\ast]+Q_t(\gamma)):=T_t(\gamma).\]
 By \eqref{InverseBound_eq}, \eqref{NonlinearEst_eq}, $\lVert T_t(\gamma_1)-T_t(\gamma_2)\rVert_{L_{-2+\delta}^{2,2}}\leqslant ct^{41/3}\rho \lVert \gamma_1-\gamma_2\rVert_{L_{-2+\delta}^{2,2}}$. Furthermore, by \eqref{ErrorEst_eq} (which is still valid when the norm is taken with respect to $g_t$),
 \[\lVert T_t(\gamma)\rVert_{L_{-2+\delta}^{2,2}}\leqslant c_Nt^{26/3-N}+ct^{41/3}\rho^2,\]
 for $\gamma\in B_\rho$. Let $\rho_t=t^{-(5+N)/2}$, where $\epsilon<c'/2$ is small. Then for $N\geqslant 14$ and $t$ large, $\rho_t<t^{-5}$, $T_t$ maps $B_{\rho_t}$ to $B_{\rho_t}$, and is a contraction mapping, hence there is a unique fixed point $\gamma_t$ in $B_{\rho_t}$. In other words, the metric $h_t$ is harmonic, which is also compatible with the parabolic structure as $h_t^{\mathrm{app}}$ is and $\gamma_t\in L_{-2+\delta}^{2,2}$.

 If $P_w=\varnothing$, then the decay of $\gamma_t$ can be improved to be exponential. By \eqref{ErrorEst_eq1},
 \[\lVert T_t(\gamma)\rVert_{L_{-2+\delta}^{2,2}}\leqslant ct^{8}\mathrm{e}^{-c' t^\sigma}+ct^{13}\rho^2,\]
 for $\gamma\in B_\rho$. Let $\rho_t=\mathrm{e}^{(-c'/2+\epsilon)t^\sigma}$, where $\epsilon<c'/2$ is small. Then as above, for $t$ large, there is a unique fixed point $\gamma_t$ of $T_t$ in $B_{\rho_t}$.
 \end{proof}

\section{Asymptotic Metric}
In this section, we study the asymptotic geometry of $\mathcal{M}$ under the framework of \cite{fredrickson_2019, fredrickson2022asymptotic}.

Elements in $\mathcal{M}$ are gauge equivalence classes of pairs $(\bar{\partial}_E,\varphi)$, with infinitesimal deformations denoted by $(\dot{\eta},\dot{\varphi})$. Linearize $\bar{\partial}_E\varphi=0$ and we get \begin{equation}\label{LinHol_eq}
  \bar{\partial}_E\dot{\varphi}+[\dot{\eta},\varphi]=0.
\end{equation}
For $[(\bar{\partial}_E,\varphi)]\in \mathcal{M}$, let $h$ be the unique harmonic metric (if exists) which induces $h_{\det E}$, and denote its infinitesimal deformation by $\dot{h}$. The infinitesimal gauge transformation $\dot{\gamma}$ acts as
\begin{equation}\label{InfGaugeTrans_eq}
  (\dot{\eta},\dot{\varphi},\dot{h})\mapsto (\dot{\eta}+\bar{\partial}_E\dot{\gamma},\dot{\varphi}+[\varphi,\dot{\gamma}],\dot{h}+\dot{\gamma}).
\end{equation}
 The linearized Hitchin equations and the Coulomb gauge condition can be unified into a single equation \cite[Prop.~2.2]{fredrickson_2019}
\begin{equation}\label{LinHitCou_eq}
  \partial_E^h\bar{\partial}_E\dot{h}-\partial_E^h\dot{\eta}-[\varphi^{\ast_h},\dot{\varphi}+[\dot{h},\varphi]]=0.
\end{equation}
Recall that $2\partial_E^h\bar{\partial}_E=F_{A_h}+\mathrm{i}\star \Delta_{A_h}=-[\varphi,\varphi^{\ast_h}]+\mathrm{i}\star \Delta_{A_h}$, the equation can be written as \[-\frac{1}{2}[[\varphi,\varphi^{\ast_h}],\dot{h}]+\frac{1}{2}\mathrm{i}\star \Delta_{A_h}\dot{h}-[\varphi^{\ast_h},[\dot{h},\varphi]]=\partial_{E}^h\dot{\eta}+[\varphi^{\ast_h},\dot{\varphi}].\]
Composing with $-2\mathrm{i}\star$ and simplifying, we obtain
\[L_{\varphi}\dot{h}:=\Delta_{A_h}\dot{h}-\mathrm{i}\star([\varphi^{\ast_h},[\varphi,\dot{h}]] -[\varphi,[\varphi^{\ast_h},\dot{h}]])=-2\mathrm{i}\star(\partial_E^h\dot{\eta}+[\varphi^{\ast_h},\dot{\varphi}]).\]
Note that $L_{\varphi}$ is exactly $L_t$ in \eqref{LinOp_eq} with $(\bar{\partial}_E,\varphi_t,h_{t}^{\mathrm{app}})$ replaced by $(\bar{\partial}_E,\varphi,h)$. Now we consider the solvability of \eqref{LinHitCou_eq}. Let $(\bar{\partial}_E,\varphi_t)$ and $L_{\delta}^{k,p}$ (depending on $g_t$, $\varphi_t$, and $h_t^{\mathrm{app}}$) be defined as before, $\nabla_t$ be the Chern connection determined by $\bar{\partial}_E$ and $h_t$.
\begin{lemma}\label{InvLin_lem}
Fix a pair $[(\bar{\partial}_E,\varphi_t)]\in\mathcal{M}$, let $t$ be sufficiently large so that the harmonic metric $h_t$ exists. Then $L_{\varphi_t}: L_{-2+\delta}^{2,2}(\mathfrak{sl}(E))\to L_{-2+\delta}^2(\mathfrak{sl}(E))$ is an isomorphism.
\end{lemma}
\begin{proof}
That $L_t$ is Fredholm of index zero follows from Lemma \ref{GlobIso_lem}, and the fact that $L_t$ preserves the decomposition $\mathfrak{sl}(E)=\mathfrak{su}(E)\oplus i\mathfrak{su}(E)$. Since $\gamma_t\in L_{-2+\delta}^{2,2}$, $L_{\varphi_t}$ differs from $L_t$ by a compact perturbation, it is also Fredholm of index zero. If $u\in \ker L_{\varphi_t}$, then as before the integration by parts yields $\nabla_t u=[\varphi_t,u]=0$. By the proof of \cite[Prop. 7.4]{fredrickson2022asymptotic}, we have $u=0$. Therefore $L_{\varphi_t}$ is an isomorphism.
\end{proof}
\begin{lemma}\label{InfsmRep_lem}
  Let $[(\bar{\partial}_E,\varphi)]\in \mathcal{M}$ and $(\dot{\eta}_1,\dot{\varphi}_1)$ be an infinitesimal variation. There exists $(\dot{\eta}_2,\dot{\varphi}_2)$ infinitesimally gauge equivalent to $(\dot{\eta}_1,\dot{\varphi}_1)$ such that the followings hold.
\begin{enumerate}[label=(\roman*)]
  \item In a neighborhood of $x\in P_s$ and a holomorphic frame where \[\bar{\partial}_E=\bar{\partial},\quad \varphi=\begin{pmatrix}
    0&1\\\zeta_x^{-1}&0
  \end{pmatrix}\,\mathrm{d}\zeta_x,\]
  we have
  \[\dot{\eta}_2=0,\quad\dot{\varphi}_2=\begin{pmatrix}
    0&0\\\frac{\dot{P}}{\zeta_x}&0
  \end{pmatrix}\,\mathrm{d}\zeta_x,\quad\dot{P}\text{ holomorphic}.\]
  \item In a neighborhood of $x\in I_u\cup P_w$ and a holomorphic frame where \[\bar{\partial}_E=\bar{\partial},\quad \varphi=fz_x^{-m_x} \sigma_3\,\mathrm{d}z_x\] with $f$ being holomorphic and nowhere vanishing, we have $\dot{\eta}_2=0, \dot{\varphi}_2=\dot{P}\sigma_3\,\mathrm{d}z$, $\dot{P}$ is holomorphic.
  \item In a neighborhood of $x\in I_t$ and a holomorphic frame where \[\bar{\partial}_E=\bar{\partial},\quad \varphi=\begin{pmatrix}
    0&\frac{1}{\zeta_{x}^{m_x-1}}\\ \frac{1}{\zeta_{x}^{m_x}}&0
  \end{pmatrix}\,\mathrm{d}\zeta_{x},\] we have \[\dot{\eta}_2=0,\quad \dot{\varphi}_2=\begin{pmatrix}
    0&\dot{P}\\\frac{\dot{P}}{\zeta_x}&0
  \end{pmatrix}\,\mathrm{d}\zeta_x,\quad \dot{P}\text{ is holomorphic}.\]
\end{enumerate}
\end{lemma}
\begin{proof}
  (\romannumeral1) This is just \cite[Lem.~7.3]{fredrickson2022asymptotic}.\\
  (\romannumeral2) By \eqref{InfGaugeTrans_eq} and the Poincaré lemma we may assume that $\dot{\eta}_1=0$ in the disk. The leading term of $\dot{\varphi}_1$ preserves the flag at $x$, then we can write $\dot{\varphi}_1$ as \[\dot{\varphi}_1=\frac{\mathrm{d}z_x}{z_x^{m_x}}\begin{pmatrix}
    \dot{P}_1&\dot{P_2}z_x\\ \dot{P}_3&-\dot{P}_1
  \end{pmatrix}.\]
  where $\dot{P}_1,\dot{P}_2,\dot{P}_3$ are holomorphic by \eqref{LinHol_eq}. Then $\dot{\varphi}_2=\dot{\varphi}_1+[\varphi,\dot{\gamma}]=\dot{P}_1z_x^{-m_x}\sigma_3\,\mathrm{d}z$ if we choose \[\dot{\gamma}=\begin{pmatrix}
    0&-\frac{\dot{P}_2 z_x}{2f}\\\frac{\dot{P}_3}{2f}&0
  \end{pmatrix},\]
  which satisfies $\bar{\partial}\dot{\gamma}=0$ and preserves the flag. Higgs fields in $\mathcal{M}$ have fixed polar parts, implying that $\dot{P}_1z_x^{-m_x}=\dot{P}$ for some holomorphic $\dot{P}$.\\
  (\romannumeral3) As above, assume that $\dot{\eta}_1=0$. The leading term of $\dot{\varphi}_1$ preserves the flag at $x$ and the leading term of the Higgs field is nilpotent, then we can write $\dot{\varphi}_1$ as \[\dot{\varphi}_1=\frac{\mathrm{d}\zeta_x}{\zeta_x^{m_x-1}}\begin{pmatrix}
    \dot{P}_1&\dot{P_2}\\ \frac{\dot{P}_3}{\zeta_x}&-\dot{P}_1
  \end{pmatrix}.\]
  where $\dot{P}_1,\dot{P}_2,\dot{P}_3$ are holomorphic by \eqref{LinHol_eq}. We choose \[\dot{\gamma}=\begin{pmatrix}
    \frac{\dot{P}_2-\dot{P}_3}{4}&\zeta_x\dot{P}_1\\0&-\frac{\dot{P}_2-\dot{P}_3}{4}
  \end{pmatrix},\]
  which satisfies $\bar{\partial}\dot{\gamma}=0$ and preserves the flag. Then \[\dot{\varphi}_2=\dot{\varphi}_1+[\varphi,\dot{\gamma}]=\begin{pmatrix}
    0&\frac{\dot{P}_2+\dot{P}_3}{2\zeta_x^{m_x-1}}\\ \frac{\dot{P}_2+\dot{P}_3}{2\zeta_x^{m_x}}&0
  \end{pmatrix}\,\mathrm{d}\zeta_x.\] Higgs fields in $\mathcal{M}$ have fixed polar parts, implying that $(\dot{P}_2+\dot{P}_3)/(2\zeta_x^{m_x-1})=\dot{P}$ for $\dot{P}$ holomorphic.
\end{proof}

Suppose $[(\dot{\eta},\dot{\varphi})]\in T_{[(\bar{\partial}_E,\varphi_t)]}\mathcal{M}$, by the above lemma we may assume that $\dot{\eta}=0$ and $\dot{\varphi}$ has the above forms near every $x\in S$. Then $\partial_E^{h_t}\dot{\eta}+[\varphi_t^{\ast_{h_t}},\dot{\varphi}]\in L_{-2+\delta}^2$, as $\partial_E^{h_t}\dot{\eta}=0$ near $x\in S$, $[\varphi_t^{\ast_{h_t^{\mathrm{app}}}},\dot{\varphi}]\in L_{-2+\delta}^2$ for $\delta$ small and $\gamma_t\in L_{-2+\delta}^{2,2}$. By Lemma \ref{InvLin_lem}, a unique $\dot{h}\in L_{-2+\delta}^{2,2}$ solving \eqref{LinHitCou_eq} can be found. Now the hyperkähler metric is given by \cite[(7.5)]{fredrickson2022asymptotic}
\begin{equation}\label{HKahlerMetric_eq}
\lVert [(\dot{\eta},\dot{\varphi})]\rVert_{g_{L^2}}^2=2\int_C |\dot{\eta}-\bar{\partial}_E\dot{h}|_{h_t}^2+|\dot{\varphi}+[\dot{h},\varphi_t]|_{h_t}^2=2\int_C \langle \dot{\eta}-\bar{\partial}_E\dot{h},\dot{\eta}\rangle_{h_t}+\langle \dot{\varphi}+[\dot{h},\varphi_t],\dot{\varphi}\rangle_{h_t},
\end{equation}
where the second equality follows from \eqref{LinHitCou_eq} and the integration by parts which is justified as in the proof of Lemma \ref{GlobIso_lem} and \cite[p.~43]{fredrickson2022asymptotic}.

For the pair $[(\bar{\partial}_E,\varphi_t)]\in\mathcal{M}$, one can associate metrics $h_t^{\mathrm{app}}$ and $h_t^{\mathrm{dh}}$, with infinitesimal deformations denoted by $\dot{h}_{\mathrm{app}}$ and $\dot{h}_{\mathrm{dh}}$. The metric $g_{\mathrm{app}}$ on $\mathcal{M}$ is defined by
\begin{equation}\label{ApproxMetric_eq}
\lVert [(\dot{\eta},\dot{\varphi})]\rVert_{g_{\mathrm{app}}}^2=2\int_C \langle \dot{\eta}-\bar{\partial}_E\dot{h}_{\mathrm{app}},\dot{\eta}\rangle_{h_t^{\mathrm{app}}}+\langle \dot{\varphi}+[\dot{h}_{\mathrm{app}},\varphi_t],\dot{\varphi}\rangle_{h_t^{\mathrm{app}}},
\end{equation}
where $\dot{h}_{\mathrm{app}}$ satisfies
\begin{equation}\label{AppLinHitCol_eq}
  \partial_E^{h_t^{\mathrm{app}}}\bar{\partial}_E\dot{h}_{\mathrm{app}}-\partial_E^{h_t^{\mathrm{app}}}\dot{\eta}-[\varphi_t^{\ast_{h_t^{\mathrm{app}}}},\dot{\varphi}+[\dot{h}_{\mathrm{app}},\varphi_t]]=0.
\end{equation}
Similarly, the metric $g_{\mathrm{dh}}$ is
\begin{equation}\label{InftyMetric_eq}
  \lVert[(\dot{\eta},\dot{\varphi})]\rVert_{g_{\mathrm{dh}}}^2=2\int_C \langle \dot{\eta}-\bar{\partial}_E\dot{h}_{\mathrm{dh}},\dot{\eta}\rangle_{h_t^{\mathrm{dh}}}+\langle \dot{\varphi}+[\dot{h}_{\mathrm{dh}},\varphi_t],\dot{\varphi}\rangle_{h_t^{\mathrm{dh}}},
\end{equation}
where $\dot{h}_{\mathrm{dh}}$ satisfies the decoupled equations
\begin{equation}\label{DecoupInftyDef_eq}
  \partial_E^{h_t^{\mathrm{dh}}}\bar{\partial}_E\dot{h}_{\mathrm{dh}}-\partial_{E}^{h_t^{\mathrm{dh}}}\dot{\eta}=0,\quad [\varphi_t^{\ast_{h_t^{\mathrm{dh}}}},\dot{\varphi}+[\dot{h}_{\mathrm{dh}},\varphi_t]]=0.
\end{equation}
Recall that for each $[(\bar{\partial}_E,\varphi_t)]\in \mathcal{M}$, we can view the associated spectral data as a pair $(\bar{\partial}_{L,t},\tau_t)$, where $L$ is a complex line bundle of degree $-|S|+N-2$ over a surface $S_t$ with genus $N-3$, $\tau_t$ is the tautological eigenvalue of $\pi^\ast \varphi_t$ which embeds $S_t$ in $\pi^\ast K(D)$. The tangent map induces an isomorphism
\[T_{[(\bar{\partial}_E,\varphi_t)]}\mathcal{M}\to H^0(S_t,K_{S_t})^\ast_{\mathrm{odd}}\oplus H^0(S_t,K_{S_t})_{\mathrm{odd}},\quad [(\dot{\eta},\dot{\varphi})]\mapsto (\dot{\xi},\dot{\tau}),\] where $H^0(S_t,K_{S_t})_{\mathrm{odd}}$ consists of elements of $H^0(S_t,K_{S_t})$ which are odd under the exchange of two sheets of $\pi:S_t\to C$, and $\dot{\xi},\dot{\tau}$ are the infinitesimal deformations of $\bar{\partial}_{L,t}$ and $\tau$. The metric $h_t^{\mathrm{dh}}$ is the orthogonal pushforward of the harmonic metric $h_{\mathcal{L}_t}$ on the spectral line bundle $\mathcal{L}_t=(L,\bar{\partial}_{L,t})\to S_t$, and the infinitesimal deformation $\dot{h}_L$ of $h_{\mathcal{L}_t}$ is determined by
\[\partial_L^{h_{\mathcal{L}_t}}\bar{\partial}_{L,t}\dot{h}_L-\partial_L^{h_{\mathcal{L}_t}}\dot{\xi}=0.\]
\begin{lemma}[{\cite[Prop.~2.14]{fredrickson_2019}}]The semiflat metric is characterized by three properties:
\begin{enumerate}[label=(\arabic*)]
  \item On horizontal deformations, the semiflat metric is $\int_{S_t}2|\dot{\tau}|^2$.
  \item On vertical deformations, the semiflat metric is $\int_{S_t} 2|\dot{\xi}-\bar{\partial}_L\dot{h}_L|^2$.
  \item Horizontal and vertical deformations are orthogonal.
\end{enumerate}
\end{lemma}
\begin{lemma}[{\cite[Cor.~7.9]{fredrickson2022asymptotic}}]
  The semiflat metric $g_{\mathrm{sf}}$ on $\mathcal{M}$ coincides with $g_{\mathrm{dh}}$.
\end{lemma}
\subsection{Comparing $g_{L^2}$ and $g_{\mathrm{app}}$}
The difference $\frac{1}{2}|\lVert [(\dot{\eta},\dot{\varphi})]\rVert_{g_{\mathrm{app}}}^2-\lVert [(\dot{\eta},\dot{\varphi})]\rVert_{g_{L^2}}^2|$ is bounded by
\[  \left|\int_C | \dot{\eta}-\bar{\partial}_E\dot{h}_{\mathrm{app}}|_{h_t^{\mathrm{app}}}^2-| \dot{\eta}-\bar{\partial}_E\dot{h}|_{h_t}^2\right|+\left|\int_C |\dot{\varphi}+[\dot{h}_{\mathrm{app}},\varphi_t]|_{h_t^{\mathrm{app}}}^2-| \dot{\varphi}+[\dot{h},\varphi_t]|_{h_t}^2\right|:=\uppercase\expandafter{\romannumeral1}+\uppercase\expandafter{\romannumeral2}.\]
We first estimate $\uppercase\expandafter{\romannumeral1}$, \begin{align}
\uppercase\expandafter{\romannumeral1}&\leqslant \int_C |\langle (\mathrm{Id}-\mathrm{e}^{-\gamma_t})(\dot{\eta}-\bar{\partial}_E\dot{h}_{\mathrm{app}}),\dot{\eta}-\bar{\partial}_E\dot{h}_{\mathrm{app}}\rangle_{h_t^{\mathrm{app}}}|\notag\\&\quad+|\langle \mathrm{e}^{-\gamma_t}(\dot{\eta}-\bar{\partial}_E\dot{h}_{\mathrm{app}}),(\mathrm{Id}-\mathrm{e}^{-\gamma_t})(\dot{\eta}-\bar{\partial}_E\dot{h}_{\mathrm{app}})\rangle_{h_t^{\mathrm{app}}}|\notag\\
&\quad+2|\langle \mathrm{e}^{-\gamma_t}\bar{\partial}_E(\dot{h}_{\mathrm{app}}-\dot{h}),\mathrm{e}^{-\gamma_t}(\dot{\eta}-\bar{\partial}_E\dot{h}_{\mathrm{app}}\rangle_{h_t^{\mathrm{app}}}|+|\mathrm{e}^{-\gamma_t}\bar{\partial}_E(\dot{h}_{\mathrm{app}}-\dot{h})|_{h_t^{\mathrm{app}}}^2.\label{L2AppCompTerm1_eq}
\end{align}
Since the $L_{-2+\delta}^{2,2}$ norm of $\gamma_t$ is $O(t^{-N})$, the first two terms are bounded by \[c_Nt^{-N}\lVert \dot{\eta}-\bar{\partial}_E\dot{h}_{\mathrm{app}}\rVert_{L^2}^2.\]
Here and in the following, unless specified otherwise, the norms are taken with respect to $h_t^{\mathrm{app}}$ as in \eqref{SobolevSpace_eq}. $N$ and $c_N$ may vary from line to line. It remains to estimate $\lVert\nabla_{A_t}\dot{\mathbf{v}}\rVert_{L^2}$, where $\dot{\mathbf{v}}=\dot{h}_{\mathrm{app}}-\dot{h}$. By \eqref{LinHitCou_eq}, \eqref{AppLinHitCol_eq}, we have
\begin{align*}
  &\partial_{E}^{h_t^{\mathrm{app}}}\bar{\partial}_E \dot{\mathbf{v}}-[\varphi_t^{\ast_{h_t^{\mathrm{app}}}},[\dot{\mathbf{v}},\varphi_t]]+T_t^1(\dot{\eta}-\bar{\partial}_E\dot{h})+T_t^2(\dot{\varphi}+[\dot{h},\varphi_t])=0,\text{ where }\\
  &T_t^1=\partial_E^{h_t}-\partial_E^{h_t^{\mathrm{app}}}=-2[\partial_E^{h_t^{\mathrm{app}}}\gamma_t,\cdot]+R_t^1,\quad T_t^2=[\varphi_t^{\ast_{h_t}}-\varphi_t^{\ast_{h_t^{\mathrm{app}}}},\cdot]=-2[[\gamma_t,\varphi_t^{\ast_{h_t^{\mathrm{app}}}}],\cdot]+R_t^2,
\end{align*}
$R_t^1,R_t^2$ are nonlinear in $\gamma_t$. As before, by composing with $-2\mathrm{i}\star$, the equation becomes
\[L_t\dot{\mathbf{v}}=2\mathrm{i}\star(T_t^1(\dot{\eta}-\bar{\partial}_E\dot{h})+T_t^2(\dot{\varphi}+[\dot{h},\varphi_t])+E_t\dot{\mathbf{v}}):=\dot{\mathbf{w}},\]
where $E_t=[F_{h_t^{\mathrm{app}}}+[\varphi_t,\varphi_t^{\ast_{h_t^{\mathrm{app}}}}],\cdot]$. Note that $\nabla_{A_t}\gamma_t, [\Psi_t,\gamma_t]\in L_{-2+\delta}^{1,2}$ and $L_{-2+\delta}^{1,2}\hookrightarrow L_{-1+\delta}^p$ for any fixed $p\in [2,\infty)$, $L_{-1+\delta}^p\cdot L^2\subset L_{-2+\delta}^{2p/(p+2)}$, similar to the proof of Theorem \ref{DefSol_thm}, we have
 \[\lVert\star T_t^1(\dot{\eta}-\bar{\partial}_E\dot{h})\rVert_{L_{-2+\delta}^{p_1}}+\lVert\star T_t^2(\dot{\varphi}+[\dot{h},\varphi_t])\rVert_{L_{-2+\delta}^{p_1}}\leqslant c_Nt^{-N}\lVert [(\dot{\eta},\dot{\varphi})]\rVert_{g_{L^2}},\]
 where $p_1=2p/(p+2)<2$. By the construction of $h_t^{\mathrm{app}}$, $E_t$ is supported on a union of annuli and the $C^0$ norm is $O(t^{-N})$, so that
 \[\lVert E_t\dot{\mathbf{v}}\rVert_{L_{-2+\delta}^{p_1}}\leqslant c_Nt^{-N}\lVert \dot{\mathbf{v}}\rVert_{L_{-2+\delta}^{p_1}}.\] Therefore,
\begin{equation}\label{p1NormEst_eq1}
  \lVert \dot{\mathbf{w}}\rVert_{L_{-2+\delta}^{p_1}}\leqslant c_Nt^{-N}(\lVert [(\dot{\eta},\dot{\varphi})]\rVert_{g_{L^2}}+\lVert \dot{\mathbf{v}}\rVert_{L_{-2+\delta}^{p_1}}).
\end{equation} By \eqref{intbyParts_eq},
\begin{equation}\label{p1NormEst_eq2}
  \lVert \nabla_{A_t}\dot{\mathbf{v}}\rVert_{L^2}^2+\lVert [\Psi_t,\dot{\mathbf{v}}]\rVert_{L^2}^2=\langle \dot{\mathbf{w}},\dot{\mathbf{v}}\rangle_{L^2}\leqslant \lVert \dot{\mathbf{w}}\rVert_{L_{-2+\delta}^{p_1}}\lVert \dot{\mathbf{v}}\rVert_{L_{-\delta}^{p_1^\ast}},
\end{equation}
 where $1/p_1+1/p_1^\ast=1$. Let $u$ be a function on the unit disk vanishing on the boundary, then we have the Sobolev inequality $\lVert u\rVert_{L^{p_1^\ast}}\leqslant c\lVert \nabla u\rVert_{L^2}$. On the other hand, by Lemma \ref{SobEmbed_lem}, $L_{-\delta-1}^{1,2}=\hat{L}_{-\delta-1}^{1,2}\hookrightarrow L_{-\delta}^{p_1^\ast}$, then \[\lVert u\rVert_{L_{-\delta}^{p_1^\ast}}\leqslant c(\lVert\nabla u\rVert_{L_{-\delta-1}^2}+\lVert u\rVert_{L_{-\delta-1}^2})\leqslant c(\lVert\nabla u\rVert_{L_{-\delta-1}^2}+\lVert\nabla u\rVert_{L_{-\delta-2}^2})\leqslant c\lVert \nabla u\rVert_{L^2}.\] These estimates can be translated to $C$ as in the proof of \eqref{PoincareIneq_eq}, and we obtain
\[  \lVert \dot{\mathbf{v}}\rVert_{L_{-\delta}^{p_1^\ast}}\leqslant ct^M(\lVert \nabla_{A_t}\dot{\mathbf{v}}\rVert_{L^2}+\lVert [\Psi_t,\dot{\mathbf{v}}]\rVert_{L^2}),\]
for some constant $M$. Combining with \eqref{p1NormEst_eq1}, \eqref{p1NormEst_eq2}, we deduce that \begin{align*}
  \lVert \nabla_{A_t}\dot{\mathbf{v}}\rVert_{L^2}+\lVert [\Psi_t,\dot{\mathbf{v}}]\rVert_{L^2}&\leqslant ct^M \lVert\dot{\mathbf{w}}\rVert_{L_{-2+\delta}^{p_1}}\leqslant c_{N}t^{-N}(\lVert [(\dot{\eta},\dot{\varphi})]\rVert_{g_{L^2}}+\lVert \dot{\mathbf{v}}\rVert_{L_{-2+\delta}^{p_1}})\\
  &\leqslant c_{N}t^{-N}(\lVert [(\dot{\eta},\dot{\varphi})]\rVert_{g_{L^2}}+\lVert\dot{\mathbf{v}}\rVert_{L_{-\delta}^{p_1^\ast}})\\
  &\leqslant c_{N}t^{-N}(\lVert [(\dot{\eta},\dot{\varphi})]\rVert_{g_{L^2}}+\lVert \nabla_{A_t}\dot{\mathbf{v}}\rVert_{L^2}+\lVert [\Psi_t,\dot{\mathbf{v}}]\rVert_{L^2}),\\
   \lVert \nabla_{A_t}\dot{\mathbf{v}}\rVert_{L^2}+\lVert [\Psi_t,\dot{\mathbf{v}}]\rVert_{L^2}&\leqslant c_{N}t^{-N}\lVert [(\dot{\eta},\dot{\varphi})]\rVert_{g_{L^2}},
\end{align*}
for $t$ large. By \eqref{L2AppCompTerm1_eq}, we conclude that
\[\uppercase\expandafter{\romannumeral1}\leqslant c_{N}t^{-N}(\lVert \dot{\eta}-\bar{\partial}_E\dot{h}_{\mathrm{app}}\rVert_{L^2}^2+\lVert [(\dot{\eta},\dot{\varphi})]\rVert_{g_{L^2}}^2).\]
Similarly, with the roll of $\lVert \nabla_{A_t}\dot{\mathbf{v}}\rVert_{L^2}$ replaced by $\lVert [\Psi_t,\dot{\mathbf{v}}]\rVert_{L^2}$, we have
\[\uppercase\expandafter{\romannumeral2}\leqslant c_Nt^{-N}(\lVert \dot{\varphi}+[\dot{h}_{\mathrm{app}},\varphi_t]\rVert_{L^2}^2+\lVert [(\dot{\eta},\dot{\varphi})]\rVert_{g_{L^2}}^2).\]
Therefore $|\lVert [(\dot{\eta},\dot{\varphi})]\rVert_{g_{\mathrm{app}}}^2-\lVert [(\dot{\eta},\dot{\varphi})]\rVert_{g_{L^2}}^2|\leqslant c_Nt^{-N}\lVert [(\dot{\eta},\dot{\varphi})]\rVert_{g_{L^2}}^2$, meaning that along the curve $[(\bar{\partial}_E,\varphi_t)]$ in $\mathcal{M}$, $g_{L^2}$ and $g_{\mathrm{app}}$ are polynomially close of any order in $t$. If $P_w=\varnothing$, according to \eqref{ErrorEst_eq1} and Theorem \ref{DefSol_thm}, this can be improved to an exponential decay.

\subsection{Comparing $g_{\mathrm{app}}$ and $g_{\mathrm{sf}}$}
By the construction, $h_t^{\mathrm{app}}$ and $h_t^{\mathrm{dh}}$ coincide on the complement of the disks $\widetilde{B}_{x,j,t}$ ($x\in I$) and $\widetilde{B}_{x,t}$ ($x\in P$), and so do the infinitesimal deformations $\dot{h}_{\mathrm{app}}$ and $\dot{h}_{\mathrm{dh}}$. Accordingly, the difference $g_{\mathrm{app}}-g_{\mathrm{sf}}$ localizes on these disks. Near zeros and strongly parabolic points, the results in \cite{fredrickson_2019, fredrickson2022asymptotic} directly apply.

  On $\widetilde{B}_{x,j,t}$, in a local holomorphic coordinate $\zeta_{x,j,t}$ and a holomorphic frame we have \[\bar{\partial}_E=\bar{\partial},\quad \varphi_t=\lambda_{x,j}(t)\begin{pmatrix}
    0&1\\\zeta_{x,j,t}&0
  \end{pmatrix}\,\mathrm{d}\zeta_{x,j,t},\quad h_t^{\mathrm{dh}}=\begin{pmatrix}
    |\zeta_{x,j,t}|^{1/2}&0\\0& |\zeta_{x,j,t}|^{-1/2}
  \end{pmatrix},\]
and $h_t^{\mathrm{app}}$ is given by \eqref{approxMetricZ_eq}. Let $\hat{\zeta}_{x,j,t}=\kappa^{-1}t^{1/j(x)}\zeta_{x,j,t}$ and $\hat{r}_{x,j,t}=|\hat{\zeta}_{x,j,t}|$, then the disk $\widetilde{B}_{x,j,t}=\{\,\hat{r}_{x,j,t}<1\,\}$. Applying the local holomorphic gauge transformation
\[  g=\begin{pmatrix}
    \kappa^{-1/4}t^{1/4j(x)}&0\\0&\kappa^{1/4}t^{-1/4j(x)}
  \end{pmatrix},\text{ then }
  g^{-1}\varphi_tg =t_{x,j}\begin{pmatrix}
    0&1\\\hat{\zeta}_{x,j,t}&0
\end{pmatrix}\,\mathrm{d}\hat{\zeta}_{x,j,t},\]
where $t_{x,j}=\lambda_{x,j}(t)\kappa^{3/2}t^{-3/2j(x)}$.
\[ g^\ast h_t^{\mathrm{dh}}g=\begin{pmatrix}
  \hat{r}_{x,j,t}^{1/2}&0\\0&\hat{r}_{x,j,t}^{-1/2}
\end{pmatrix},\quad g^\ast h_t^{\mathrm{app}}g=\begin{pmatrix}
  \hat{r}_{x,j,t}^{1/2}\mathrm{e}^{l_{t_{x,j}(\hat{r}_{x,j,t})}\chi(\hat{r}_{x,j,t})}&0\\0&\hat{r}_{x,j,t}^{-1/2}\mathrm{e}^{-l_{t_{x,j}}(\hat{r}_{x,j,t})\chi(\hat{r}_{x,j,t})}
\end{pmatrix},\]
where $l_{t_{x,j}}(\hat{r}_{x,j,t})=\psi_1(8t_{x,j}\hat{r}_{x,j,t}^{3/2}/3)$. In this local holomorphic frame, the following result holds.
\begin{proposition}[{\cite[Th.~3.1]{fredrickson_2019}}]\label{MetCompZero_Prop}
Let $[(\dot{\eta},\dot{\varphi})]\in T_{[(\bar{\partial}_E,\varphi_t)]}\mathcal{M}$, then on the disk $\{\,\hat{r}_{x,j,t}<1\,\}$ there is a unique representative in the equivalence class where
\[\dot{\eta}=0,\quad \dot{\varphi}=t_{x,j}\begin{pmatrix}
  0&0\\\dot{P}&0
\end{pmatrix}\,\mathrm{d}\hat{\zeta}_{x,j,t},\quad \dot{h}_{\mathrm{dh}}=\frac{\dot{P}}{4\hat{\zeta}_{x,j,t}}\begin{pmatrix}
  1&0\\0&-1
\end{pmatrix}\quad\text{for }\bar{\partial}\dot{P}=0.\]
These infinitesimal deformations satisfy \eqref{DecoupInftyDef_eq}, and $\dot{h}_{\mathrm{app}}$ solving \eqref{AppLinHitCol_eq} is diagonal. For any $[(\dot{\eta},\dot{\varphi})]\in T_{[(\bar{\partial}_E,\varphi_t)]}\mathcal{M}$, there exists a constant $c>0$ such that \[\lVert [(\dot{\eta},\dot{\varphi})]\rVert_{g_{\mathrm{app}}(\widetilde{B}_{x,j,t})}^2-\lVert [(\dot{\eta},\dot{\varphi})]\rVert_{g_{\mathrm{sf}}(\widetilde{B}_{x,j,t})}^2=O(\mathrm{e}^{-ct_{x,j}}).\]
\end{proposition}
On $\widetilde{B}_{x,t}$ for $x\in P_s$, in a local holomorphic coordinate $\zeta_{x,t}$ and holomorphic frame we have
\[\bar{\partial}_E=\bar{\partial},\quad \varphi_t=\lambda_{x}(t)\begin{pmatrix}
  0&1\\\frac{1}{\zeta_{x,t}}&0
\end{pmatrix}\,\mathrm{d}\zeta_{x,t},\quad h_t^{\mathrm{dh}}=\begin{pmatrix}
  |\zeta_{x,t}|^{1/2}&0\\0& |\zeta_{x,t}|^{3/2}
\end{pmatrix},\]
and $h_t^{\mathrm{app}}$ is given by \eqref{approxMetricT_eq}. Let $\hat{\zeta}_{x,t}=\kappa^{-1}\zeta_{x,t}$ and $\hat{r}_{x,t}=|\hat{\zeta}_{x,t}|$, then the disk $\widetilde{B}_{x,t}=\{\,\hat{r}_{x,t}<1\,\}$. Applying the local holomorphic gauge transformation
\[  g=\kappa^{-1/2}\begin{pmatrix}
    \kappa^{1/4}&0\\0&\kappa^{-1/4}
  \end{pmatrix},\text{ then }
  g^{-1}\varphi_tg =t_x\begin{pmatrix}
    0&1\\\hat{\zeta}_{x,t}&0
\end{pmatrix}\,\mathrm{d}\hat{\zeta}_{x,t},\]
where $t_x=\lambda_x(t)\kappa^{1/2}$.
\[ g^\ast h_t^{\mathrm{dh}}g=\begin{pmatrix}
  \hat{r}_{x,t}^{1/2}&0\\0&\hat{r}_{x,t}^{3/2}
\end{pmatrix},\quad g^\ast h_t^{\mathrm{app}}g=\begin{pmatrix}
  \hat{r}_{x,t}^{1/2}\mathrm{e}^{m_{t_{x}(\hat{r}_{x,t})}\chi(\hat{r}_{x,t})}&0\\0&\hat{r}_{x,t}^{3/2}\mathrm{e}^{-m_{t_{x}}(\hat{r}_{x,t})\chi(\hat{r}_{x,t})}
\end{pmatrix},\]
where $m_{t_x}(\hat{r}_{x,t})=\psi_{2,x}(8t_x\hat{r}_{x,t}^{1/2})$. In this local holomorphic frame, the following result holds.
\begin{proposition}[{\cite[Lem.~7.10,~Prop.~7.11]{fredrickson2022asymptotic}}]\label{MetCompSPara_Prop}
Let $[(\dot{\eta},\dot{\varphi})]\in T_{[(\bar{\partial}_E,\varphi_t)]}\mathcal{M}$, then on the disk $\{\,\hat{r}_{x,t}<1\,\}$ there is a unique representative in the equivalence class where
\[\dot{\eta}=0,\quad \dot{\varphi}=t_x\begin{pmatrix}
  0&0\\\frac{\dot{P}}{\hat{\zeta}_{x,t}}&0
\end{pmatrix}\,\mathrm{d}\hat{\zeta}_{x,t},\quad \dot{h}_{\mathrm{dh}}=\frac{\dot{P}}{4}\begin{pmatrix}
  1&0\\0&-1
\end{pmatrix}\quad\text{for }\bar{\partial}\dot{P}=0.\]
These infinitesimal deformations satisfy \eqref{DecoupInftyDef_eq}. For any $[(\dot{\eta},\dot{\varphi})]\in T_{[(\bar{\partial}_E,\varphi_t)]}\mathcal{M}$, there exists a constant $c>0$ such that \[\lVert [(\dot{\eta},\dot{\varphi})]\rVert_{g_{\mathrm{app}}(\widetilde{B}_{x,t})}^2-\lVert [(\dot{\eta},\dot{\varphi})]\rVert_{g_{\mathrm{sf}}(\widetilde{B}_{x,t})}^2=O(\mathrm{e}^{-ct_{x}}).\]
\end{proposition}
 On $\widetilde{B}_{x,t}$ for $x\in P_w$, as in Proposition \ref{NormalForm_prop}.(\romannumeral2), choose a holomorphic coordinate $\zeta_{x,t}$ and a holomorphic frame in which \[    \bar{\partial}_E=\bar{\partial},\quad \varphi_t=\frac{1}{\zeta_{x,t}}\begin{pmatrix}
        0&-\mu_{x,2}(1-\epsilon_{x,0,t}^{-1}\zeta_{x,t})\\1&0
      \end{pmatrix}\,\mathrm{d}\zeta_{x,t}.\] We have the following standard representatives of infinitesimal deformations.
\begin{lemma}
There exists a representative $(\dot{\eta},\dot{\varphi})$ of $[(\dot{\eta},\dot{\varphi})]\in T_{[(\bar{\partial}_E,\varphi_t)]}\mathcal{M}$ on $\widetilde{B}_{x,t}$ where
\[\dot{\eta}=0,\quad \dot{\varphi}=\begin{pmatrix}
0& \dot{P}\\ 0&0
\end{pmatrix}\,\mathrm{d}\zeta_{x,t},\quad\dot{P}\text{ is holomorphic}.\]
\end{lemma}
\begin{proof}
  The proof is similar to Lemma \ref{InfsmRep_lem}. By the Poincaré lemma, we can make $\dot{\eta}=0$. Since the residue of $\dot{\varphi}$ preserves the flag at $\zeta_{x,t}=0$, we have
  \[\dot{\varphi}=
    \frac{\mathrm{d}\zeta_{x,t}}{\zeta_{x,t}}\begin{pmatrix}
      \dot{P}_1& 2\dot{P}_1+\dot{P}_3+\dot{P}_2\zeta_{x,t}\\ \dot{P}_3&-\dot{P}_1
  \end{pmatrix}.
  \]
  Similarly, since $\dot{\gamma}$ preserves the flag, we have
  \[\dot{\gamma}=\begin{pmatrix}
    \dot{\gamma}_1&\dot{\gamma}_1+\dot{\gamma}_3-\dot{\gamma}_4+\dot{\gamma}_2\zeta_{x,t}\\ \dot{\gamma}_3&\dot{\gamma}_4
  \end{pmatrix}, \]
   and then
   \[[\varphi_t,\dot{\gamma}]=\frac{\mathrm{d}\zeta_{x,t}}{\zeta_{x,t}}\begin{pmatrix}
     -\mu_{x,2}(1-\epsilon_{x,0,t}^{-1}\zeta_{x,t})\dot{\gamma}_3-(\dot{\gamma_1}+\dot{\gamma}_3-\dot{\gamma_4}+\dot{\gamma}_2\zeta_{x,t})&-\mu_{x,2}(1-\epsilon_{x,0,t}^{-1}\zeta_{x,t})(\dot{\gamma}_4-\dot{\gamma}_1)\\\dot{\gamma}_1-\dot{\gamma}_4& \ast
   \end{pmatrix},\]
   where the $(2,2)$-entry is the negative of the $(1,1)$-entry. Choose $\dot{\gamma_1}=\dot{P}_1$, and $\dot{\gamma}_k=0$ ($k=2,3,4$), then $\dot{\varphi}$ becomes off-diagonal, and can be written as
   \[\dot{\varphi}=\frac{\mathrm{d}\zeta_{x,t}}{\zeta_{x,t}}\begin{pmatrix}
     0&\dot{P}_3+\dot{P}_2\zeta_{x,t}\\ \dot{P}_3&0
   \end{pmatrix}.\]
   Since the singularity data at $x$ is fixed, $\dot{P}_3$ must vanish at $0$, so that $\dot{P}_3=\zeta_{x,t}\dot{P}'_{3}$ for $\dot{P}'_{3}$ holomorphic. Next we choose $\dot{\gamma}_1=-\dot{P}_3$, $\dot{\gamma}_2=\dot{P}'_3$, $\dot{\gamma}_3=\dot{\gamma}_4=0$, then $\dot{\varphi}$ has the desired form.
\end{proof}
To prove the decay of the difference between $g_{\mathrm{app}}$ and $g_{\mathrm{sf}}$, we adapt the proof of \cite[Th.~3.1]{fredrickson_2019} and \cite[Prop.~7.11]{fredrickson2022asymptotic} to weakly parabolic points, which is based on the method developed in \cite{dumas2019asymptotics}. To apply this method, it is crucial that the (decoupled) harmonic metric and its infinitesimal deformation are diagonal. We make the following assumption.
\begin{assumption}
  The parabolic data at weakly parabolic points are trivial, i.e., $\alpha_{x,1}=\alpha_{x,2}=1/2$.
\end{assumption}
Let $(\hat{\mathcal{E}},\hat{\varphi})$ be the model Higgs bundle defined in section \ref{ModSol_subsubsec} with trivial parabolic structure at $0$, $\hat{h}$ be a compatible harmonic metric. Define $g:\hat{\mathcal{E}}\to\hat{\mathcal{E}}$ by $(e_1,e_2)\mapsto (e_1,-e_2)$. Then $g^\ast \hat{h}$ is also a compatible harmonic metric. By uniqueness, there is $c>0$ such that $g\cdot h=ch$. Then $-\hat{h}(e_1,e_2)=ch(e_1,e_2)$, implying that $\hat{h}(e_1,e_2)=0$, i.e., $\hat{h}$ is diagonal. $\hat{h}^{\mathrm{dh}}$ is also diagonal, and is explicitly given by
\begin{equation}\label{HolVectF_eq}
  \hat{h}^{\mathrm{dh}}=|z|\begin{pmatrix}
    |\mu_{x,2}(1-z)|^{-1/2}&0\\ 0&  |\mu_{x,2}(1-z)|^{1/2}
  \end{pmatrix}.
\end{equation}
\begin{proposition}\label{MetCompWPara_Prop}
  Let $(0,\dot{\varphi})$ be the local representative in the previous lemma, then
  \begin{equation}\label{WParaComp_eq}
      \lVert (0,\dot{\varphi}, \dot{h}_{\mathrm{app}})\rVert_{g_{\mathrm{app}}(\widetilde{B}_{x,t})}^2-\lVert (0,\dot{\varphi},\dot{h}_{\mathrm{dh}})\rVert_{g_{\mathrm{dh}}(\widetilde{B}_{x,t})}^2=O(t^{-N}),
  \end{equation}
  for any $N>0$.
\end{proposition}
\begin{proof}
By the proof of \cite[Lem.~7.10]{fredrickson2022asymptotic}, $\dot{h}_{\mathrm{dh}}$ is diagonal, possibly after applying an infinitesimal gauge transformation preserving $(0,\dot{\varphi})$. Moreover, $\dot{h}_{\mathrm{app}}-\dot{h}_{\mathrm{dh}}=O(t^{-N})$ on $\widetilde{B}_{x,t}\backslash \frac{1}{2}\widetilde{B}_{x,t}$, where $\frac{1}{2}\widetilde{B}_{x,t}$ is the disk $\{|\zeta_{x,t}|<\frac{1}{2}\kappa t^{-2/3}\}$. Then the difference in \eqref{WParaComp_eq} is $O(t^{-N})$ on $\widetilde{B}_{x,t}\backslash \frac{1}{2}\widetilde{B}_{x,t}$. Introduce the notation \[\delta\big((0,\dot{\varphi}),(h_1,\dot{h}_1),(h_2,\dot{h}_2)\big)=(\langle \dot{\varphi}+[\dot{h}_1,\varphi_t],\dot{\varphi}\rangle_{h_1}-\langle \dot{\varphi}+[\dot{h}_2,\varphi_t],\dot{\varphi}\rangle_{h_2})\,\mathrm{dvol}.\]
It remains to prove that the following difference is $O(t^{-N})$:
\[\lVert (0,\dot{\varphi}, \dot{h}_{\mathrm{app}})\rVert_{g_{\mathrm{app}}(\frac{1}{2}\widetilde{B}_{x,t})}^2-\lVert (0,\dot{\varphi},\dot{h}_{\mathrm{dh}})\rVert_{g_{\mathrm{dh}}(\frac{1}{2}\widetilde{B}_{x,t})}^2=\int_{\frac{1}{2}\widetilde{B}_{x,t}}\delta\big((0,\dot{\varphi}),(h_t^{\mathrm{app}},\dot{h}_{\mathrm{app}}),(h_t^{\mathrm{dh}},\dot{h}_{\mathrm{dh}}))\big).\]
Note that $h_t^{\mathrm{app}}=h_t^{\mathrm{model}}:=f_{x,t}^{\ast}\hat{\sigma}_{x,t}^{\ast}\hat{h}$ on $\frac{1}{2}\widetilde{B}_{x,t}$, we can decompose the integrand as
\begin{equation}\label{DecompDiffAppDh_eq}
  \delta\big((0,\dot{\varphi}),(h_t^{\mathrm{model}},\dot{h}_X),(h_t^{\mathrm{dh}},\dot{h}_{\mathrm{dh}}))\big)+\delta\big((0,\dot{\varphi}),(h_t^{\mathrm{model}},\dot{h}_{\mathrm{app}}),(h_t^{\mathrm{model}},\dot{h}_{X}))\big),
\end{equation}
where, analogous to \eqref{LinHitCou_eq}, $\dot{h}_X$ solves
\begin{equation}\label{LinHitCouPwModel_eq}
  \partial_E^{h_t^\mathrm{model}}\bar{\partial}_E\dot{h}_X-[\varphi_t^{\ast_{h_t^{\mathrm{model}}}},\dot{\varphi}+[\dot{h}_X,\varphi_t]]=0.
\end{equation}
We construct $\dot{h}_X$ following the method of Dumas-Neitzke \cite{dumas2019asymptotics}. First we write
\begin{align*}
  \varphi_t&=\frac{\mathrm{d}\zeta_{x,t}}{\zeta_{x,t}}\begin{pmatrix}
    0&P\\1&0
  \end{pmatrix},\quad \dot{\varphi}=\begin{pmatrix}
    0&\dot{P}\\0&0
  \end{pmatrix}\mathrm{d}\zeta_{x,t},\\\quad h_t^{\mathrm{model}}&=|\zeta_{x,t}|\begin{pmatrix}
    \mathrm{e}^{-u}&0\\0&\mathrm{e}^u
  \end{pmatrix},\quad \dot{h}_X=\begin{pmatrix}
    -\frac{1}{2}F_X&0\\0&\frac{1}{2}F_X
\end{pmatrix},\\
h_t^{\mathrm{dh}}&=|\zeta_{x,t}|\begin{pmatrix}
  |P|^{-1/2}&0\\0&|P|^{1/2}
\end{pmatrix},\quad \dot{h}_{\mathrm{dh}}=\begin{pmatrix}
  -\frac{1}{2}F_{\mathrm{dh}}&0\\0&\frac{1}{2}F_{\mathrm{dh}}.
\end{pmatrix}.
\end{align*}
Here $P=-\mu_{x,2}(1-\epsilon_{x,0,t}^{-1}\zeta_{x,t})$. Now \eqref{LinHitCouPwModel_eq} is equivalent to
\[(\Delta-8(\mathrm{e}^{2u}+ \mathrm{e}^{-2u}|P|^2)|\zeta_{x,t}|^{-2})F_X+8|\zeta_{x,t}|^{-2}\zeta_{x,t}\mathrm{e}^{-2u}\widebar{P}\dot{P}=0.\]
Let $P=\zeta_{x,t}^2 P_1, \dot{P}=\zeta_{x,t}\dot{P}_1, \mathrm{e}^{2u}=|\zeta_{x,t}|^2 \mathrm{e}^{2u_1}$, the equation becomes
\[(\Delta-8(\mathrm{e}^{2u_1}+ \mathrm{e}^{-2u}|P_1|^2))F_X+8\mathrm{e}^{-2u_1}\widebar{P_1}\dot{P}_1=0,\]
which is exactly \cite[Eq.~(4.14)]{dumas2019asymptotics}. By \cite[Th.~11]{dumas2019asymptotics}, we have \[F_X=\mathcal{X}'+2\mathcal{X} u_1',\]
where $\mathcal{X}'=\frac{\partial}{\partial \zeta_{x,t}}\mathcal{X}$ (similar for $u_1'$), $X=\mathcal{X}(\zeta_{x,t})\frac{\partial}{\partial \zeta_{x,t}}$ is the holomorphic vector field such that
\[\mathcal{L}_X (P_1\,\mathrm{d}\zeta_{x,t}^2)=(\mathcal{X} P_1'+2P_1 \mathcal{X}' )\mathrm{d}\zeta_{x,t}^2=\dot{P}_1\,\mathrm{d}\zeta_{x,t}^2.\]
In other words, $\mathcal{X}$ solves the ODE
\[\mathcal{X}'+\frac{2-\epsilon_{x,0,t}^{-1}\zeta_{x,t}}{-2\zeta_{x,t}(1-\epsilon_{x,0,t}^{-1}\zeta_{x,t})}\mathcal{X}=\frac{\zeta_{x,t}^2\dot{P}_1}{-2\mu_{x,2}(1-\epsilon_{x,0,t}^{-1}\zeta_{x,t})}.\]
$\mathcal{X}$ is explicitly given by
\[\mathcal{X}=\frac{\zeta_{x,t}}{\sqrt{1-\epsilon_{x,0,t}^{-1}\zeta_{x,t}}}\int \frac{\zeta_{x,t}\dot{P}_1}{-2\mu_{x,2}\sqrt{1-\epsilon_{x,0,t}^{-1}\zeta_{x,t}}},\]
where the integral is the term-by-term integration of the power series expansion of the integrand at $\zeta_{x,t}=\epsilon_{x,0,t}$, which is convergent to $\sqrt{1-\epsilon_{x,0,t}^{-1}\zeta_{x,t}}\widetilde{P}(\zeta_{x,t})$ for $\widetilde{P}$ holomorphic in $\widetilde{B}_{x,t}$. Then $\mathcal{X}=\zeta_{x,t}\widetilde{P}$ is holomorphic in $\widetilde{B}_{x,t}$ and vanishes at $0$. Replacing $u_1$ by $\frac{1}{2}\log |P_1|$, we obtain
\[F_{\mathrm{dh}}=\mathcal{X}'+\mathcal{X} (\log |P_1|)'=\mathcal{X}'+\mathcal{X} (2P_1)^{-1}P_1'.\]
The term $\delta\big((0,\dot{\varphi}),(h_t^{\mathrm{model}},\dot{h}_X),(h_t^{\mathrm{dh}},\dot{h}_{\mathrm{dh}}))\big)$ in \eqref{DecompDiffAppDh_eq} is
\begin{align*}
  &\phantom{=}\Big((\mathrm{e}^{-2u}-|P|^{-1})|\dot{P}|^2-\zeta_{x,t}^{-1}\mathrm{e}^{-2u}F_X P\widebar{\dot{P}}+\zeta_{x,t}^{-1}|P|^{-1}F_{\mathrm{dh}}P\widebar{\dot{P}}\Big)\mathrm{i}\,\mathrm{d}\zeta_{x,t}\mathrm{d}\widebar{\zeta_{x,t}}\\
  &=\Big((\mathrm{e}^{-2u_1}-|P_1|^{-1})|\dot{P}_1|^2-\mathrm{e}^{-2u_1}((2P_1)^{-1}(\dot{P}_1-\mathcal{X}P_1')+2\mathcal{X}u_1')P_1\widebar{\dot{P_1}}\\
  &\qquad +|P_1|^{-1}((2P_1)^{-1}(\dot{P}_1-\mathcal{X}P_1')+\mathcal{X}(2P_1)^{-1}P_1')P_1\widebar{\dot{P_1}}\Big)\mathrm{i}\,\mathrm{d}\zeta_{x,t}\mathrm{d}\widebar{\zeta_{x,t}}\\
  &=\Big((\mathrm{e}^{-2u_1}-|P_1|^{-1})|\dot{P}_1|^2/2+\mathrm{e}^{-2u_1}\mathcal{X}P_1'\widebar{\dot{P_1}}/2-2\mathrm{e}^{-2u_1}u_1'\mathcal{X}P_1\widebar{\dot{P_1}}\Big)\mathrm{i}\,\mathrm{d}\zeta_{x,t}\mathrm{d}\widebar{\zeta_{x,t}}.
\end{align*}
On the other hand, let $\beta=(\mathrm{e}^{-2u_1}-|P_1|^{-1})\mathcal{X}P_1\widebar{\dot{P_1}} \mathrm{i}\,\mathrm{d}\widebar{\zeta_{x,t}}$, we have
\begin{align*}
  \mathrm{d}\beta&=\Big(-2\mathrm{e}^{-2u_1}u_1'\mathcal{X}P_1\widebar{\dot{P_1}}+|P_1|^{-1}P_1'\mathcal{X}\widebar{\dot{P_1}}/2+(\mathrm{e}^{-2u_1}-|P_1|^{-1})((2P_1)^{-1}(\dot{P}_1-\mathcal{X}P_1'))\widebar{\dot{P_1}}P_1\\&\qquad+(\mathrm{e}^{-2u_1}-|P_1|^{-1})\mathcal{X}P_1'\widebar{\dot{P_1}}\Big)\mathrm{i}\,\mathrm{d}\zeta_{x,t}\mathrm{d}\widebar{\zeta_{x,t}}\\
  &=\delta\big((0,\dot{\varphi}),(h_t^{\mathrm{model}},\dot{h}_X),(h_t^{\mathrm{dh}},\dot{h}_{\mathrm{dh}}))\big).
\end{align*}
Let $B(r):=\{ |\zeta_{x,t}-\epsilon_{x,0,t}|<r\kappa t^{-2/3}\}$, then for $t$ large, $\frac{1}{2}\widetilde{B}_{x,t}\subset B(3/5)\subset B(4/5)\subset \widetilde{B}_{x,t}$.
By Stokes' theorem,
\[\int_{B(r)}\delta\big((0,\dot{\varphi}),(h_t^{\mathrm{model}},\dot{h}_X),(h_t^{\mathrm{dh}},\dot{h}_{\mathrm{dh}}))\big)=\int_{\partial B(r)}\beta.\]
For $3/5\leqslant r \leqslant 4/5$, by Lemma \ref{Fastdecay_lem} we have
\[\int_{\frac{1}{2}\widetilde{B}_{x,t}}\delta\big((0,\dot{\varphi}),(h_t^{\mathrm{model}},\dot{h}_X),(h_t^{\mathrm{dh}},\dot{h}_{\mathrm{dh}}))\big)=O(t^{-N})+\int_{B(r)}\delta\big((0,\dot{\varphi}),(h_t^{\mathrm{model}},\dot{h}_X),(h_t^{\mathrm{dh}},\dot{h}_{\mathrm{dh}}))\big),\]
which is then bounded by (up to $O(t^{-N})$)
\begin{equation}\label{DiffAppSfBd1_eq}
    5\int_{3/5}^{4/5}\int_{\partial B(r)}\beta \mathrm{d}r\leqslant 5\int_{3/5}^{4/5} \Big(\int_{\partial B(r)}|\mathcal{X}|^2 \mathrm{dvol}\Big)^{1/2} \Big(\int_{\partial B(r)}|(\mathrm{e}^{-2u_1}-|P_1|^{-1})P_1|^2|\dot{P}_1|^2 \mathrm{dvol}\Big)^{1/2}\mathrm{d}r.
\end{equation}
By considering the rescaling $\epsilon_{x,0,t}^{-1}\zeta_{x,t}\mapsto \zeta$, $\mathcal{X}$ in \eqref{HolVectF_eq} becomes $t^M$ times a function independent of $t$, and depends linearly on $\dot{P}$. Then similar to \cite[Eq.~(10.17)]{dumas2019asymptotics}, we have
\[|\chi|\leqslant ct^M \int_{\widetilde{B}_{x,t}} \frac{|\dot{P}|^2}{|P|}\,\mathrm{dvol}=ct^M \lVert (0,\dot{\varphi},\dot{h}_{\mathrm{dh}})\rVert_{g_{\mathrm{dh}}(\widetilde{B}_{x,t})}^2\text{ on }\{ |\zeta_{x,t}-\epsilon_{x,0,t}|<2|\epsilon_{x,0,t}|\}.\]
Consequently, for $3/5\leqslant r\leqslant 4/5$,
\begin{align*}
  \int_{\partial B(r)} |\mathcal{X}|^2 \mathrm{dvol}&\leqslant ct^M \lVert (0,\dot{\varphi},\dot{h}_{\mathrm{dh}})\rVert_{g_{\mathrm{dh}}(\widetilde{B}_{x,t})}^2+ c\int_{\partial B(r)}\Big(\int_{2|\epsilon_{x,0,t}|}^{r\kappa t^{-2/3}}\frac{|\dot{P}|}{|P|^{1/2}} \mathrm{d}\mathbf{r}\Big )^2 \mathrm{dvol}\\
  &\leqslant ct^M \Big(\lVert (0,\dot{\varphi},\dot{h}_{\mathrm{dh}})\rVert_{g_{\mathrm{dh}}(\widetilde{B}_{x,t})}^2+ \int_{\partial B(r)}\int_{2|\epsilon_{x,0,t}|}^{r\kappa t^{-2/3}}\frac{|\dot{P}|^2}{|P|}\mathbf{r}\, \mathrm{d}\mathbf{r}\mathrm{d}\boldsymbol{\theta}  \Big)\\
  &\leqslant ct^M \lVert (0,\dot{\varphi},\dot{h}_{\mathrm{dh}})\rVert_{g_{\mathrm{dh}}(\widetilde{B}_{x,t})}^2,
\end{align*}
where we used the Cauchy-Schwarz inequality, and denoted $\zeta_{x,t}-\epsilon_{x,0,t}=\mathbf{r}\mathrm{e}^{\mathrm{i} \boldsymbol{\theta}}$. By Lemma \ref{Fastdecay_lem} and \eqref{DiffAppSfBd1_eq}, we have
\begin{equation}\label{DiffIntEst_eq1}
  \int_{\frac{1}{2}\widetilde{B}_{x,t}}\delta\big((0,\dot{\varphi}),(h_t^{\mathrm{model}},\dot{h}_X),(h_t^{\mathrm{dh}},\dot{h}_{\mathrm{dh}}))\big)\leqslant c_Nt^{-N}\lVert (0,\dot{\varphi},\dot{h}_{\mathrm{dh}})\rVert_{g_{\mathrm{dh}}(\widetilde{B}_{x,t})}^2.
\end{equation}
By the proof of \cite[Prop.~7.11]{fredrickson2022asymptotic}, we have
\begin{equation}\label{DiffIntEst_eq2}
  \int_{\frac{1}{2}\widetilde{B}_{x,t}} \delta\big((0,\dot{\varphi}),(h_t^{\mathrm{model}},\dot{h}_{\mathrm{app}}),(h_t^{\mathrm{model}},\dot{h}_{X}))\big)=O(t^{-N}).
\end{equation}
Combining \eqref{DecompDiffAppDh_eq}, \eqref{DiffIntEst_eq1}, \eqref{DiffIntEst_eq2}, the proposition is proved.
\end{proof}

The Propositions \ref{MetCompZero_Prop}, \ref{MetCompSPara_Prop}, \ref{MetCompWPara_Prop} immediately imply that $g_{\mathrm{app}}-g_{\mathrm{sf}}=O(t^{-N})$ along the curve $[(\bar{\partial}_E,\varphi_t)]$ in $\mathcal{M}$. If moreover $P_w=\varnothing$, the difference is $O(\mathrm{e}^{-ct^{\sigma}})$. Hence the Theorem \ref{Main_thm} is established.
\begin{remark}
 $g_{\mathrm{L^2}}-g_{\mathrm{app}}$ and $g_{\mathrm{app}}-g_{\mathrm{sf}}$ are still polynomially decaying of arbitrary order when $I_u=\varnothing$ and $I_t\neq \varnothing$, only slight modifications of the previous proofs are needed.
\end{remark}

\subsection{Comparing $g_{\mathrm{sf}}$ and $g_{\mathrm{model}}$}
In this subsection we assume that $\mathrm{dim}_{\mathbb{C}}\,\mathcal{M}=2$, or equivalently $N=\mathrm{deg}\,D=4$. The Hitchin base now has complex dimension one, parametrized by $t\in \mathbb{C}$. The constants in this subsection will not depend on $\mathrm{Arg}(t)$ and the choice of the Higgs bundle in the Hitchin fiber $H^{-1}(t)$.
\subsubsection{ALG metrics}\label{ALG_subsubsec}
We first consider the case $D=4\cdot\{0\}$, i.e., there is an untwisted order four pole. The Hitchin base is given by
\begin{equation}\label{HitBaseU4_eq}
  \mathcal{B}=\left\{\,\left(\sum_{k=5}^8\frac{\mu_k}{z^k}+\frac{t}{z^4}\right)\,\mathrm{d}z^2\,\middle|\,t\in \mathbb{C}\,\right\},
\end{equation}
where $\mu_5,\ldots,\mu_8$ are constants determined by the irregular type at $0$, and $\mu_8\neq 0$. We may assume that $\mu_8=-1$, using a rescaling of $z$ if necessary.
Fix \[\nu=\nu(z)\,\mathrm{d}z^2=\left(\sum_{k=5}^8 \frac{\mu_k}{z^k}+\frac{t}{z^4}\right)\,\mathrm{d}z^2\in\mathcal{B}',~\dot{\nu}=\frac{\dot{t}}{z^4}\,\mathrm{d}z^2\in T_\nu\mathcal{B}'.\]
Then the special Kähler metric $g_{\mathrm{sK}}$, which is the restriction of $g_{\mathrm{sf}}$ on $\mathcal{B}'$, is given by \cite[Prop.~8.3]{fredrickson2022asymptotic} \begin{align*}
  g_{\mathrm{sK}}(\dot{\nu},\dot{\nu})&=\int_{C}\frac{|\dot{\nu}|^2}{|\nu|}\,\mathrm{dvol}_C=\frac{|\dot{t}|^2}{|t|}\int_{\mathbb{C}}\frac{1}{\left|\prod_{k=1}^4(z-z_k(t))\right|}~\mathrm{i}\,\mathrm{d}z\mathrm{d}\bar{z},
\end{align*}
where $z_k(t)=t^{-1/4}\mathrm{e}^{\mathrm{i}(k-1)\pi/2}+t^{-1/2}(-1)^k \mu_7 / 4 + t^{-3/4} \mathrm{e}^{-\mathrm{i}(k-1)\pi/2} (-\mu_6/4 - \mu_7^2/32)+O(|t|^{-1})$ ($k=1,\ldots,4$) are the four roots of $z^8\nu(z):=\tilde{\nu}(z)$.  Let $z=t^{-1/4}\xi$ and $z_k=t^{-1/4}\xi_k$, then the integral becomes
\[
|t|^{1/2}\int_{\mathbb{C}}\frac{1}{\left|\prod_{k=1}^4(\xi-\xi_k)\right|}~\mathrm{i}\,\mathrm{d}\xi \mathrm{d}\bar\xi.
\]
We define a funtion \[F(\xi_1,...\xi_4,\bar\xi_1,...\bar\xi_4) := \int_{\mathbb{C}}\frac{1}{\left|\prod_{k=1}^4(\xi-\xi_k)\right|}~\mathrm{i}\,\mathrm{d} \xi \mathrm{d} \bar\xi.\]
Then we can prove that
\[\frac{\partial F}{\partial \xi_k}= \lim_{\epsilon \to 0}\int_{|\xi-\xi_k|>\epsilon}\frac{\bar\xi-\bar\xi_k}{2\left|\xi-\xi_k\right|^3} \frac{1}{\left|\prod_{l\not =k}(\xi-\xi_l)\right|}~\mathrm{i}\,\mathrm{d} \xi \mathrm{d} \bar{\xi},\,\quad\frac{\partial F}{\partial \bar\xi_k}=\overline{\frac{\partial F}{\partial \xi_k}}.\]
Then, for $l\not= k$,
\begin{align*}
\frac{\partial^2 F}{\partial \xi_k\partial \xi_l}= &\lim_{\epsilon\to 0}\int_{|\xi-\xi_k|>\epsilon, |\xi-\xi_l|>\epsilon}\frac{\bar\xi-\bar\xi_k}{2\left|\xi-\xi_k\right|^3}\frac{\bar\xi-\bar\xi_l}{2\left|\xi-\xi_l\right|^3} \frac{1}{\left|\prod_{m\not =k,l}(\xi-\xi_m)\right|}~\mathrm{i}\,\mathrm{d} \xi \mathrm{d} \bar{\xi},\\
\frac{\partial^2 F}{\partial \xi_k\partial \bar\xi_l}= &\lim_{\epsilon\to 0}\int_{|\xi-\xi_k|>\epsilon, |\xi-\xi_l|>\epsilon}\frac{\bar\xi-\bar\xi_k}{2\left|\xi-\xi_k\right|^3}\frac{\xi-\xi_l}{2\left|\xi-\xi_l\right|^3} \frac{1}{\left|\prod_{m\not =k,l}(\xi-\xi_m)\right|}~\mathrm{i}\,\mathrm{d} \xi \mathrm{d} \bar{\xi},\\
\frac{\partial^2 F}{\partial \bar\xi_k\partial \bar\xi_l}= &\lim_{\epsilon\to 0}\int_{|\xi-\xi_k|>\epsilon, |\xi-\xi_l|>\epsilon}\frac{\xi-\xi_k}{2\left|\xi-\xi_k\right|^3}\frac{\xi-\xi_l}{2\left|\xi-\xi_l\right|^3} \frac{1}{\left|\prod_{m\not =k,l}(\xi-\xi_m)\right|}~\mathrm{i}\,\mathrm{d} \xi \mathrm{d} \bar{\xi}.
\end{align*}
For any $k=1,2,3,4$,
\[
\frac{\partial^2 F}{\partial \xi_k\partial \xi_k}=-\sum_{l\not=k}\frac{\partial^2 F}{\partial \xi_k\partial \xi_l},\,
\frac{\partial^2 F}{\partial \xi_k\partial \bar\xi_k}=-\sum_{l\not=k}\frac{\partial^2 F}{\partial \xi_k\partial \bar\xi_l},\,
\frac{\partial^2 F}{\partial \bar\xi_k\partial \bar\xi_k}=-\sum_{l\not=k}\frac{\partial^2 F}{\partial \bar\xi_k\partial \bar\xi_l}\]
because we can shift $\xi$ by any variation of $\xi_k$ in $\frac{\partial F}{\partial \xi_k}$ or $\frac{\partial F}{\partial \bar{\xi}_k}$.

Now we use these to calculate
\begingroup
\allowdisplaybreaks
\begin{align*}
  &g_{\mathrm{sK}}(\dot{\nu},\dot{\nu}) =|t|^{1/2}\int_{\mathbb{C}}\frac{1}{\left|\prod_{k=1}^4(\xi-\xi_k)\right|}~\mathrm{i}\,\mathrm{d}\xi \mathrm{d}\bar{\xi} = |t|^{1/2}F(\xi_k, \bar\xi_k)\\
 =&|t|^{1/2}\left( F(\mathrm{e}^{\mathrm{i}(k-1)\pi/2}, \mathrm{e}^{-\mathrm{i}(k-1)\pi/2}) + \sum_{k=1}^{4}2\mathrm{Re}\left(\frac{\partial F}{\partial \xi_k}\frac{(-1)^k\mu_7t^{-1/4}}{4}\right)\right.\\
 &+\sum_{k=1}^{4}2\mathrm{Re}\left(\frac{\partial F}{\partial \xi_k}(t^{-1/2} \mathrm{e}^{-\mathrm{i}(k-1)\pi/2})(-\mu_6/4 - \mu_7^2/32)\right)\\
&+\mathrm{Re}\left(\sum_{k=1}^{4}\sum_{l=1}^{4}\frac{\partial^2 F}{\partial \xi_k\partial \xi_l}(-1)^{k+l}\frac{\mu_7^2 t^{-1/2}}{16}\right)\\& +\left. \mathrm{Re}\left(\sum_{k=1}^{4}\sum_{l=1}^{4}\frac{\partial^2 F}{\partial \xi_k\partial \bar\xi_l}\right)(-1)^{k+l}\frac{|\mu_7|^2 |t|^{-1/2}}{16}+O(|t|^{-3/4})\right).
\end{align*}
\endgroup
The terms involving $\frac{\partial F}{\partial \xi_k}$ cancel with each other. So in polar coordinates, as $|t|=r\to\infty$, the special Kähler metric is \[(C_0+\mathrm{Re} (C_1 \mu_7^2 t^{-1/2}) + C_2 |\mu_7|^2 |t|^{-1/2} + O(r^{-3/4}))r^{-1/2}g_{\mathrm{Euc}},\]
where $g_{\mathrm{Euc}}=\mathrm{d}r^2+r^2 \mathrm{d}\theta^2$. This metric is asymptotic to a conic metric with cone angle $3\pi/2$.

Next we consider $g_{\mathrm{sf}}$ along the fiber direction. By the discussion in Section \ref{Prelim_sec}, the Hitchin fiber $H^{-1}(t)$ can be written as
\begin{align}
  H^{-1}(t)&=\left\{\,\left[\left(\bar{\partial}_{E}^0,\varphi\right)\right]\,\middle|\,\varphi=\frac{\mathrm{d}z}{z^4}\begin{pmatrix}
    \pm (-t)^{1/2}z^2&1-\mu_7z-\mu_6z^2-\mu_5z^3\\1&\mp(-t)^{1/2}z^2
\end{pmatrix}\,\right\}\notag\\
&\cup\left\{\,\left[\left(\bar{\partial}_{E}^0,\varphi\right)\right]\,\middle|\,\varphi=\frac{\mathrm{d}z}{z^4}\begin{pmatrix}
    a_0&\begin{smallmatrix}
      -tz^3+(-\mu_5+tc_0)z^2+(-\mu_6+\mu_5c_0-tc_0^2)z  \\ {}-
      \mu_7+\mu_6c_0-\mu_5c_0^2+tc_0^3
    \end{smallmatrix}\\ c_0+z&-a_0
  \end{pmatrix},~a_0^2=-\tilde{\nu}(-c_0)\right\},\label{HitFibU4_eq}
\end{align}
where $(E,\bar{\partial}_E^0)\cong \mathcal{O}\oplus\mathcal{O}(-1)$. Therefore $H^{-1}(t)$ is the closure of the second set, which is the elliptic curve \[a_0^2=-\tilde{\nu}(-c_0)=-t\prod_{k=1}^4(-c_0-z_k(t)),~(a_0,c_0)\in \mathbb{C}^2,\]
with modulus $\tau(t)=\lambda^{-1}(l(t))$. Here $\lambda$ is the modular lambda function, and
\begin{align*}
  l(t)&=\frac{(z_4(t)-z_1(t))(z_3(t)-z_2(t))}{(z_3(t)-z_1(t))(z_4(t)-z_2(t))}=\frac{1}{2}+\frac{3\mathrm{i}\mu_7^2}{32}t^{-1/2}+O(|t|^{-3/4}), \\
  \tau(t)&=\mathrm{i}+\frac{3\pi \mu_7^2}{32K(\sqrt{1/2})^2}t^{-1/2}+O(|t|^{-3/4}),
\end{align*}
where $K(x)=\int_0^{\pi/2} (1-x^2\sin^2\theta)^{-1/2}\,\mathrm{d}\theta$ is the complete elliptic integral of the first kind. The semiflat metric induces a constant multiple of the Euclidean metric on $H^{-1}(t)$:
\begin{proposition}[{\cite[Proposition 8.4]{fredrickson2022asymptotic}}]
  The area of $H^{-1}(t)$ in the semiflat metric is $4\pi^2$, so $H^{-1}(t)\cong \mathbb{C}\,/\,c_t(\mathbb{Z}\oplus\tau(t)\mathbb{Z})$ with Euclidean metric $\mathrm{d}x^2+\mathrm{d}y^2$ and $c_t=2\pi/\sqrt{\mathrm{Im}\,\tau(t)}$.
\end{proposition}
Let $T_{\hat{\tau}}^2=\mathbb{C}\,/\,\hat{c}(\mathbb{Z}\oplus \hat{\tau}\mathbb{Z})$, where $\hat{\tau}=i, \hat{c}=2\pi/\sqrt{\mathrm{Im}\,\hat{\tau}}$. Define $\mu_t: T_{\hat{\tau}}^2\to H^{-1}(t)$ as \[x+\mathrm{i}y\mapsto \frac{c_t}{\hat{c}}\Big(x+\frac{\mathrm{Re}\,(\tau(t)-\hat{\tau})}{\mathrm{Im}\,\hat{\tau}}y+\mathrm{i}\frac{\mathrm{Im}\,\tau(t)}{\mathrm{Im}\,\hat{\tau}}y\Big)=\frac{c_t}{\hat{c}}(x+y\mathrm{Re}\,\tau(t)+\mathrm{i}y \mathrm{Im}\,\tau(t)).\]
Then we have
\begin{align*}
\mu_t^\ast (g_{\mathrm{sf}}|_{H^{-1}(t)})&=\frac{c_t^2}{\hat{c}^2}\left((\mathrm{d}x+\mathrm{Re}\,\tau(t)\,\mathrm{d}y)^2+(\mathrm{Im}\,\tau(t))^2\,\mathrm{d}y^2\right)
\\&=\mathrm{d}x^2+\mathrm{d}y^2+\mathrm{Im}\Big(\frac{3\pi \mu_7^2}{32K(\sqrt{1/2})^2}t^{-1/2}\Big)(-\mathrm{d}x^2+\mathrm{d}y^2)\\
&\quad+2\mathrm{Re}\Big(\frac{3\pi \mu_7^2}{32K(\sqrt{1/2})^2}t^{-1/2}\Big)\,\mathrm{d}x\mathrm{d}y+O(|t|^{-3/4}).
\end{align*}
Recall the $\beta$ and $\tau$ given in Table \ref{ALG_tab}. We conclude that the semiflat metric is asymptotic to an ALG-$A_1$ model metric. As shown in Table \ref{ALG_tab}, this case has 5 parameters, which come from $\mu_5,\mu_6,\mu_7,\mu_8$ and $\alpha_1$, after the reduction by Möbius transformations fixing the pole (2 complex parameters are reduced).

We summarize similar results for other ALG types in the following table.
\begin{table}
\begin{longtable}{| c | c |}
\hline
Type&$\uppercase\expandafter{\romannumeral2}^\ast~(A_0)$ \\
$D$& $4\cdot\{\tilde{0}\}$\\
$\mathcal{B}$&$\left(-z^{-7}+\mu_6z^{-6}+\mu_5z^{-5}+tz^{-4}\right)\,\mathrm{d}z^2$  \\
$g_{\mathrm{sK}}$&$(C_0+\mathrm{Re}(C_1 \mu_5 t^{-2/3})+\mathrm{Re}(C_2 \mu_6^2t^{-2/3})+C_3|\mu_6|^2r^{-2/3}+O(r^{-1}))r^{-1/3}g_{\mathrm{Euc}}$\\
$H^{-1}(t)$& $\begin{aligned}\bar{\partial}_E^0,~ \varphi&=\frac{\mathrm{d}z}{z^4}\begin{pmatrix}
\pm (-t)^{1/2}z^2&z-\mu_6z^2-\mu_5z^3\\1&\mp(-t)^{1/2}z^2
\end{pmatrix};\\
\bar{\partial}_E^0,~\varphi&=\frac{\mathrm{d}z}{z^4}\begin{pmatrix}
 a_0&\begin{smallmatrix}
   -tz^3+(-\mu_5+tc_0)z^2+(-\mu_6+\mu_5c_0-tc_0^2)z  \\{}+
   1+\mu_6c_0-\mu_5c_0^2+tc_0^3
 \end{smallmatrix}\\ c_0+z&-a_0
\end{pmatrix},~a_0^2=-\tilde{\nu}(-c_0).\\
(E,\bar{\partial}_E^0)&\cong \mathcal{O}\oplus\mathcal{O}(-1).
\end{aligned}$\\
$\tau(t)$&$\mathrm{e}^{2\pi \mathrm{i}/3}+\dfrac{\mathrm{e}^{-\pi \mathrm{i}/3}\pi (3\mu_5+\mu_6^2)}{12\sqrt{3}K(\mathrm{e}^{-\pi \mathrm{i}/6})^2}t^{-2/3}+O(|t|^{-4/3})$\\
\hline
Type&$\uppercase\expandafter{\romannumeral3}^\ast~(A_1)$\\
$D$& $3\cdot\{\tilde{0}\}+\{\infty\}$\\
$\mathcal{B}$&$\left(-w+\mu_4+tw^{-1}+\mu_2w^{-2}\right)\,\mathrm{d}w^2$  \\
$g_{\mathrm{sK}}$&$(C_0+\mathrm{Re}(C_1\mu_4t^{-1/2})+O(r^{-1}))r^{-1/2}g_{\mathrm{Euc}}$\\
$H^{-1}(t)$& $\begin{aligned}\bar{\partial}_E^1,~ \varphi&=\frac{\mathrm{d}w}{w}\begin{pmatrix}
0&w^3-\mu_4w^2-tw-\mu_2\\1&0
\end{pmatrix};\\
\bar{\partial}_E^0,~\varphi&=\frac{\mathrm{d}w}{w}\begin{pmatrix}
a_0&\begin{smallmatrix}
w^2-(\mu_4+c_0)w\\{}+\mu_4c_0+c_0^2-t
\end{smallmatrix}\\c_0+w & -a_0
\end{pmatrix},~a_0^2=-\tilde{\nu}(-c_0).\\
(E,\bar{\partial}_E^0)&\cong \mathcal{O}(-1)\oplus \mathcal{O}(-1),~ (E,\bar{\partial}_E^1)\cong \mathcal{O}\oplus \mathcal{O}(-2).
\end{aligned}$\\
$\tau(t)$&$\mathrm{i}+\dfrac{\mu_4\pi \mathrm{i} }{4K(\sqrt{1/2})^2}t^{-1/2}+O(|t|^{-1})$\\
\hline
Type&$\uppercase\expandafter{\romannumeral4}^\ast~(A_2)$\\
$D$& $3\cdot\{0\}+\{\infty\}$\\
$\mathcal{B}$&$\left(-w^2+\mu_5w+\mu_4+tw^{-1}+\mu_2w^{-2}\right)\,\mathrm{d}w^2$  \\
$g_{\mathrm{sK}}$&$(C_0+\mathrm{Re}(C_1t^{-1/3})+O(r^{-2/3}))r^{-2/3}g_{\mathrm{Euc}}$\\
$H^{-1}(t)$& $\begin{aligned}\bar{\partial}_E^1,~ \varphi&=\frac{\mathrm{d}w}{w}\begin{pmatrix}
0&w^4-\mu_5w^3-\mu_4w^2-tw-\mu_2\\1&0
\end{pmatrix};\\
\bar{\partial}_{E}^1,~\varphi&=\frac{\mathrm{d}w}{w}\begin{pmatrix}
-\mu_5w/2+w^2&-(\mu_5^2/4+\mu_4)w^2-tw-\mu_2\\1&\mu_5w/2-w^2
\end{pmatrix};\\
\bar{\partial}_{E}^0,~\varphi&=\frac{\mathrm{d}w}{w}\begin{pmatrix}
a_0+w^2&\begin{smallmatrix}
-\mu_5w^2+(c_0\mu_5-\mu_4-2a_0)w\\{}-c_0^2\mu_5+c_0\mu_4+2a_0c_0-t
\end{smallmatrix}\\c_0+w & -a_0-w^2
\end{pmatrix},~(a_0+c_0^2)^2=-\tilde{\nu}(-c_0).\\
(E,\bar{\partial}_E^0)&\cong \mathcal{O}(-1)\oplus \mathcal{O}(-1),~ (E,\bar{\partial}_E^1)\cong \mathcal{O}\oplus \mathcal{O}(-2).
\end{aligned}$\\
$\tau(t)$&$\mathrm{e}^{2\pi \mathrm{i}/3}+\dfrac{\sqrt{3}\pi}{4K(\mathrm{e}^{-\pi/6})^2}t^{-1/3}+O(|t|^{-2/3})$\\
\hline
\caption{other ALG cases}
\label{tab:other ALG}
\end{longtable}
\end{table}

\subsubsection{ALG$^\ast$ metrics} When there exists at least one pole of order $2$, the semiflat metric is asymptotic to an ALG$^\ast$ metric. We compute in detail for the case of two untwisted order two poles.

Suppose these poles are located at $0$ and $\infty$ respectively, then \[\mathcal{B}=\left\{\,\left(\frac{\mu_4}{z^4}+\frac{\mu_3}{z^3}+\frac{t}{z^2}+\frac{\mu_1}{z}+\mu_0\right)\,\mathrm{d}z^2\,\middle|\,t\in\mathbb{C}\,\right\},\]
where $\mu_0,\mu_1,\mu_3,\mu_4$ are fixed and $\mu_0,\mu_4\neq 0$, we may assume that $\mu_4=-1$. Fix \[\nu=\left(\frac{-1}{z^4}+\frac{\mu_3}{z^3}+\frac{t}{z^2}+\frac{\mu_1}{z}+\mu_0\right)\,\mathrm{d}z^2\in\mathcal{B}',~\dot{\nu}=\frac{\dot{t}}{z^2}\,\mathrm{d}z^2\in T_\nu\mathcal{B}'.\]
Then $g_{\mathrm{sK}}$ on $\mathcal{B}'$ is
\begin{align*}
  g_{\mathrm{sK}}(\dot{\nu},\dot{\nu})&=\int_{C}\frac{|\dot{\nu}|^2}{|\nu|}\,\mathrm{dvol}_C=|\dot{t}|^2\int_{\mathbb{C}}\frac{1}{\left|\mu_0\prod_{k=1}^4 (z-z_k(t))\right|}~\mathrm{i}\,\mathrm{d}z\mathrm{d}\bar{z}\\
  &=\frac{|\dot{t}|^2}{|\mu_0|}\int_{|z|\leqslant 1}\frac{1}{\left|\prod_{k=1}^4(z-z_k(t))\right|}\,\mathrm{i}\,\mathrm{d}z\mathrm{d}\bar{z}+|\dot{t}|^2\int_{|w|\leqslant 1}\frac{1}{\left|\prod_{k=1}^4(w-z_k(t)^{-1})\right|}~\mathrm{i}\,\mathrm{d}w\mathrm{d}\bar{w},
\end{align*}
where $z_k(t)=(-1)^{k-1}t^{-1/2}-\mu_3t^{-1}/2+O(|t|^{-3/2})$ ($k=1,2$), $z_k(t)=(-1)^{k-1}(-\mu_0)^{-1/2}t^{1/2}-\mu_1/2\mu_0+O(|t|^{-1/2})$ ($k=3,4$) are the four roots of $z^4\nu(z)=:\tilde{\nu}(z)$. In $\{|z|\leqslant 1\}$, $|z-z_k(t)|=|\mu_0^{-1/2}t^{1/2}|(1+O(|t|^{-1/2}))$ for $k=3,4$, then the first integral is
\begin{align*}
  &\int_{|z|\leqslant 1}\frac{|\mu_0t^{-1}|}{|z-z_1(t)||z-z_2(t)|}(1+O(|t|^{-1/2}))~\mathrm{i}\,\mathrm{d}z\mathrm{d}\bar{z}\quad (\text{let }\xi=t^{1/2}z)\\=&\int_{|\xi|\leqslant |t^{1/2}|}\frac{1}{|\xi-t^{1/2}z_1(t)||\xi-t^{1/2}z_2(t)|}~\mathrm{i}\,\mathrm{d}\xi \mathrm{d}\bar{\xi}\cdot|\mu_0t^{-1}|(1+O(|t|^{-1/2}))\\
  =&\int_{|\xi|\leqslant |t^{1/2}|}\frac{1}{|\xi-1||\xi+1|}~\mathrm{i}\,\mathrm{d}\xi \mathrm{d}\bar{\xi}\cdot|\mu_0t^{-1}|(1+O(|t|^{-1/2}))\\
  =&2\pi |\mu_0||t|^{-1}\log|t|+O(|t|^{-1}).
\end{align*}
Similarly, the second integral is $2\pi |t|^{-1}\log|t|+O(|t|^{-1})$.
Therefore in polar coordinates, as $|a|=r\to\infty$, the special Kähler metric is \[(4\pi\log\,r+O(1))r^{-1}g_{\mathrm{Euc}}.\]

The fiber $H^{-1}(t)$ can be written as
\begin{align*}
  H^{-1}(t)&=\left\{\,\left[\left(\bar{\partial}_{E}^1,\varphi\right)\right]\,\middle|\,\varphi=\frac{\mathrm{d}z}{z^2}\begin{pmatrix}
    0&1-\mu_3z-tz^2-\mu_1z^3-\mu_0z^4\\1&0
\end{pmatrix}\,\right\}\\
  &\cup\left\{\,\left[\left(\bar{\partial}_{E}^0,\varphi\right)\right]\,\middle|\,\varphi=\frac{\mathrm{d}z}{z^2}\begin{pmatrix}
  -\frac{\mu_1z}{2(-\mu_0)^{1/2}}+(-\mu_0)^{1/2}z^2&-(t+\mu_1^2/4\mu_0)z^2-\mu_3z+1\\1&  \frac{\mu_1z}{2(-\mu_0)^{1/2}}-(-\mu_0)^{1/2}z^2
\end{pmatrix}\,\right\}\\
  &\cup\left\{\,\left[\left(\bar{\partial}_{E}^0,\varphi\right)\right]\,\middle|\,\varphi=\frac{\mathrm{d}z}{z^2}\begin{pmatrix}
   a_0+(-\mu_0)^{1/2}z^2&\begin{smallmatrix}
     -\mu_1z^2+(\mu_1c_0-2a_0(-\mu_0)^{1/2}-t)z\\{}-
     \mu_1c_0^2+2a_0c_0(-\mu_0)^{1/2}+tc_0-\mu_3
   \end{smallmatrix}\\c_0+z&-a_0-(-\mu_0)^{1/2}z^2
 \end{pmatrix},\right.\notag\\&\left.\hspace{8cm}\vphantom{\begin{pmatrix}
     1+a_1z&\begin{smallmatrix}
       \mu_1c_0-2a_0(-\mu_0)^{1/2}\\a_0c_0(-\mu_0)^{1/2}
     \end{smallmatrix}\\1+c_1z&-(-\mu_0)^{1/2}z^2
  \end{pmatrix}}  (a_0+(-\mu_0)^{1/2}c_0^2)^2=-\tilde{\nu}(-c_0)\,\right\},
\end{align*}
where $(E,\bar{\partial}_E^0)\cong \mathcal{O}(-1)\oplus \mathcal{O}(-1)$, $(E,\bar{\partial}_E^1)\cong \mathcal{O}\oplus \mathcal{O}(-2)$. $H^{-1}(t)$ is the closure of the last set, which is the elliptic curve
\[(a_0+(-\mu_0)^{1/2}c_0^2)^2=-\tilde{\nu}(-c_0)=-\mu_0\prod_{k=1}^4(-c_0-z_k(t)),~(a_0,c_0)\in \mathbb{C}^2,\]
with modulus $\tau(t)=\lambda^{-1}(l(t))$, where
\begin{align*}
  l(t)&=\frac{(z_1(t)-z_2(t))(z_3(t)-z_4(t))}{(z_3(t)-z_2(t))(z_1(t)-z_4(t))}\\&=4(-\mu_0)^{1/2}t^{-1}+(8\mu_0+(-\mu_0)^{-1/2}(\mu_1^2-\mu_0\mu_3^2)/2)t^{-2}+O(|t|^{-3}),\\
  \tau(t)&=\frac{\mathrm{i}}{\pi}\log(4(-\mu_0)^{-1/2}t)+\frac{\mathrm{i}}{\pi}(\mu_0^{-1}(\mu_1^2-\mu_0\mu_3^2)/8)t^{-1}+O(|t|^{-2}).
\end{align*}
 $H^{-1}(t)\cong \mathbb{C}\,/\,c_t(\mathbb{Z}\oplus\tau(t)\mathbb{Z})$ with Euclidean metric and $c_t=2\pi/\sqrt{\mathrm{Im}\,\tau(t)}$.

 We summarize other ALG$^\ast$ cases below.
 \begin{longtable}{| c | c |}
   \hline
   Type&$\uppercase\expandafter{\romannumeral1}_4^\ast~(D_0)$ \\
   $D$& $2\cdot\{\tilde{0}\}+2\cdot\{\tilde{\infty}\}$\\
   $\mathcal{B}$&$\left(-z^{-3}+tz^{-2}+\mu_1z^{-1}\right)\,\mathrm{d}z^2$  \\
   $g_{\mathrm{sK}}$&$(8\pi\log\,r+O(1))r^{-1}g_{\mathrm{Euc}}$\\
   $H^{-1}(t)$& $\begin{aligned}\bar{\partial}_E^1,~ \varphi&=\frac{\mathrm{d}z}{z^2}\begin{pmatrix}
     0&z(1-tz+\mu_1z^2)\\1&0
  \end{pmatrix};\\
   \bar{\partial}_E^0,~\varphi&=\frac{\mathrm{d}z}{z^2}\begin{pmatrix}
    a_0&\begin{smallmatrix}
      -\mu_1z^2+(\mu_1c_0-a)z\\{}+1+tc_0-\mu_1c_0^2
    \end{smallmatrix} \\c_0+z&-a_0
  \end{pmatrix},~a_0^2=c_0\tilde{\nu}(-c_0).\\
   (E,\bar{\partial}_E^0)&\cong \mathcal{O}(-1)\oplus \mathcal{O}(-1),~ (E,\bar{\partial}_E^1)\cong \mathcal{O}\oplus \mathcal{O}(-2).
   \end{aligned}$\\
   $\tau(t)$&$(\mathrm{i}/\pi)(\log(16\mu_1^{-1}t^2)+5\mu_1t^{-2}/2)+O(|t|^{-4})$\\
 \hline
 Type&$\uppercase\expandafter{\romannumeral1}_3^\ast~(D_1)$ \\
 $D$& $2\cdot\{0\}+2\cdot\{\tilde{\infty}\}$\\
 $\mathcal{B}$&$\left(-z^{-4}+\mu_3z^{-3}+tz^{-2}+\mu_1z^{-1}\right)\,\mathrm{d}z^2$  \\
 $g_{\mathrm{sK}}$&$(6\pi\log\,r+O(1))r^{-1}g_{\mathrm{Euc}}$\\
 $H^{-1}(t)$& $\begin{aligned}\bar{\partial}_E^1,~ \varphi&=\frac{\mathrm{d}z}{z^2}\begin{pmatrix}
   0&1-\mu_3z-tz^2-\mu_1z^3\\1&0
\end{pmatrix};\\
 \bar{\partial}_E^0,~\varphi&=\frac{\mathrm{d}z}{z^2}\begin{pmatrix}
  a_0&\begin{smallmatrix}
    -\mu_1z^2+(\mu_1c_0-t)z\\{}-
    \mu_3+tc_0-\mu_1c_0^2
  \end{smallmatrix}\\c_0+z&-a_0
\end{pmatrix},~a_0^2=-\tilde{\nu}(-c_0).\\
 (E,\bar{\partial}_E^0)&\cong \mathcal{O}(-1)\oplus \mathcal{O}(-1),~ (E,\bar{\partial}_E^1)\cong \mathcal{O}\oplus \mathcal{O}(-2).
 \end{aligned}$\\
 $\tau(t)$&$(\mathrm{i}/\pi)\left(\log(8\mu_1^{-1}t^{3/2})-\mu_3^2t^{-1}/8-\mu_1t^{-3/2}\right)+O(|t|^{-2})$\\
 \hline
 Type&$\uppercase\expandafter{\romannumeral1}_2^\ast~(D_2)$ \\
 $D$& $\{0\}+\{1\}+2\cdot\{\tilde{\infty}\}$\\
 $\mathcal{B}$&$\left(\mu_0z^{-2}+\mu_1(z-1)^{-2}+\mu_3z^{-1}+t(z(z-1))^{-1}\right)\,\mathrm{d}z^2$  \\
 $g_{\mathrm{sK}}$&$(4\pi\log\,r+O(1))r^{-1}g_{\mathrm{Euc}}$\\
 $H^{-1}(t)$& $\begin{aligned}\bar{\partial}_E^0,~ \varphi&=\frac{\mathrm{d}z}{z(z-1)}\begin{pmatrix}
   0&-\tilde{\nu}(z)\\1&0
\end{pmatrix};\\
 \bar{\partial}_E^0,~\varphi&=\varphi=\frac{\mathrm{d}z}{z(z-1)}\begin{pmatrix}
  a_0&\begin{smallmatrix}
   -\mu_3z^2+(\mu_3c_0+2\mu_3-\mu_0-\mu_1-t)z\\{}+
   t+2\mu_0-\mu_3\\{}+(t+\mu_0+\mu_1-2\mu_3)c_0-\mu_3c_0^2
  \end{smallmatrix}\\c_0+z&-a_0
 \end{pmatrix}, a_0^2=-\tilde{\nu}(-c_0).\\
 (E,\bar{\partial}_E^0)&\cong \mathcal{O}(-1)\oplus \mathcal{O}(-2).
 \end{aligned}$\\
 $\tau(t)$&$(\mathrm{i}/\pi)(\log(-16\mu_3^{-1}t)+(\mu_3/2-1)t^{-1})+O(|t|^{-2})$\\
 \hline
 Type&$\uppercase\expandafter{\romannumeral1}_1^\ast~(D_3=A_3)$ \\
 $D$& $\{0\}+\{1\}+2\cdot\{\infty\}$\\
 $\mathcal{B}$&$\left(\mu_0z^{-2}+\mu_1(z-1)^{-2}+\mu_3z^{-1}+\mu_4+t(z(z-1))^{-1}\right)\,\mathrm{d}z^2$  \\
 $g_{\mathrm{sK}}$&$(2\pi\log\,r+O(1))r^{-1}g_{\mathrm{Euc}}$\\
 $H^{-1}(t)$& $\begin{aligned}\bar{\partial}_E^0,~ \varphi&=\frac{\mathrm{d}z}{z(z-1)}\begin{pmatrix}
   \pm (-\mu_4)^{1/2}z^2&\mu_4z^4-\tilde{\nu}(z)\\1&\mp(-\mu_4)^{1/2}z^2
\end{pmatrix};\\
 \bar{\partial}_E^0,~\varphi&=\varphi=\frac{\mathrm{d}z}{z(z-1)}\begin{pmatrix}
  a_0&\begin{smallmatrix}
    -\mu_4z^3+(\mu_4c_0+2\mu_4-\mu_3)z^2\\{}+
    (2\mu_3-\mu_4-\mu_1-t-\mu_4c_0^2-2\mu_4c_0+\mu_3c_0)z\\{}+2\mu_0+t-
   \mu_3+(\mu_4-2\mu_3+\mu_1+t)c_0\\{}+(2\mu_4-\mu_3)c_0^2+\mu_4c_0^3
  \end{smallmatrix}\\c_0+z&-a_0
 \end{pmatrix},~a_0^2=-\tilde{\nu}(-c_0).\\
 (E,\bar{\partial}_E^0)&\cong \mathcal{O}(-1)\oplus \mathcal{O}(-2).
 \end{aligned}$\\
 $\tau(t)$&$(\mathrm{i}/\pi)(\log(-8(-\mu_4)^{-1/2}t^{1/2})+O(|t|^{-1})$\\
 \hline
 \caption{other ALG$^\ast$ cases}
 \label{tab:other ALGst}
 \end{longtable}

\subsection{The uniformness}
In Theorem \ref{Main_thm}, the constants in the decay rate of $g_{L^2}-g_{\mathrm{sf}}$ depend on the choice of the curve. Now we aim to find constants uniform in these choices for a four-dimensional moduli space $\mathcal{M}$.

We first consider the moduli space of Higgs bundles with an untwisted order four pole. As in \eqref{HitBaseU4_eq}, $t\in \mathbb{C}$ parametrizes the quadratic differential $\nu_t$ in the Hitchin base, and is no longer a positive real number as in the previous sections. Let $(\bar{\partial}_E,\varphi_t)$ be any Higgs bundle in the Hitchin fiber $H^{-1}(t)$. By Lemma \ref{HiggsDetRoot_lem}, the four roots of $\tilde{\nu}_t(z)$ have asymptotics
\[z_k(t)=t^{-1/4}\mathrm{e}^{\mathrm{i}(k-1)\pi/2}+O(|t|^{-1/2}),\]
as $|t|\to\infty$. The local mass at $z_k(t)$ $(k=1,2,3,4)$ is
\[\lambda_k(t)=\Big|t\prod_{j\neq k}(z_k(t)-z_j(t))z_k(t)^{-8}\Big|^{1/2}=4t^{9/4}(1+O(|t|^{-1/4})).\]
Here and below, the constants are independent of $\mathrm{Arg}(t)$ and the choice of $(\bar{\partial}_E,\varphi_t)$ in $H^{-1}(t)$. Proposition \ref{NormalForm_prop} hold on disks of uniform size. Part (\romannumeral5) is easy to check. For part (\romannumeral1), we can find $\zeta_{k,t}$ such that $\det\varphi_t=-\lambda_k(t)^2\zeta_{k,t}\,\mathrm{d}\zeta_{k,t}^2$ on $\widetilde{B}_{k,t}:=\{\,|\zeta_{k,t}|<\kappa |t|^{-1/4}\,\}$, we can also make $\kappa$ independent of $\mathrm{Arg}(t)$. Next, we need to check that the gauge transformations in the proof is well-defined on the whole disk $\widetilde{B}_{k,t}$ for some uniform $\kappa$. This is clear for $\varphi_t$ in the first subset of \eqref{HitFibU4_eq}. For $\varphi_t$ in the second subset, previously we fixed $a_0,c_0$, which now can depend on $t$. Suppose $|c_0+z_k(t)|\geqslant |t|^{-1/4}/2$, the first gauge transformation is $\left(\begin{smallmatrix}
  a_0&1\\ z_k(t)+c_0&0
\end{smallmatrix}\right)$, which is constant and invertible over $\widetilde{B}_{k,t}$. The second gauge transformation is well-defined since \[\left|\frac{z+c_0}{z_k(t)+c_0}-1\right|= \left|\frac{z-z_k(t)}{z_k(t)}\right|\left|\frac{z_k(t)}{c_0+z_k(t)}\right|\leqslant c\kappa\leqslant 1/2,\]
for $\kappa$ small. The rest of the proof remains the same. If $|c_0+z_k(t)|\leqslant |t|^{-1/4}/2$ for some $k$, then $|c_0+z_j(t)|\geqslant |t|^{-1/4}/2$ for $j\neq k$ and $|t|$ large, and the previous discussion applies for the disks $\widetilde{B}_{j,t}$. On $\widetilde{B}_{k,t}$,
\[\varphi_t=\frac{\mathrm{d}z}{z^4}\begin{pmatrix}
    a_0&\begin{smallmatrix}
      -tz^3+(-\mu_5+tc_0)z^2+(-\mu_6+\mu_5c_0-tc_0^2)z  \\-{}
      \mu_7+\mu_6c_0-\mu_5c_0^2+tc_0^3
    \end{smallmatrix}\\ c_0+z&-a_0
  \end{pmatrix}:=\begin{pmatrix}
    a_t(\zeta_{k,t})&b_t(\zeta_{k,t})\\c_t(\zeta_{k,t})&-a_t(\zeta_{k,t})
  \end{pmatrix}\,\mathrm{d}\zeta_{k,t},\]
  where \begin{align*}
    b_t(\zeta_{k,t})&=\frac{\mathrm{d}z}{\mathrm{d}\zeta_{k,t}}(\zeta_{k,t})\frac{-t}{z^4}(z^3-c_0z^2+c_0^2z-c_0^3)(1+O(|t|^{-1/4}))\\
    &=\frac{\mathrm{d}z}{\mathrm{d}\zeta_{k,t}}(\zeta_{k,t})\frac{-t}{z^4}(z-c_0)(z^2+c_0^2)(1+O(|t|^{-1/4})).
\end{align*}
Then $\sqrt{b_t(\zeta_{k,t})}$ is well-defined on $\widetilde{B}_{k,t}$ by the assumption that $|c_0+z_k(t)|\leqslant |t|^{-1/4}/2$, and we have \[g^{-1}\varphi_t g=\begin{pmatrix}
  0&1\\\lambda_k(t)^2\zeta_{k,t}&0
\end{pmatrix}\,\mathrm{d}\zeta_{k,t}, \text{ for }g=\frac{1}{\sqrt{b_t(\zeta_{k,t})}}\begin{pmatrix}
  b_t(\zeta_{k,t})&0\\-a_t(\zeta_{k,t})&1
\end{pmatrix}.\]
The remaining proof is the same. The approximate metric $h_t^{\mathrm{app}}$ is constructed as in Definition \ref{ApproxMet_def}, and the proof of the error estimate in \eqref{ErrorEst_eq} works uniformly. Then there are uniform constants $t_0,c,c'$, such that for any $|t|\geqslant t_0$,
\[
\big\lVert F_{h_t^\mathrm{app}}+\big[\varphi_t,\varphi_t^{\ast_{h_t^\mathrm{app}}}\big] \big\rVert_{L^2}\leqslant c\mathrm{e}^{-c' |t|^{3/4}}.
\]
The computations in Sections 4, 5, 6 only use the expression of $q_t=\det\,\varphi_t$ and the local forms of $\varphi_t$ provided by Proposition \ref{NormalForm_prop}, so the estimates remain valid uniformly for the choice of $\varphi_t$ in $H^{-1}(t)$. Moreover, the constants can also be chosen to be independent of $\mathrm{Arg}(t)$. Therefore, we find uniform constants $t_0,c,c'$ so that for $[(\dot{\eta},\dot{\varphi})]\in T_{[(\bar{\partial}_E,\varphi_t)]}\mathcal{M}$ and $|t|\geqslant t_0$,
  \[|\lVert [(\dot{\eta},\dot{\varphi})]\rVert_{g_{L^2}}/\lVert [(\dot{\eta},\dot{\varphi})]\rVert_{g_{\mathrm{sf}}}-1|\leqslant c\mathrm{e}^{-c'|t|^{3/4}}.\]

Next we consider the case when $D=3\cdot \{\tilde{0}\}+\{\infty\}$. As before, Proposition \ref{NormalForm_prop} (\romannumeral1) holds uniformly, and we can verify (\romannumeral4) similarly. When the pole at $\infty$ is weakly parabolic, by the proof of Proposition \ref{NormalForm_prop} (\romannumeral2), we can find a holomorphic coordinate $\zeta_{t}$ around $\infty$ such that $\det\varphi_t=\zeta_{t}^{-2}\mu_{x,2}(1-\epsilon_{x,0,t}^{-1}\zeta_{t})$ on $\widetilde{B}_{t}:=\{\,|\zeta_{t}|<\kappa t^{-2/3}\,\}$. By a rotation, we may assume $\kappa$ is independent of $\mathrm{Arg}(t)$. From Table \ref{ALG_tab}, the first subset is the Hitchin section, where the Higgs field is already off-diagonal. After the coordinate change $w\mapsto\zeta_t$, by a diagonal gauge transformation defined on $\widetilde{B}_t$ we can make the $(2,1)$-entry to be $1$. For the second subset of $H^{-1}(t)$, on $\widetilde{B}_{t}$,
\[\varphi_t=\frac{\mathrm{d}w}{w}\begin{pmatrix}
  a_0& w^2-(\mu_4+c_0)w+\mu_4 c_0+c_0^2-t\\c_0+w&-a_0
\end{pmatrix}:= \frac{\mathrm{d}\zeta_t}{\zeta_{t}}\begin{pmatrix}
  a_t(\zeta_t)&b_t(\zeta_t)\\ c_t(\zeta_t)&-a_t(\zeta_t)
\end{pmatrix}.\]
Again we consider two cases $|c_0|\geqslant \kappa |t|^{-2/3}/2$ and $|c_0|< \kappa |t|^{-2/3}/2$. Note that $\sqrt{b_t(\zeta_t)}$ is well-defined on $\widetilde{B}_t$ for the second case, then the previous arguments can be applied. When the pole at $\infty$ is strongly parabolic, we can also adjust the proof of Proposition \ref{NormalForm_prop} (\romannumeral3) uniformly.

The remaining cases can be analyzed in the same way.

\bibliographystyle{amsalpha}
\bibliography{Irregular-Higgs}
\end{document}